\documentclass[10pt,reqno]{amsart}
\usepackage[numbers]{natbib}
\usepackage{hyperref}

\usepackage[margin=1.0in]{geometry}
\usepackage{multicol}            
\usepackage{adjustbox}           

\numberwithin{equation}{section} 

\usepackage[utf8]{inputenc}
\usepackage[T1]{fontenc}
\usepackage{lmodern}     
\usepackage{libertinus}  
\usepackage{microtype}   
\usepackage{xspace}      

\usepackage{amsmath}
\usepackage{amssymb}
\usepackage{amsfonts}
\usepackage{amsthm}
\usepackage{faktor}
\usepackage{mathtools}   
\usepackage{mathrsfs}    
\usepackage{newtxtext}
\usepackage{newtxmath}   
\usepackage{mathdots}    
\usepackage{stackrel}    
\usepackage{bbm}         
\usepackage{latexsym}    
\usepackage{cases}
\usepackage{csquotes}
\usepackage{etoolbox}
\usepackage{mathrsfs}
\usepackage{anyfontsize}
\usepackage{enumitem}    
\usepackage{graphicx}    
\usepackage{array}       
\usepackage{comment}     
\usepackage{xcolor}      

\definecolor{cherryred}{HTML}{D2042D}

\usepackage[nameinlink,capitalize]{cleveref}

\usepackage{tikz}
\usepackage{tikz-3dplot}
\usepackage{quiver}      
\usepackage{tikz-cd}

\usetikzlibrary{cd}
\usetikzlibrary{nfold}
\usetikzlibrary{math}
\usetikzlibrary{fit}
\usetikzlibrary{calc}
\usetikzlibrary{arrows}
\usetikzlibrary{patterns}
\usetikzlibrary{positioning}
\usetikzlibrary{intersections}
\usetikzlibrary{matrix}
\usetikzlibrary{decorations.pathmorphing}  
\usetikzlibrary{decorations.pathreplacing} 
\usetikzlibrary{decorations.markings}      

\tikzset{
    Rightarrow/.style={double equal sign distance,>={Implies},->},
    triple/.style={-,preaction={draw,Rightarrow}},
    quadruple/.style={preaction={draw,Rightarrow,shorten >=0pt},shorten >=1pt,-,double,double distance=0.2pt}
}

\tikzset{
    symbol/.style={
        draw=none,
        every to/.append style={
            edge node={node [sloped, allow upside down, auto=false]{$#1$}}}
    }
}

\newcommand{\adjunction}[4]{%
  \begin{tikzcd}[
    ampersand replacement=\&,
    baseline=(current bounding box.center),
    column sep=5em
  ]
  #1
    \arrow[r, shift left=1.2ex, "#3"]
    \arrow[r, phantom, "\bot" description]
  \&
  #2
    \arrow[l, shift left=1.2ex, "#4"]
  \end{tikzcd}%
}

\tikzset{
    scalearrow/.style n args={3}{
        decoration={
            markings,
            mark=at position (1-#1)/2*\pgfdecoratedpathlength with {\coordinate (#2);},
            mark=at position (1+#1)/2*\pgfdecoratedpathlength with {\coordinate (#3);},
        },
        postaction=decorate,
    }
}

\makeatletter
\newcommand{\@bbify}[1]{
  \ifcsname b#1\endcsname
    \message{WARNING: Skipping b#1 to prevent infinite loop conflict!}
  \else
    \expandafter\edef\csname b#1\endcsname{\noexpand\ensuremath{\noexpand\mathbb #1}\noexpand\xspace}
  \fi}

\newcommand{\@calify}[1]{
  \ifcsname c#1\endcsname
    \message{WARNING: Skipping c#1 to prevent infinite loop conflict!}
  \else
    \expandafter\edef\csname c#1\endcsname{\noexpand\ensuremath{\noexpand\mathcal #1}\noexpand\xspace}
  \fi}

\newcommand{\@bfify}[1]{
  \ifcsname bf#1\endcsname
    \message{WARNING: Skipping bf#1 to prevent infinite loop conflict!}
  \else
    \expandafter\edef\csname bf#1\endcsname{\noexpand\ensuremath{\noexpand\mathbf #1}\noexpand\xspace}
  \fi}

\newcounter{@letter}\stepcounter{@letter}
\loop\@bbify{\Alph{@letter}}\@calify{\Alph{@letter}}\@bfify{\Alph{@letter}}
\ifnum\the@letter<26\stepcounter{@letter}\repeat
\makeatother
\theoremstyle{plain}
\newtheorem{theorem}{Theorem}[section]
\newtheorem{proposition}[theorem]{Proposition}
\newtheorem{lemma}[theorem]{Lemma}
\newtheorem{corollary}[theorem]{Corollary}

\newtheorem*{unnumberedtheorem}{Theorem}

\theoremstyle{definition}
\newtheorem{definition}[theorem]{Definition}
\newtheorem{construction}[theorem]{Construction}
\newtheorem{notation}[theorem]{Notation}

\newtheorem{convention}[theorem]{Convention}

\newtheorem{rmk}[theorem]{Remark}
\newtheorem{ex}[theorem]{Example}

\newtheorem*{unnumberedtheoremdefinition}{Definition}

\theoremstyle{remark}

\newtheorem{hypothesis}[theorem]{Hypothesis}

\makeatletter
\let\c@equation\c@theorem
\numberwithin{equation}{section}
\makeatother

\crefname{lem}{Lemma}{Lemmas}
\crefname{thm}{Theorem}{Theorems}
\crefname{defn}{Definition}{Definitions}
\crefname{notn}{Notation}{Notations}
\crefname{note}{Note}{Note}
\crefname{const}{Construction}{Constructions}
\crefname{prop}{Proposition}{Propositions}
\crefname{rmk}{Remark}{Remarks}
\crefname{goal}{Goal}{Goals}
\crefname{cor}{Corollary}{Corollaries}
\crefname{equation}{Display}{Displays}
\crefname{ex}{Example}{Examples}
\crefname{nonex}{Non-Example}{Non-Examples}
\crefname{thmalph}{Theorem}{Theorems}
\crefname{answ}{Answer}{Answers}
\crefname{question}{Question}{Questions}

\newcommand{\sset}{\mathbf{sSet}}

\newcommand{\reedyslice}{\bigl((\sset_{(\infty,1)})^{\Delta^{\mathrm{op}}}_{/\bisimFin}\bigr)}
\newcommand{\genreedyslice}{\mathbf M^{\Delta^{\mathrm{op}}}_{/\mathbb B}}
\newcommand{\cocart}{\bigl(\sset^+_{/(\simFin)^\sharp}\bigr)_{\mathrm{coCart}}}
\newcommand{\cocartsp}{\bigl(\sset^+_{/(\simFin)^\sharp}\bigr)_{\mathrm{coCart}}^{\mathrm{sp}}}
\newcommand{\genldblfib}{\bigl((\Minfn)^{\Delta^{\mathrm{op}}}_{/\bisimFin}\bigr)_{\mathrm{LDbl}}}
\newcommand{\genldblfibsp}{\bigl((\Minfn)^{\Delta^{\mathrm{op}}}_{/\bisimFin}\bigr)_{\mathrm{LDbl}}^{\mathrm{sp}}}
\newcommand{\ldblfib}{\bigl((\sset_{(\infty,1)})^{\Delta^{\mathrm{op}}}_{/\bisimFin}\bigr)_{\mathrm{LDbl}}}
\newcommand{\ldblfibsp}{\bigl((\sset_{(\infty,1)})^{\Delta^{\mathrm{op}}}_{/\bisimFin}\bigr)_{\mathrm{LDbl}}^{\mathrm{sp}}}

\newcommand{\Set}{\mathbf{Set}}

\newcommand{\Fin}{\mathrm{Fin}_*}

\newcommand{\RMap}{\mathbf{R}\mathrm{Map}}

\newcommand{\simFin}{\mathcal{F}\mathrm{in}_*}
\newcommand{\bisimFin}{\mathbb{F}\mathrm{in}_*}

\DeclareMathOperator{\Ho}{Ho}

\DeclareMathOperator{\Map}{Map}

\DeclareMathOperator{\Cocart}{Cocart}

\DeclareMathOperator{\id}{id}

\DeclareMathOperator{\pr}{pr}

\newcommand{\Minfn}{\mathbf M_{(\infty,n)}}
\newcommand{\HGen}{\mathsf G}
\newcommand{\Stair}{\mathcal{S}\mathrm{tair}}
\newcommand{\Int}{\mathbb{I}\mathrm{nt}}
\newcommand{\Diag}{\mathcal{D}\mathrm{iag}}
\newcommand{\diagset}{\mathrm{diag}}
\crefname{construction}{Construction}{Constructions}
\crefname{notation}{Notation}{Notations}
\crefname{convention}{Convention}{Conventions}
\crefname{definition}{Definition}{Definitions}
\crefname{remark}{Remark}{Remarks}
\crefname{hypothesis}{Hypothesis}{Hypotheses}

\patchcmd{\section}{\normalfont}{\normalfont\large}{}{}
\patchcmd{\section}{\scshape}{}{}{}

\makeatletter
\patchcmd{\@settitle}{\uppercasenonmath\@title}{}{}{}
\patchcmd{\@settitle}{\bfseries}{\Large\bfseries}{}{}
\makeatother

\title{An internal approach to symmetric monoidal $(\infty,n)$-categories}
\author{Arghan Dutta}
\address{Department of Mathematics and Statistics, University of Ottawa, Ottawa, Canada}
\email{adutt095@uottawa.ca}
\begin{document}
\flushbottom

\begin{abstract}
We propose a definition of symmetric monoidal $(\infty,n)$-categories
as special left double $(\infty,n)$-fibrations over the category of
finite pointed sets. The definition makes sense in a general model of
$(\infty,n)$-categories, under a mild condition on its discrete
objects, and when the model is left proper and combinatorial we
construct model structures whose fibrant objects are exactly these
fibrations. For $n=1$ we show that the definition recovers Lurie's
symmetric monoidal $(\infty,1)$-categories: the homotopy theory of
special left double $(\infty,1)$-fibrations is equivalent to that of
special coCartesian fibrations over the category of finite pointed
sets. The equivalence is realized by an explicit Quillen equivalence,
and two constructions, the staircase and coCartesian strings, translate
between the two descriptions.
\end{abstract}

\maketitle

\tableofcontents

\section*{Introduction}

\subsection*{Symmetric monoidal structures on higher categories}
A symmetric monoidal structure on a category is a tensor product of
objects and of morphisms, with a unit, that is associative, unital and
commutative up to coherent isomorphism. Vector spaces under the tensor
product and sets under the cartesian product are standard examples. In
a higher category, the tensor product has to act on morphisms of every
dimension, and the coherence data never stop: for instance,
associativity holds up to an equivalence, that equivalence satisfies
the pentagon up to a higher one, and so on.

Higher categories of bordisms provide one of the motivating examples.
Gluing bordisms is composition, disjoint union is the tensor product,
and both are already built into Atiyah's axioms for a topological
quantum field theory \cite{Atiyah1988}. Once manifolds of every
codimension are allowed, together with bordisms between bordisms, one
is led to higher categories \cite{BaezDolan1995}. A fully extended
framed field theory of dimension $n$ is then a symmetric monoidal
functor out of a bordism $(\infty,n)$-category, and the cobordism
hypothesis identifies the space of such theories with the space of
fully dualizable objects of the target \cite[\S2.4]{lurie2009tft}.
Complete Segal models of bordism categories, with their symmetric
monoidal structures, are constructed in \cite{calaque2019cobordisms};
locality and geometric structures are treated in
\cite{gradypavlov2020locality,GradyPavlov2026Bordisms}. There are
algebraic examples as well. Higher Morita categories of algebras and
bimodules \cite{Vitale1992,dutta2025morita} inherit a symmetric
monoidal structure from the category they are built in
\cite{HaugsengMorita,johnsonfreyd2017oplax}, and iterated spans in a
category with finite limits are symmetric monoidal under the cartesian
product \cite{Haugseng2018Spans}.

\subsection*{Symmetric monoidal $(\infty,n)$-categories in the literature}
In the bicategorical range, the tensor product, the unit and their
coherence data are written down explicitly
\cite{SchommerPries2011,GurskiOsorno2013}, and double categories have
been used to construct symmetric monoidal bicategories
\cite{Shulman2010}. In general, there are several ways to define a
symmetric monoidal $(\infty,n)$-category. One can ask for a
commutative algebra object in the $(\infty,1)$-category of
$(\infty,n)$-categories with its cartesian product
\cite[\S4.1]{lurie2009tft}, which is close in spirit to the iterated
enrichment of \cite{GepnerHaugseng2015}; for $n=1$, there are also
operadic and diagrammatic point-set models
\cite{KodjabachevSagave2015,NikolausSagave2017}. One can ask for a
compatible tower of deloopings
\cite{schommerpries2014dualizability,calaque2019cobordisms,Haugseng2018Spans}.

Or one can follow Segal \cite{Seg74} and consider Segal copresheaves
indexed by the category $\Fin$ of finite pointed sets
$\langle k\rangle=\{*,1,\ldots,k\}$, i.e., functors $X$ from
$\Fin$ to $(\infty,n)$-categories such that, for every $k\geq0$,
the $k$ \emph{inert} maps $\langle k\rangle\to\langle1\rangle$,
each of which keeps a single non-basepoint, induce an equivalence of
$(\infty,n)$-categories
\[
X\langle k\rangle\xrightarrow{\ \simeq\ }
\bigl(X\langle1\rangle\bigr)^{k}.
\]
Conditions of this form are studied systematically in
\cite{ChuHaugseng2021}. In specific models of
$(\infty,n)$-categories, symmetric monoidal $(\infty,n)$-categories
are defined as Segal copresheaves on $\Fin$ with values in Segal
$n$-categories in \cite{schommerpries2014dualizability}, and with
values in complete $n$-fold Segal spaces in
\cite{calaque2019cobordisms,johnsonfreyd2017oplax}. Model structures
obtained by localizing at the categorical and the symmetric monoidal
Segal conditions are constructed in \cite{gradypavlov2020locality},
building on \cite{Rezk2001,Rezk2010} and on the homotopy theory of
$\Gamma$-spaces of Bousfield and Friedlander
\cite{BousfieldFriedlander1978}; the models are compared in
\cite{GradyPavlov2026Bordisms,BarwickSchommerPries2021}.

Copresheaves on $\Fin$ correspond to coCartesian fibrations over
$\Fin$ through an appropriate straightening--unstraightening
correspondence \cite[\S3.2]{lurie2009higher}. Hence, one can also
define a symmetric monoidal higher category as a coCartesian fibration
$p:\mathcal E\to\simFin$ over the nerve $\simFin$ of $\Fin$ that
satisfies the corresponding version of the Segal condition; we call
such a coCartesian fibration \emph{special}. For $n=1$, this is
Lurie's definition \cite[Definition~2.0.0.7]{LurieHA}, and Sharma
compares diagrams of quasi-categories on finite pointed sets with a
marked fibrational presentation \cite{Sharma2021}. A
higher-dimensional version of this approach is pursued in
\cite{Nuiten2023}; see also \cite{GagnaHarpazLanari2024} for
fibrations of $(\infty,2)$-categories.

\subsection*{Symmetric monoidal $(\infty,n)$-categories as left double fibrations}
The goal of this paper is to provide one way of packaging all the
coherence data of a symmetric monoidal structure on an
$(\infty,n)$-category, in a form that is not tied to a particular
model of $(\infty,n)$-categories. For $n=1$, we show that it
recovers exactly Lurie's established framework of symmetric monoidal
$(\infty,1)$-categories. To this end, we take the fibrational point
of view in every dimension, but we replace coCartesian fibrations of
$(\infty,n)$-categories $p:\mathcal E\to\simFin$ by left fibrations
of double $(\infty,n)$-categories $P:\mathbb E\to\bisimFin$, which
make sense in a general model of $(\infty,n)$-categories.

Recall that a \emph{double $(\infty,n)$-category} is a simplicial
object $\mathbb E$ in $(\infty,n)$-categories that satisfies the
Segal condition \cite{Moser2026InftyN}, i.e., such that, for every
$m\geq2$, the Segal map
\[
\mathbb E_m\longrightarrow
\underbrace{\mathbb E_1\times_{\mathbb E_0}\cdots
\times_{\mathbb E_0}\mathbb E_1}_{m\text{ factors}}
\]
is an equivalence of $(\infty,n)$-categories. We refer to
$\mathbb E_0$ as the $(\infty,n)$-category of objects and
\emph{vertical} morphisms, and to $\mathbb E_1$ as the
$(\infty,n)$-category of \emph{horizontal} morphisms and squares. For
example, the nerve of $\Fin$, with each of its levels regarded as a
discrete $(\infty,n)$-category, is a double $(\infty,n)$-category
$\bisimFin$ whose horizontal morphisms are the maps of $\Fin$ and
whose vertical morphisms are all trivial.

A \emph{left double $(\infty,n)$-fibration} is a map of double
$(\infty,n)$-categories $P:\mathbb E\to\bisimFin$, subject to a
fibrancy condition, such that every horizontal morphism of
$\mathbb E$ is essentially determined by its source and by its image
in $\Fin$, i.e., the canonical map
\[
\mathbb E_1\longrightarrow
\mathbb E_0\times_{(\bisimFin)_0}(\bisimFin)_1
\]
is an equivalence of $(\infty,n)$-categories; see
\cref{def:left-double-infinity-n-fibration}. Like coCartesian
fibrations, left double fibrations over $\bisimFin$ encode
copresheaves on $\Fin$ with values in $(\infty,n)$-categories,
thanks to the straightening--unstraightening construction of
\cite{moser2023inftyncategoricalstraighteningunstraighteningconstruction};
the value of the copresheaf at $\langle k\rangle$ is the fiber
$\mathbb E_0(\langle k\rangle)$ of $\mathbb E_0$ over
$\langle k\rangle$. We are then led to introduce as
\cref{def:special_left_double} the following notion.

\begin{unnumberedtheoremdefinition}
A left double $(\infty,n)$-fibration $P:\mathbb E\to\bisimFin$ is
\emph{special} if, for every $k\geq0$, the comparison between the
fiber $\mathbb E_0(\langle k\rangle)$ and the $k$-fold product of
the fiber $\mathbb E_0(\langle1\rangle)$ induced by the inert maps
$\langle k\rangle\to\langle1\rangle$ is an equivalence of
$(\infty,n)$-categories. A \emph{symmetric monoidal
$(\infty,n)$-category} is a special left double
$(\infty,n)$-fibration.
\end{unnumberedtheoremdefinition}

\subsection*{Overview of the results}
In this paper, we show as a first consistency check that the
definition is sound for $n=1$. Precisely, we show as
\cref{thm:special_quillen_equivalence} that:

\begin{unnumberedtheorem}
There is an equivalence of homotopy theories between symmetric monoidal
$(\infty,1)$-categories presented by special coCartesian fibrations
$p:\mathcal E\to\simFin$ and symmetric monoidal
$(\infty,1)$-categories presented by special left double fibrations
$P:\mathbb E\to\bisimFin$.
\end{unnumberedtheorem}

We prove this result through a Quillen equivalence between the two
appropriate model structures. On fibrant objects, the equivalence is
realized by the following two constructions.
\[
\begin{tikzcd}[ampersand replacement=\&,column sep=8.5em]
\left\{\begin{tabular}{@{}c@{}}special left double fibrations\\
$P:\mathbb E\to\bisimFin$\end{tabular}\right\}
\arrow[r,shift left=1.2ex,"\text{staircase}"]
\&
\left\{\begin{tabular}{@{}c@{}}special coCartesian fibrations\\
$p:\mathcal E\to\simFin$\end{tabular}\right\}
\arrow[l,shift left=1.2ex,"\text{coCartesian strings}"]
\end{tikzcd}
\]
The \emph{staircase} of a left double fibration
$P:\mathbb E\to\bisimFin$ is the coCartesian fibration over
$\simFin$ whose objects are the objects of $\mathbb E$ and whose
morphisms are given by a horizontal morphism of $\mathbb E$ followed
by a vertical morphism of $\mathbb E$. Conversely, the
\emph{coCartesian strings} of a coCartesian fibration
$p:\mathcal E\to\simFin$ form the left double fibration whose
horizontal $m$-simplices are the strings of $m$ composable
coCartesian morphisms of $\mathcal E$. Moreover, the comparison
between left double fibrations and coCartesian fibrations does not use
any particular feature of $\Fin$, and it holds over the nerve of any
small category; see \cref{rem:general_base}.

\subsection*{Relation to other work}
The equivalence between left double fibrations and coCartesian
fibrations over $\Fin$ also follows from Nuiten's framework
\cite{Nuiten2023}. While his approach is formulated for double
$\infty$-categories and proceeds through a reflection of their
fibrations, we implement the comparison through a Quillen equivalence
between explicit model structures; see \cref{rem:nuiten_dictionary} for
more details, and \cite{BoavidaDeBrito2018} for the Grothendieck
construction over Segal objects. 

Several ingredients are taken from the literature. Double higher
categories as Segal objects are those of \cite{Moser2026InftyN}, and
the left fibration condition is that of \cite{Rasekh2021Yoneda}. The
staircase is the total simplicial set of Artin and Mazur
\cite{ArtinMazur1966VanKampen}, and we rely on
Stevenson's theorem that the total simplicial set is equivalent to the
diagonal in the Joyal model structure
\cite{Stevenson2025TotalSimplicialSet}, which sharpens the classical
statement for weak homotopy equivalences \cite{CegarraRemedios2005}.
The adjunction underlying our Quillen equivalence is a marked, fibered
relative of the adjunctions of Joyal and Tierney between simplicial
sets and bisimplicial sets \cite{JoyalTierney2007}; compare
\cite{Rasekh2021QuasiCategories}.

\subsection*{Organization of the paper}
In \cref{sec:general_preliminaries} we recall preliminaries on finite
pointed sets and on model categories. In
\cref{sec:internal_description,sec:enriched_description} we present the
homotopy theories of (special) left double fibrations and of (special)
coCartesian fibrations, respectively. In \cref{sec:comparisons} we
establish the equivalence between them, postponing a technical part of
the argument to \cref{app:staircase_extension}.

\subsection*{Acknowledgements}
The author is grateful to his advisor, Martina Rovelli, for her support, for many insightful conversations, and for helpful suggestions throughout this project.


\section{General preliminaries}
\label{sec:general_preliminaries}

In this section, we collect the preliminaries used throughout the
paper. After fixing conventions in \cref{subsec:global_conventions}, we
introduce the inert fans in \cref{subsec:prelim_finite_pointed_sets},
the evaluation of a simplicial object at a simplicial set in
\cref{subsec:prelim_reedy}, and the base $\bisimFin$ in
\cref{subsec:discrete_base}. We then recall in
\cref{subsec:prelim_model_categories} the facts about Reedy model
structures and left Bousfield localizations that we need, and we
specialize to simplicial sets in
\cref{subsec:prelim_simplicial_models}.

\subsection{Conventions}
\label{subsec:global_conventions}

We write $\mathbf M$ and $\mathbf N$ for model categories,
$[0]$ for the terminal object, $\varnothing$ for the initial
object, and $!_A:A\to[0]$ for the terminal map; an empty product is
$[0]$. Unless another model structure is named, weak equivalences,
cofibrations, and fibrations are those of the ambient model category,
and we decorate them as
\[
\begin{tikzcd}[ampersand replacement=\&]
A\arrow[r,two heads]\&B
\end{tikzcd}
\ \text{a fibration,}
\quad
\begin{tikzcd}[ampersand replacement=\&]
A\arrow[r,hook]\&B
\end{tikzcd}
\ \text{a cofibration,}
\quad
\begin{tikzcd}[ampersand replacement=\&]
A\arrow[r,"\sim"]\&B
\end{tikzcd}
\ \text{a weak equivalence.}
\]
All limits, colimits, and fibers are the ordinary ones unless called
homotopical; a \emph{point} of an object $X$ is a map $[0]\to X$, and a fiber is a pullback over a point. We write
$\amalg$ for binary coproducts and pushouts and $\coprod$ for
indexed ones.

Two model structures on simplicial sets occur throughout:
$\sset_{(\infty,0)}$, the Kan--Quillen structure, whose fibrant
objects are the Kan complexes, and $\sset_{(\infty,1)}$, the Joyal
structure, whose fibrant objects are the quasi-categories. Both are
recalled in \cref{subsec:prelim_simplicial_models}.
Simplicial sets are written $\mathcal K,\mathcal L,\ldots$ and
simplicial objects of $\mathbf M$ are written
$\mathbb X,\mathbb Y,\ldots$. The simplicial direction of such an
object is called \emph{horizontal}. As a rule $m$ is a horizontal
degree and $r$ is the vertical degree of a simplicial set occurring
as one of its levels; $n$ is always the dimension in
$(\infty,n)$, and $k$ is always an arity.

The notation $\RMap(X,Y)$ means the derived mapping space, or
homotopy function complex, of $X$ and $Y$ in the model category
containing them. It is defined up to weak equivalence in
$\sset_{(\infty,0)}$ and can always be taken to be a Kan complex
\cite[Definition~17.1.1 and Proposition~17.1.3]{hirschhorn2003model}.
We write $\mathbf M(X,Y)$ for the set of morphisms in the underlying
category. We never attach the ambient category as a subscript to
$\RMap$, $\Map$, $\Map^\flat$, or $\Map^\sharp$; the last
three are defined where they are needed. This convention causes no ambiguity
under left Bousfield localization: the cofibrations of
$L_{\mathsf S}\mathbf N$ are those of $\mathbf N$ and every
weak equivalence of $\mathbf N$ is one of $L_{\mathsf S}\mathbf N$,
so a cosimplicial resolution in $\mathbf N$ is still one in
$L_{\mathsf S}\mathbf N$, and the two structures compute the same
$\RMap(X,Y)$ whenever $X$ is cofibrant and $Y$ is fibrant in the
localization \cite[Definitions~16.1.2 and~17.1.1 and
Proposition~17.1.3]{hirschhorn2003model}.

\subsection{Finite pointed sets and inert maps}
\label{subsec:prelim_finite_pointed_sets}

We use the covariant convention for finite pointed sets, and denote by
$\Fin$ the skeletal category of finite pointed sets
$\langle k\rangle=\{*,1,\ldots,k\}$, for $k\geq0$, and basepoint
preserving maps. The nonnegative integer $k$ is the \emph{arity} of
$\langle k\rangle$. A morphism
$f:\langle k\rangle\to\langle\ell\rangle$ of $\Fin$ is \emph{inert}
if each fiber $f^{-1}(j)$, for $1\leq j\leq\ell$, has exactly one
element, and it is \emph{active} if $f^{-1}(*)=\{*\}$; see
\cite[Definitions~2.1.1.8 and~2.1.2.1]{LurieHA}. We will use the
following two families of maps.

\begin{notation}
\label{not:inert_and_fold}
For $k\geq0$, let $\mu_k:\langle k\rangle\to\langle1\rangle$ be
the active map sending every non-basepoint to $1$; for $k=0$ it is
the unique pointed map $\langle0\rangle\to\langle1\rangle$.
For $k\geq 1$ and $1\leq i\leq k$, let $\rho_i:\langle k\rangle\rightarrow\langle 1\rangle$
be the inert map defined by
\[
\rho_i(j)=
\begin{cases}
1, & j=i,\\
*, & j\neq i.
\end{cases}
\]
\end{notation}

We denote by $\simFin=N(\Fin)$ the nerve of $\Fin$, so that an
$m$-simplex of $\simFin$ is a string of $m$ composable maps in
$\Fin$. To encode the $k$ inert maps $\rho_i$ with common source
$\langle k\rangle$, we introduce the following simplicial set over
$\simFin$.

\begin{definition}
\label{def:inert_fan}
For $k\geq1$, the \emph{inert fan} $\mathcal J_k$ is the
simplicial set obtained by identifying the initial vertices of $k$
copies of $\Delta[1]$. Denote the common vertex by $c$ and the
terminal vertices by $v_1,\ldots,v_k$, and write
\[
\mathcal V_k=\coprod_{i=1}^k\Delta\{v_i\}\subseteq\mathcal J_k
\]
for the subcomplex they span.
The map $\beta_k:\mathcal J_k\to\simFin$ takes $c$ to
$\langle k\rangle$, each $v_i$ to $\langle1\rangle$, and
the edge $c\to v_i$ to $\rho_i$. For $k=0$, set
$
\mathcal J_0=\Delta[0],\ \mathcal V_0=\varnothing,
\ \text{and} \ \beta_0(c)=\langle0\rangle.
$
\end{definition}

Note that the vertices of $\mathcal V_k$ are distinct, although
$\beta_k$ sends all of them to $\langle1\rangle$. The inert fans
$\mathcal J_0$, $\mathcal J_1$, and $\mathcal J_3$, and the image
of $\mathcal J_3$ in $\simFin$, can be depicted as follows.
\[
\begin{tikzcd}[ampersand replacement=\&]c\end{tikzcd}
\qquad
\begin{tikzcd}[ampersand replacement=\&]c\arrow[r]\&v_1\end{tikzcd}
\qquad
\begin{tikzcd}[ampersand replacement=\&]
\&v_1\\c\arrow[ur]\arrow[r]\arrow[dr]\&v_2\\\&v_3
\end{tikzcd}
\qquad \qquad
\begin{tikzcd}[ampersand replacement=\&]
\&\langle1\rangle\\
\langle3\rangle\arrow[ur,"\rho_1"]\arrow[r,"\rho_2"]
\arrow[dr,"\rho_3"']\&\langle1\rangle\\\&\langle1\rangle.
\end{tikzcd}
\]
We will consider the restrictions along the inclusions of simplicial
sets over $\simFin$
\[
\begin{tikzcd}[ampersand replacement=\&]
\Delta\{c\}\arrow[r,hook]\&\mathcal J_k
\&\mathcal V_k\arrow[l,hook],
\end{tikzcd}
\]
where the left-hand inclusion is the identity when $k=0$.

\subsection{Horizontal shapes and Reedy fibrancy}
\label{subsec:prelim_reedy}

In this section, we study the evaluation of a simplicial object at a
simplicial set and its compatibility with the Reedy model structure.
Throughout, $\mathbf M$ denotes a model category. A simplicial object
$\mathbb X$ in $\mathbf M$ is a functor
$\Delta^{\mathrm{op}}\to\mathbf M$, with face and degeneracy maps
written $d_i^h,s_i^h$. Given a simplicial set $\mathcal K$,
regarded as a \emph{horizontal shape}, we define the object of diagrams
of shape $\mathcal K$ in $\mathbb X$ to be
\begin{equation}
\label{eq:global_shape_mapping}
\mathbb X(\mathcal K)\coloneqq
\lim_{(\Delta[m]\to \mathcal K)\in(\Delta/\mathcal K)^{\mathrm{op}}}\mathbb X_m,
\end{equation}
where $\Delta/\mathcal K$ denotes the category of simplices of
$\mathcal K$, and a morphism from $\sigma:\Delta[m]\to \mathcal K$
to $\tau:\Delta[\ell]\to \mathcal K$, i.e., a map
$\alpha:[m]\to[\ell]$ with $\tau\alpha=\sigma$, acts by
$\alpha_h^*:\mathbb X_\ell\to\mathbb X_m$. Given an object $A$ of
$\mathbf M$, we denote by $\mathcal K\boxtimes A$ the simplicial
object of $\mathbf M$ given by
\[
(\mathcal K\boxtimes A)_m=\coprod_{\mathcal K_m}A,
\]
with structure maps induced by those of $\mathcal K$. The universal
properties of the coproducts and of the limit give a natural
isomorphism
\[
\mathbf M(A,\mathbb X(\mathcal K))\cong
\mathbf M^{\Delta^{\mathrm{op}}}(\mathcal K\boxtimes A,\mathbb X).
\]
In particular, the functor $\mathbb X(-)$ sends colimits of
simplicial sets to limits in $\mathbf M$, and we have
$\mathbb X(\varnothing)\cong[0]$ and
$\mathbb X(\Delta[m])\cong\mathbb X_m$.

For $m\geq1$, we denote by
\[
I[m]\coloneqq\bigcup_{j=1}^{m}\Delta\{j-1,j\}\subseteq\Delta[m]
\]
the \emph{spine} of $\Delta[m]$, where $\Delta\{j-1,j\}$ is the
face spanned by the indicated vertices, and we set
$I[0]\coloneqq\Delta[0]$. Restriction along the inclusion
$I[m]\hookrightarrow\Delta[m]$ gives the \emph{Segal map}
\begin{equation}
\label{eq:horizontal_segal_map}
\varphi_m:\mathbb X_m\longrightarrow\mathbb X(I[m])
\cong\underbrace{\mathbb X_1\times_{\mathbb X_0}\cdots
\times_{\mathbb X_0}\mathbb X_1}_{m\text{ factors}},
\end{equation}
where, in each pullback, the target $d_0^h$ of an edge is matched
with the source $d_1^h$ of the next edge.

We always equip $\mathbf M^{\Delta^{\mathrm{op}}}$ with the Reedy
model structure for the usual Reedy structure on
$\Delta^{\mathrm{op}}$, whose weak equivalences are the degreewise
ones \cite[Theorem~15.3.4\textup{(1)} and
Example~15.3.5]{hirschhorn2003model}. For a simplicial object
$\mathbb X$ of $\mathbf M$, the matching objects are given by
\[
M_m\mathbb X\cong\mathbb X(\partial\Delta[m]),
\qquad
M_0\mathbb X\cong[0],
\]
and we write $L_m\mathbb X$ for the latching objects, so that
$L_0\mathbb X=\varnothing$. A map
$P:\mathbb E\to\mathbb B$ is a Reedy fibration exactly when every
relative matching map
\[
\mathbb E_m\longrightarrow
M_m\mathbb E\times_{M_m\mathbb B}\mathbb B_m
\]
is a fibration, and a map $\mathbb X\to\mathbb Y$ is a Reedy
cofibration exactly when every relative latching map
$\mathbb X_m\amalg_{L_m\mathbb X}L_m\mathbb Y\to\mathbb Y_m$ is a
cofibration. The next lemma extends the fibrancy of the relative
matching maps to every inclusion of simplicial sets.
\begin{lemma}[{cf.~\cite[Lemmas~C.5.9 and~C.2.3]{RiehlVerity2022}}]
\label{lem:cellular_restriction}
Let $P:\mathbb E\to\mathbb Y$ be a Reedy fibration in
$\mathbf M^{\Delta^{\mathrm{op}}}$. For every monomorphism
$i:\mathcal K\hookrightarrow\mathcal L$ of simplicial sets,
the relative restriction map
\[
R_i:\mathbb E(\mathcal L)\longrightarrow
\mathbb E(\mathcal K)\times_{\mathbb Y(\mathcal K)}\mathbb Y(\mathcal L)
\]
is a fibration in $\mathbf M$. In particular, if $\mathbb X$
is Reedy fibrant, then $\mathbb X(\mathcal K)$ is fibrant for every
simplicial set $\mathcal K$, and every inclusion
$\mathcal K\hookrightarrow\mathcal L$ induces a fibration
$\mathbb X(\mathcal L)\to\mathbb X(\mathcal K)$.
\end{lemma}
\begin{proof}
The evaluation $\mathbb X(\mathcal K)$ is the limit of $\mathbb X$
weighted by $\mathcal K$, and $R_i$ is the corresponding relative,
or Leibniz, limit map. A monomorphism of simplicial sets has
monomorphic relative latching maps by \cite[Example~6.9]{RV14}. Hence,
by the dual of \cite[Proposition~9.1]{RV14} applied to the Reedy
category $\Delta^{\mathrm{op}}$, the map $R_i$ is a fibration. The
last assertions follow by applying the first one to the map
$\mathbb X\to[0]$ and, for the fibrancy of
$\mathbb X(\mathcal K)$, to the inclusion
$\varnothing\hookrightarrow\mathcal K$.
\end{proof}

We conclude this section with a lemma about spans with a common
source, which will be applied to the spans induced by the inert maps
in \cref{prop:cleavage_free_segal,prop:marked_inert_span_fibrations}.

\begin{lemma}
\label{lem:common_source_fibrations}
Let $k\geq 0$, let $A$ be a fibrant object of $\mathbf M$, and let $D_1,\ldots,D_k$ be objects of $\mathbf M$. Suppose that, for each $1\leq i\leq k$, there is a span
\[
A\xleftarrow{\ s_i\ }H_i\xrightarrow{\ t_i\ }D_i
\]
such that $s_i$ is a trivial fibration and
\[
(s_i,t_i):H_i\twoheadrightarrow A\times D_i
\]
is a fibration. Define
\[
R=H_1\times_A\cdots\times_A H_k,
\qquad
\pi:R\longrightarrow A,
\qquad
t=(t_i\pr_i)_{i=1}^k:
R\longrightarrow\prod_{i=1}^k D_i.
\]
Then $\pi$ is a trivial fibration, $R$ is fibrant, and $t$ is a
fibration. For $k=0$, this statement is to be read with $R=A$,
$\pi=\id_A$, and $t={!_A}:A\rightarrow[0]$.
\end{lemma}

\begin{proof}
For $k\geq1$, put $R^{(0)}=A$ and
$R^{(i)}=R^{(i-1)}\times_A H_i$, so that each square
\[
\begin{tikzcd}[ampersand replacement=\&, sep=large]
R^{(i)} \arrow[r,"\pr_2"] \arrow[d,two heads,"\pr_1"',"\sim"]
  \arrow[dr,phantom,"\lrcorner",very near start]
\& H_i \arrow[d,two heads,"s_i","\sim"']\\
R^{(i-1)} \arrow[r]
\& A
\end{tikzcd}
\]
is a pullback. Since $s_i$ is a trivial fibration, so is the map
$R^{(i)}\to R^{(i-1)}$. The composite of these maps is
$\pi:R\to A$, and so $\pi$ is a trivial fibration and $R$ is
fibrant. Next, let
\[
\delta:A\times\prod_{i=1}^kD_i
  \longrightarrow\prod_{i=1}^k(A\times D_i)
\]
be the map whose $i$-th component is
$\langle\pr_A,\pr_{D_i}\rangle$. There is a pullback square
\[
\begin{tikzcd}[ampersand replacement=\&, sep=large]
R \arrow[r,"{(\pr_i)_i}"] \arrow[d,two heads,"{(\pi,t)}"']
  \arrow[dr,phantom,"\lrcorner",very near start]
\& \prod_{i=1}^kH_i
  \arrow[d,two heads,"{\prod_{i=1}^k(s_i,t_i)}"]\\
A\times\prod_{i=1}^kD_i \arrow[r,"\delta"']
\& \prod_{i=1}^k(A\times D_i),
\end{tikzcd}
\]
and so the map $(\pi,t)$ is a fibration, as fibrations are stable
under products and pullbacks. Moreover, the projection
$A\times\prod_{i=1}^kD_i\to\prod_{i=1}^kD_i$ is a fibration, as it
is a pullback of the fibration $!_A:A\to[0]$. Hence the map $t$,
which is the composite of $(\pi,t)$ with this projection, is a
fibration. For $k=0$, the assertions reduce to the fibrancy of
$A$.
\end{proof}

\subsection{The discrete finite-pointed-set base}
\label{subsec:discrete_base}

In this section, we introduce the simplicial object $\bisimFin$ of
$\mathbf M$ over which the fibrations of
\cref{sec:internal_description} are defined. We denote by
\[
(-)^{\mathrm{disc}}:\Set\longrightarrow\mathbf M,
\qquad
S\longmapsto S^{\mathrm{disc}}\coloneqq\coprod_{s\in S}[0],
\]
the \emph{discrete realization} functor, which sends a function of
sets $S\to T$ to the map $S^{\mathrm{disc}}\to T^{\mathrm{disc}}$
that maps the summand indexed by $s$ identically onto the summand
indexed by the image of $s$. Given a simplicial set $\mathcal K$,
we write $\mathcal K^h\coloneqq\mathcal K\boxtimes[0]$, and we define
$\bisimFin\coloneqq(\simFin)^h$, so that
\[
(\bisimFin)_m=((\simFin)_m)^{\mathrm{disc}},
\]
with structure maps induced by those of $\simFin$.

Recall that small coproducts in $\mathbf M$ are \emph{disjoint} if
every coproduct injection is a monomorphism and the pullback of two
distinct coproduct injections is the initial object, and
\emph{universal} if they are preserved by pullback along every map.
Both properties hold in every presheaf category.

\begin{proposition}
\label{prop:discrete_base}
Suppose that small coproducts in $\mathbf M$ are disjoint and
universal, and that every map
$S^{\mathrm{disc}}\to T^{\mathrm{disc}}$ induced by a function of
sets is a fibration. Then the discrete realization functor preserves
finite limits, the simplicial object $\bisimFin$ is Reedy fibrant,
and its Segal maps $\varphi_m$ from \eqref{eq:horizontal_segal_map}
are isomorphisms.
\end{proposition}
\begin{proof}
Given functions of sets $S\to U\leftarrow T$, universality
decomposes the pullback
$S^{\mathrm{disc}}\times_{U^{\mathrm{disc}}}T^{\mathrm{disc}}$ as the
coproduct, indexed by the pairs $(s,t)\in S\times T$, of the
pullbacks of the corresponding summands. By disjointness, the pullback
indexed by $(s,t)$ is $[0]$ if $s$ and $t$ have the same image
in $U$, and it is the initial object otherwise. Hence
\[
S^{\mathrm{disc}}\times_{U^{\mathrm{disc}}}T^{\mathrm{disc}}
\cong\coprod_{(s,t)\in S\times_U T}[0]
=(S\times_U T)^{\mathrm{disc}},
\]
and so discrete realization preserves finite limits, as it also
preserves the terminal object. Since the matching category of
$\Delta^{\mathrm{op}}$ in each degree is finite, we get an
isomorphism
\[
M_m\bisimFin\cong(M_m\simFin)^{\mathrm{disc}},
\]
under which the matching map of $\bisimFin$ is induced by a function
of sets. Hence the matching map is a fibration by hypothesis, and
$\bisimFin$ is Reedy fibrant. Finally, the nerve $\simFin$ has the
Segal bijections
\[
(\simFin)_m\cong(\simFin)_1
\times_{(\simFin)_0}\cdots\times_{(\simFin)_0}(\simFin)_1,
\]
which discrete realization sends to the Segal isomorphisms of
$\bisimFin$.
\end{proof}

\begin{rmk}
\label{rem:discrete_base_hypotheses}
The fibration hypothesis in \cref{prop:discrete_base} does not follow
from the axioms of a model category. For instance, in the left
Bousfield localization of $\sset_{(\infty,0)}$ at the fold map
$\partial\Delta[1]\to\Delta[0]$, a Kan complex $\mathcal Z$ is
local if and only if its diagonal
$\mathcal Z\to\mathcal Z\times\mathcal Z$ is a weak equivalence,
i.e., if and only if $\mathcal Z$ is empty or weakly contractible.
Hence the object $\coprod_{\mathbb N}\Delta[0]=(\bisimFin)_0$ is not
fibrant in this localization. However, the hypothesis holds in the
models considered in this paper; see
\cref{prop:quasicategorical_standing_assumptions} and the examples
following \cref{prop:presheaf_discrete_criterion}.
\end{rmk}

\subsection{Fixed-shape fibers, derived mapping spaces, and localization}
\label{subsec:prelim_model_categories}

In this section, we fix a simplicial object $\mathbb B$ of
$\mathbf M$, which will be $\bisimFin$ from
\cref{sec:internal_description} onwards. We introduce the fibers of
the evaluations $\mathbb E(\mathcal K)\to\mathbb B(\mathcal K)$ of a
simplicial object $\mathbb E$ over $\mathbb B$, and we recall the
facts about derived mapping spaces and left Bousfield localizations
needed to construct the model structures of
\cref{sec:internal_description}. Recall that the weak equivalences,
cofibrations, and fibrations of the slice model structure
$\genreedyslice$ are created in the Reedy model structure
\cite[Theorem~7.6.5\textup{(2)}]{hirschhorn2003model}. In particular,
an object $P:\mathbb E\to\mathbb B$ of $\genreedyslice$ is fibrant
if and only if $P$ is a Reedy fibration.

\begin{definition}
\label{def:fixed_shape_fiber}
Let $P:\mathbb E\to\mathbb B$ be a map of simplicial objects of
$\mathbf M$, let $\mathcal K$ be a simplicial set, and let
$\beta:[0]\to\mathbb B(\mathcal K)$ be a point of
$\mathbb B(\mathcal K)$. The \emph{fixed-shape fiber}
$\mathbb E(\mathcal K)_\beta$ of $\mathbb E$ at $\beta$ is the
pullback
\[
\begin{tikzcd}[ampersand replacement=\&, sep=large]
\mathbb E(\mathcal K)_\beta\arrow[r,"\pr_1"]\arrow[d]
  \arrow[dr,phantom,"\lrcorner",very near start]
\&\mathbb E(\mathcal K)\arrow[d,"P(\mathcal K)"]\\
{[0]}\arrow[r,"\beta"']\&\mathbb B(\mathcal K).
\end{tikzcd}
\]
\end{definition}

\begin{lemma}
\label{lem:fixed_shape_fiber_fibrant}
Let $P:\mathbb E\to\mathbb B$ be a Reedy fibration. Then every
fixed-shape fiber of $\mathbb E$ is fibrant. Moreover, given a
monomorphism $i:\mathcal K\hookrightarrow\mathcal L$ of simplicial
sets and a point $\gamma:[0]\to\mathbb B(\mathcal L)$ with
$\beta=\mathbb B(i)\gamma$, the restriction map
$i^*:\mathbb E(\mathcal L)_\gamma\to\mathbb E(\mathcal K)_\beta$ is
a fibration.
\end{lemma}
\begin{proof}
By \cref{lem:cellular_restriction} applied to the inclusion
$\varnothing\hookrightarrow\mathcal K$, the map $P(\mathcal K)$ is
a fibration, and so its fiber $\mathbb E(\mathcal K)_\beta$ is
fibrant. For the second assertion, the map $i^*$ is a pullback of
the fibration $R_i$ of \cref{lem:cellular_restriction}, as displayed
in the following diagram.
\[
\begin{tikzcd}[ampersand replacement=\&, sep=large]
\mathbb E(\mathcal L)_\gamma\arrow[r,"\pr_1"]
  \arrow[d,two heads,"i^*"']
  \arrow[dr,phantom,"\lrcorner",very near start]
\&\mathbb E(\mathcal L)\arrow[d,two heads,"R_i"]\\
\mathbb E(\mathcal K)_\beta\arrow[r,"{(\pr_1,\gamma !)}"']
\&\mathbb E(\mathcal K)\times_{\mathbb B(\mathcal K)}\mathbb B(\mathcal L).
\end{tikzcd}
\]
\end{proof}

\begin{notation}
\label{not:fixed_simplex_fibers}
A map of simplicial sets $\beta:\mathcal K\to\simFin$ determines a
point $[0]\to\bisimFin(\mathcal K)$, whose component at a simplex
$\tau:\Delta[m]\to\mathcal K$ is the inclusion of the summand of
$(\bisimFin)_m$ indexed by $\beta\tau$. Given a map
$P:\mathbb E\to\bisimFin$, we denote the fixed-shape fiber of
$\mathbb E$ at this point by $\mathbb E(\mathcal K)_\beta$. In
particular, for an $m$-simplex $\sigma$ of $\simFin$, we write
\[
\mathbb E_m(\sigma)\coloneqq\mathbb E(\Delta[m])_\sigma
=\mathbb E_m\times_{(\bisimFin)_m}[0].
\]
\end{notation}

The fixed-shape fiber functor admits a left adjoint, which will be
used to construct the maps at which we localize.
\begin{lemma}
\label{lem:general_shape_quillen}
Given a simplicial set $\mathcal K$ and a point
$\beta:[0]\to\mathbb B(\mathcal K)$, there is a Quillen adjunction
\[
\adjunction{\mathbf M}{\genreedyslice}
  {(\mathcal K\boxtimes-)_{\beta}}
  {\mathbb E\mapsto\mathbb E(\mathcal K)_\beta}
\]
where, for an object $A$ of $\mathbf M$, the structure map of
$(\mathcal K\boxtimes A)_\beta$ is given on the summand indexed by
$\sigma\in\mathcal K_m$ by the composite
$A\xrightarrow{!_A}[0]\xrightarrow{\beta_\sigma}\mathbb B_m$. In
particular, if $A$ is cofibrant, then
$(\mathcal K\boxtimes A)_\beta$ is cofibrant in $\genreedyslice$.
\end{lemma}
\begin{proof}
Under the natural isomorphism
$\mathbf M(A,\mathbb E(\mathcal K))\cong
\mathbf M^{\Delta^{\mathrm{op}}}(\mathcal K\boxtimes A,\mathbb E)$
of \cref{subsec:prelim_reedy}, a map $a:A\to\mathbb E(\mathcal K)$
factors through $\mathbb E(\mathcal K)_\beta$ if and only if
$P(\mathcal K)a=\beta\,!_A$, i.e., if and only if the adjunct of
$a$ is a map over $\mathbb B$. Hence we get a natural isomorphism
\[
\genreedyslice((\mathcal K\boxtimes A)_\beta,\mathbb E)
\cong\mathbf M(A,\mathbb E(\mathcal K)_\beta).
\]
To show that the left adjoint is left Quillen, let
$\mathsf D_m\mathcal K$ and $\mathsf N_m\mathcal K$ denote the sets
of degenerate and nondegenerate $m$-simplices of $\mathcal K$,
respectively. Since every simplex is uniquely a degeneracy of a
nondegenerate simplex, we have
$L_m(\mathcal K\boxtimes A)\cong\coprod_{\mathsf D_m\mathcal K}A$;
see \cite[Lemma~15.8.4]{hirschhorn2003model}. Hence, given a map
$i:A\to A'$ in $\mathbf M$, the source of the relative latching map
of $\mathcal K\boxtimes i$ in degree $m$ is
\[
(\mathcal K\boxtimes A)_m
  \amalg_{L_m(\mathcal K\boxtimes A)}L_m(\mathcal K\boxtimes A')
\cong
\left(\coprod_{\mathsf D_m\mathcal K}A'\right)
  \amalg\left(\coprod_{\mathsf N_m\mathcal K}A\right),
\]
and the relative latching map to $(\mathcal K\boxtimes A')_m$ is the
identity on the summands indexed by $\mathsf D_m\mathcal K$ and the
map $i$ on the summands indexed by $\mathsf N_m\mathcal K$. Since
cofibrations and trivial cofibrations are stable under coproducts, if
$i$ is a cofibration, then $\mathcal K\boxtimes i$ is a Reedy
cofibration, and if $i$ is moreover a weak equivalence, then so is
each level $\coprod_{\mathcal K_m}i$ of $\mathcal K\boxtimes i$.
As the cofibrations and weak equivalences of $\genreedyslice$ are
created in the Reedy model structure, the left adjoint preserves
cofibrations and trivial cofibrations, and hence cofibrant objects.
\end{proof}

\begin{lemma}
\label{lem:shape_inclusion_cofibrant}
Let $i:\mathcal K\hookrightarrow\mathcal L$ be a monomorphism
of simplicial sets and let $A$ be a cofibrant object of
$\mathbf M$. Then the map $i\boxtimes A$ is a Reedy cofibration
between Reedy cofibrant objects. Moreover, given a point
$\gamma:[0]\to\mathbb B(\mathcal L)$ with
$\beta=\mathbb B(i)\gamma$, the map
$(\mathcal K\boxtimes A)_\beta\to(\mathcal L\boxtimes A)_\gamma$
induced by $i$ is a cofibration between cofibrant objects of
$\genreedyslice$.
\end{lemma}
\begin{proof}
We identify $\mathcal K$ with its image in $\mathcal L$. A simplex
of $\mathcal K$ is degenerate in $\mathcal L$ if and only if it is
degenerate in $\mathcal K$, as a simplex $x$ of $\mathcal L$ with
$s_jx\in\mathcal K$ is the face $d_js_jx$ of a simplex of
$\mathcal K$. Hence
$\mathcal K_m\cap\mathsf D_m\mathcal L=\mathsf D_m\mathcal K$ and,
by the computation of the latching objects in the proof of
\cref{lem:general_shape_quillen}, the relative latching map of
$i\boxtimes A$ in degree $m$ is the map
\[
\coprod_{\mathcal K_m\cup\mathsf D_m\mathcal L}A
\longrightarrow\coprod_{\mathcal L_m}A,
\]
which is a coproduct of identity maps and copies of the cofibration
$\varnothing\to A$. Hence $i\boxtimes A$ is a Reedy cofibration.
The same argument applied to the inclusions
$\varnothing\hookrightarrow\mathcal K$ and
$\varnothing\hookrightarrow\mathcal L$ shows that
$\mathcal K\boxtimes A$ and $\mathcal L\boxtimes A$ are Reedy
cofibrant. The last assertion follows from the fact that the
cofibrations of $\genreedyslice$ are created in the Reedy model
structure.
\end{proof}

We compute derived mapping spaces using cosimplicial resolutions
\cite[Definition~17.1.1 and Theorem~17.4.14]{hirschhorn2003model}.
Recall that, given a cofibrant object $X$ and a fibrant object $Y$
of a model category $\mathbf N$, the derived mapping space
$\RMap(X,Y)$ is modeled by the Kan complex
$\mathbf N(X^\bullet,Y)$, where
$X^\bullet\to\operatorname{const}(X)$ is a cosimplicial resolution
of $X$, i.e., a Reedy cofibrant cosimplicial object $X^\bullet$
with a levelwise weak equivalence to the constant cosimplicial object
at $X$ \cite[Definition~16.1.2\textup{(1)} and
Proposition~17.1.3]{hirschhorn2003model}. The adjunction of
\cref{lem:general_shape_quillen} is compatible with derived mapping
spaces as follows.
\begin{lemma}
\label{lem:derived_shape_adjunction}
Let $A$ be a cofibrant object of $\mathbf M$ and let
$P:\mathbb E\to\mathbb B$ be a Reedy fibration.
\begin{enumerate}[label=(\arabic*)]
\item Given a simplicial set $\mathcal K$ and a point
$\beta:[0]\to\mathbb B(\mathcal K)$, there is a weak equivalence
in $\sset_{(\infty,0)}$
\[
\RMap((\mathcal K\boxtimes A)_\beta,\mathbb E)
\simeq\RMap(A,\mathbb E(\mathcal K)_\beta).
\]
\item Given a map $i:\mathcal K\to\mathcal L$ of simplicial sets and
a point $\gamma:[0]\to\mathbb B(\mathcal L)$ with
$\beta=\mathbb B(i)\gamma$, the weak equivalences of \textup{(1)}
identify precomposition with $i\boxtimes A$ with postcomposition
with
$i^*:\mathbb E(\mathcal L)_\gamma\to\mathbb E(\mathcal K)_\beta$.
\end{enumerate}
\end{lemma}
\begin{proof}
Let $A^\bullet\to\operatorname{const}(A)$ be a cosimplicial
resolution of $A$. The left Quillen functor
$(\mathcal K\boxtimes-)_\beta$ of \cref{lem:general_shape_quillen}
sends $A^\bullet$ to a cosimplicial resolution of
$(\mathcal K\boxtimes A)_\beta$, and the adjunction gives an
isomorphism of Kan complexes
\[
\genreedyslice
 ((\mathcal K\boxtimes A^\bullet)_\beta,\mathbb E)
\cong
\mathbf M(A^\bullet,\mathbb E(\mathcal K)_\beta);
\]
see \cite[Proposition~17.4.16\textup{(1)}]{hirschhorn2003model}. This
isomorphism proves \textup{(1)}. For \textup{(2)}, by naturality of
the adjunction in each cosimplicial degree, the square
\[
\begin{tikzcd}[ampersand replacement=\&, sep=large]
\RMap((\mathcal L\boxtimes A)_\gamma,\mathbb E)
  \arrow[r,"\cong"]\arrow[d,"(i\boxtimes A)^*"']
\&
\RMap(A,\mathbb E(\mathcal L)_\gamma)
  \arrow[d,"{\RMap(A,i^*)}"]\\
\RMap((\mathcal K\boxtimes A)_\beta,\mathbb E)
  \arrow[r,"\cong"']
\&\RMap(A,\mathbb E(\mathcal K)_\beta)
\end{tikzcd}
\]
commutes when the derived mapping spaces are computed using the
resolution $A^\bullet$ for both $\mathcal K$ and $\mathcal L$.
Since derived mapping spaces are independent of the chosen
resolutions, this proves \textup{(2)}.
\end{proof}

We next identify the Reedy cofibrations when $\mathbf M$ is a
presheaf category whose cofibrations are the monomorphisms.
\begin{lemma}
\label{lem:presheaf_reedy_cofibrancy}
Suppose that the underlying category of $\mathbf M$ is a presheaf
category $\mathbf{Set}^{\mathrm C^{\mathrm{op}}}$, with
$\mathrm C$ small, and that the cofibrations of $\mathbf M$ are
precisely the monomorphisms. Then the Reedy cofibrations in
$\mathbf M^{\Delta^{\mathrm{op}}}$ are precisely the levelwise
monomorphisms. In particular, every object of
$\mathbf M^{\Delta^{\mathrm{op}}}$ and every object of
$\genreedyslice$ is cofibrant.
\end{lemma}
\begin{proof}
Since limits and colimits of presheaves are computed objectwise, for
every object $c$ of $\mathrm C$, the latching map
$(L_m\mathbb X)(c)\to\mathbb X_m(c)$ is the inclusion of the set of
degenerate simplices; compare
\cite[Lemma~15.8.4 and Proposition~15.8.6]{hirschhorn2003model}.
Let $\mathbb X\hookrightarrow\mathbb Y$ be a levelwise
monomorphism. As in the proof of \cref{lem:shape_inclusion_cofibrant},
a simplex of $\mathbb X$ that is degenerate in $\mathbb Y$ is
degenerate in $\mathbb X$, and so
\[
\mathbb X_m(c)\cap(L_m\mathbb Y)(c)=(L_m\mathbb X)(c).
\]
Hence the source of the relative latching map is, at $c$, the union
$\mathbb X_m(c)\cup(L_m\mathbb Y)(c)$ inside $\mathbb Y_m(c)$, and
so the relative latching map is a monomorphism. Conversely, a Reedy
cofibration is levelwise a cofibration by
\cite[Proposition~15.3.11\textup{(1)}]{hirschhorn2003model}. The last
assertion follows by applying this characterization to the maps from
the initial object, as the cofibrations of $\genreedyslice$ are
created in the Reedy model structure.
\end{proof}

We now recall the left Bousfield localization of a model category
$\mathbf N$ at a small set $\mathsf S$ of maps in $\mathbf N$.
A fibrant object $Z$ of $\mathbf N$ is \emph{$\mathsf S$-local}
if, for every map $s:A\to D$ in $\mathsf S$, the induced map
\[
s^*:\RMap(D,Z)
\longrightarrow\RMap(A,Z)
\]
is a weak equivalence in $\sset_{(\infty,0)}$, and a map
$f:X\to Y$ in $\mathbf N$ is an
\emph{$\mathsf S$-local equivalence} if the induced map
$f^*:\RMap(Y,Z)\to\RMap(X,Z)$ is a weak equivalence in
$\sset_{(\infty,0)}$ for every $\mathsf S$-local object $Z$; see
\cite[Definition~3.1.4\textup{(1)}]{hirschhorn2003model}.
\begin{theorem}[{\cite[Theorem~4.7]{barwick2010bousfield}}]
\label{thm:left_bousfield_localization}
Let $\mathbf N$ be a left proper combinatorial model category and
let $\mathsf S$ be a small set of maps in $\mathbf N$. Then the
left Bousfield localization $L_{\mathsf S}\mathbf N$ of $\mathbf N$
at $\mathsf S$ exists, and it is left proper and combinatorial. Its
cofibrations are the cofibrations of $\mathbf N$, its weak
equivalences are the $\mathsf S$-local equivalences, and its fibrant
objects are the $\mathsf S$-local objects.
\end{theorem}

\begin{proposition}
\label{prop:reedy_localization_existence}
If $\mathbf M$ is left proper and combinatorial, then so is
$\genreedyslice$. In particular, the left Bousfield localization of
$\genreedyslice$ at every small set of maps exists.
\end{proposition}
\begin{proof}
The Reedy model structure on $\mathbf M^{\Delta^{\mathrm{op}}}$ is
left proper by \cite[Theorem~15.3.4\textup{(2)}]{hirschhorn2003model}
and combinatorial by \cite[Lemma~3.33]{barwick2010bousfield}. The
underlying category of $\genreedyslice$ is locally presentable by
\cite[Proposition~1.57]{adamek_rosicky1994}, and the slice model
structure is cofibrantly generated and left proper by
\cite[Theorem~1.7\textup{(1), (3)}]{Hirschhorn2015Overcategories}.
The last assertion then follows from
\cref{thm:left_bousfield_localization}.
\end{proof}

We collect here the standard facts about model categories and left
Bousfield localizations that are used in the proofs below.

\begin{proposition}[Standard facts]
\label{prop:model_category_recollections}
Let $\mathbf M$ be a model category. Then the following statements hold.
\begin{enumerate}[label=\textup{(\arabic*)}]
\item Fibrations and trivial fibrations are closed under pullback,
composition, and retract. A map is a fibration if and only if it has
the right lifting property against every trivial cofibration; it is
a trivial fibration if and only if it has the right lifting property
against every cofibration
\cite[Lemma~1.1.10 and Corollary~1.1.11]{hovey1999model}.
\item A right Quillen functor preserves weak equivalences between
fibrant objects, and a left Quillen functor preserves weak
equivalences between cofibrant objects
\cite[Lemma~1.1.12]{hovey1999model}.
\item A morphism becomes an isomorphism in $\Ho(\mathbf M)$ if and
only if it is a weak equivalence, and morphisms from a cofibrant to a
fibrant object are represented by maps modulo homotopy
\cite[Theorem~1.2.10]{hovey1999model};
see also \cite[Theorem~8.3.10]{hirschhorn2003model}.
\item Derived mapping spaces are invariant under weak equivalences in
either variable \cite[Theorem~17.6.3]{hirschhorn2003model} and are
independent of the chosen resolutions
\cite[Theorem~17.4.14]{hirschhorn2003model}.
\item Finite products preserve weak equivalences between fibrant
objects, and the ordinary product of finitely many fibrant objects
represents their product in $\Ho(\mathbf M)$. Indeed, the diagonal
$\mathbf M\to\mathbf M^{\mathsf I}$ into the objectwise model
structure on a finite family preserves cofibrations and trivial
cofibrations, so the product functor is right Quillen, and
\textup{(2)} applies.
\item For a map $b:B'\to B$, pullback $b^*$ between the slice
model categories is right Quillen: its left adjoint is composition
with $b$, and cofibrations and trivial cofibrations of a slice are
detected in $\mathbf M$. Hence $b^*$ preserves weak equivalences
between fibrant objects of the slice, by \textup{(2)}.
\end{enumerate}
\end{proposition}

\begin{proposition}[{\cite[Theorem~3.2.13 and
Propositions~3.3.4, 3.3.16, 3.3.18]{hirschhorn2003model}}]
\label{prop:localization_recollections}
Suppose that $L_{\mathsf S}\mathbf N$ exists.
\begin{enumerate}[label=\textup{(\arabic*)}]
\item The identity functor $\mathbf N\to L_{\mathsf S}\mathbf N$
is left Quillen. A map between fibrant local objects is a weak
equivalence, respectively a fibration, in the localization if and only
if it is one in $\mathbf N$.
\item Let $F:\mathbf N\to\mathbf M$ be left Quillen whose derived
functor inverts every element of $\mathsf S$. Then
$F:L_{\mathsf S}\mathbf N\to\mathbf M$ is left Quillen. When the
domains and codomains of $\mathsf S$ are cofibrant, it is enough
that $F$ send each element of $\mathsf S$ to a weak equivalence.
\end{enumerate}
\end{proposition}

Recall that a \emph{cartesian model category} is a cartesian closed
model category satisfying the pushout-product and unit axioms
\cite[\S4.2]{hovey1999model}, where the pushout-product of maps
$i:A\to B$ and $j:C\to D$ is the map
\[
i\mathbin{\square}j:
(B\times C)\amalg_{A\times C}(A\times D)\longrightarrow B\times D.
\]
In a cartesian model category, the product with a cofibrant object
$A$ is a left Quillen functor, as it is the pushout-product with the
cofibration $\varnothing\to A$. Hence its right adjoint $(-)^A$
preserves fibrations, trivial fibrations, and weak equivalences
between fibrant objects.

\begin{notation}[{\cite[Corollary~4.33]{barwick2010bousfield}}]
\label{not:homotopy_generators}
A set of \emph{homotopy generators} of a model category $\mathbf N$
is a small set $\mathsf H$ of objects of $\mathbf N$ such that
every object of $\mathbf N$ is weakly equivalent to a small homotopy
colimit of objects of $\mathsf H$. Every combinatorial model
category admits a set of cofibrant homotopy generators, as a
consequence of Dugger's presentation theorem \cite{Dug01}. In
particular, so do its slice model categories, which are again
combinatorial.
\end{notation}

\begin{lemma}
\label{lem:homotopy_generators_detect}
Let $\mathsf H$ be a set of cofibrant homotopy generators in a
model category $\mathbf N$. A map $v:U\to V$ between fibrant
objects is a weak equivalence if and only if
\[
\RMap(A,v):\RMap(A,U)\longrightarrow\RMap(A,V)
\]
is a weak equivalence in $\sset_{(\infty,0)}$ for every
$A\in\mathsf H$.
\end{lemma}
\begin{proof}
If $v$ is a weak equivalence, then so is $\RMap(A,v)$ for every
$A\in\mathsf H$, by \cref{prop:model_category_recollections}(4).
Conversely, suppose that $\RMap(A,v)$ is a weak equivalence for
every $A\in\mathsf H$, and let $X$ be an object of $\mathbf N$.
Given a presentation
$X\simeq\operatorname{hocolim}_{d\in\mathrm D}A_d$ with
$A_d\in\mathsf H$, the map $\RMap(X,v)$ is identified with the
induced map between the homotopy limits of the diagrams
$\RMap(A_d,U)$ and $\RMap(A_d,V)$ by
\cite[Theorem~19.4.4\textup{(1)}]{hirschhorn2003model}, and hence it
is a weak equivalence by
\cite[Theorem~19.4.2\textup{(2)}]{hirschhorn2003model}. Taking
connected components, we get natural bijections
\[
\Ho(\mathbf N)(X,U)\longrightarrow\Ho(\mathbf N)(X,V),
\]
and so $v$ is an isomorphism in $\Ho(\mathbf N)$ by the Yoneda
lemma. Hence $v$ is a weak equivalence by
\cref{prop:model_category_recollections}(3).
\end{proof}

\begin{lemma}
\label{lem:shape_locality}
Let $\mathsf H$ be a set of cofibrant homotopy generators of
$\mathbf M$, let $i:\mathcal K\to\mathcal L$ be a map of
simplicial sets, and let $\gamma:[0]\to\mathbb B(\mathcal L)$ be a
point with $\beta=\mathbb B(i)\gamma$. A Reedy fibration
$P:\mathbb E\to\mathbb B$ is local with respect to the maps
\[
(i\boxtimes A):
(\mathcal K\boxtimes A)_\beta
\longrightarrow(\mathcal L\boxtimes A)_\gamma,
\qquad A\in\mathsf H,
\]
if and only if the restriction map
\[
i^*:\mathbb E(\mathcal L)_\gamma
\longrightarrow\mathbb E(\mathcal K)_\beta
\]
is a weak equivalence in $\mathbf M$. The same criterion holds for
a fibrant object of any left Bousfield localization of
$\genreedyslice$.
\end{lemma}
\begin{proof}
By \cref{lem:derived_shape_adjunction}, the map induced by
$i\boxtimes A$ on derived mapping spaces into $\mathbb E$ is
identified with the map $\RMap(A,i^*)$. Since the fixed-shape fibers
$\mathbb E(\mathcal L)_\gamma$ and $\mathbb E(\mathcal K)_\beta$
are fibrant by \cref{lem:fixed_shape_fiber_fibrant}, the first
assertion follows from \cref{lem:homotopy_generators_detect}. The last
assertion holds since the derived mapping spaces into a fibrant object
of a left Bousfield localization of $\genreedyslice$ can be computed
in $\genreedyslice$, as explained in
\cref{subsec:global_conventions}.
\end{proof}

\subsection{The quasi-categorical specialization}
\label{subsec:prelim_simplicial_models}

In this section, we recall the main facts about the Kan--Quillen and
Joyal model structures on simplicial sets, and we specialize the
constructions of the previous sections to the case
$\mathbf M=\sset_{(\infty,1)}$.

The Kan--Quillen model structure $\sset_{(\infty,0)}$ has the
monomorphisms as cofibrations, the maps with the right lifting
property against every horn inclusion as fibrations, and the Kan
complexes as fibrant objects
\cite[Theorem~7.10.12]{hirschhorn2003model}.
A map of simplicial sets is an \emph{inner fibration} if it has the
right lifting property against every inner horn inclusion
$\Lambda^a[r]\hookrightarrow\Delta[r]$, with $r\geq2$ and
$0<a<r$, and a \emph{quasi-category} is a simplicial set whose map
to $\Delta[0]$ is an inner fibration
\cite[Definitions~1.1.2.4 and~2.0.0.3]{lurie2009higher}. The Joyal
model structure $\sset_{(\infty,1)}$ has the monomorphisms as
cofibrations and the quasi-categories as fibrant objects, and it is
left proper and combinatorial
\cite[Theorems~2.2.5.1 and~2.4.6.1]{lurie2009higher}. Between
quasi-categories, its weak equivalences are the fully faithful and
essentially surjective maps
\cite[Tag~01JX, Theorem~4.8.4.1]{kerodon}, and its fibrations are the
inner fibrations that lift equivalences with specified source
\cite[Corollary~2.4.6.5]{lurie2009higher}. The two model structures
have the same trivial fibrations, namely the maps with the right
lifting property against every boundary inclusion
$\partial\Delta[r]\hookrightarrow\Delta[r]$, which we call
\emph{trivial fibrations of simplicial sets}. In both model
structures, every object is cofibrant and the terminal object is
$\Delta[0]$.

Moreover, the model category $\sset_{(\infty,1)}$ is cartesian, as
the functor $\mathcal A\times(-)$ preserves monomorphisms and weak
equivalences in $\sset_{(\infty,1)}$ for every simplicial set
$\mathcal A$ \cite[Corollary~2.2.5.4]{lurie2009higher}. Hence the
functor $(-)^{\mathcal A}$ preserves fibrations, trivial fibrations,
and weak equivalences between fibrant objects; see also
\cite[Proposition~8.5.7\textup{(2)}]{hirschhorn2003model}.

For $\mathbf M=\sset_{(\infty,1)}$, the discrete realization
$S^{\mathrm{disc}}$ of a set $S$ is the constant simplicial set at
$S$, and a simplicial object $\mathbb E$ of $\mathbf M$ is a
bisimplicial set, with $\mathbb E_{m,r}=(\mathbb E_m)_r$. We write
the horizontal degree first and denote the vertical face and
degeneracy maps by $d_i^v,s_i^v$. Given simplicial sets
$\mathcal K$ and $\mathcal A$, we have
\[
(\mathcal K\boxtimes\mathcal A)_{m,r}=\mathcal K_m\times\mathcal A_r,
\qquad\mathcal K^h=\mathcal K\boxtimes\Delta[0],
\qquad\mathbb E(\mathcal K)_r
=\sset(\mathcal K,\mathbb E_{\bullet,r}).
\]
Let $\mathcal B$ be a simplicial set and let
$\mathbb B=\mathcal B^h$. Every map between discrete simplicial sets
is a fibration in $\sset_{(\infty,1)}$: it is an inner fibration, as
inner horns are connected, and it lifts equivalences with specified
source, as every equivalence in a discrete simplicial set is an
identity. Since discrete realization preserves limits, the matching
maps of $\mathbb B$ are maps between discrete simplicial sets, and
so $\mathbb B$ is Reedy fibrant. Moreover, we have
\[
\mathbb B(\mathcal K)\cong
\sset(\mathcal K,\mathcal B)^{\mathrm{disc}},
\]
so that a point of $\mathbb B(\mathcal K)$ is a map
$\beta:\mathcal K\to\mathcal B$, and the structure map of
$(\mathcal K\boxtimes\mathcal A)_\beta$ is the composite
\[
\mathcal K\boxtimes\mathcal A\longrightarrow\mathcal K^h
\xrightarrow{\beta^h}\mathcal B^h.
\]
Given a map $\mathbb E\to\mathbb B$, we obtain a decomposition
\[
\mathbb E(\mathcal K)\cong
\coprod_{\beta:\mathcal K\to\mathcal B}\mathbb E(\mathcal K)_\beta.
\]

\begin{lemma}
\label{lem:fixed_index_coproducts}
Let $\mathsf I$ be a small set and let
$f_i:\mathcal X_i\to\mathcal Y_i$, for $i\in\mathsf I$, be maps of
simplicial sets. The map
\[
\coprod_{i\in\mathsf I}f_i:
\coprod_{i\in\mathsf I}\mathcal X_i
\longrightarrow\coprod_{i\in\mathsf I}\mathcal Y_i
\]
is a weak equivalence in $\sset_{(\infty,1)}$ if and only if every
$f_i$ is one. Moreover, if every $f_i$ is a fibration in
$\sset_{(\infty,1)}$, respectively a trivial fibration of simplicial
sets, then so is $\coprod_{i\in\mathsf I}f_i$.
\end{lemma}
\begin{proof}
The first assertion follows from \cite[Tags~06SZ and~06T1,
Remarks~4.8.8.7 and~4.8.8.9]{kerodon}: a map is a weak equivalence in
$\sset_{(\infty,1)}$ if and only if it is categorically
$n$-connective for every $n$, and a coproduct of maps is
categorically $n$-connective if and only if each of its summands is.
Suppose now that every $f_i$ is a fibration in
$\sset_{(\infty,1)}$, and consider a lifting problem as below left,
where $j$ is a trivial cofibration in $\sset_{(\infty,1)}$.
\[
\begin{tikzcd}[ampersand replacement=\&, sep=large]
\mathcal A\arrow[r]\arrow[d,hook,"j"',"\sim"]
\&\coprod_i\mathcal X_i\arrow[d,"{\coprod_i f_i}"]\\
\mathcal D\arrow[r,"b"']\arrow[ur,dashed]
\&\coprod_i\mathcal Y_i
\end{tikzcd}
\qquad \qquad
\begin{tikzcd}[ampersand replacement=\&, sep=large]
\mathcal A_i\arrow[r]\arrow[d,hook,"j_i"',"\sim"]
\&\mathcal X_i\arrow[d,two heads,"f_i"]\\
\mathcal D_i\arrow[r]\arrow[ur,dashed]\&\mathcal Y_i.
\end{tikzcd}
\]
The composite $\mathcal D\to\coprod_i\mathcal Y_i\to
\mathsf I^{\mathrm{disc}}$ decomposes the left-hand square as a
coproduct of squares as above right, with
$\mathcal D=\coprod_i\mathcal D_i$,
$\mathcal A=\coprod_i\mathcal A_i$, and $j=\coprod_i j_i$.
By the first assertion, each monomorphism $j_i$ is a weak
equivalence in $\sset_{(\infty,1)}$, and so the right-hand squares
admit lifts, which assemble into a lift in the left-hand square.
Hence $\coprod_if_i$ is a fibration in $\sset_{(\infty,1)}$. The
statement for trivial fibrations follows by combining the statement
for fibrations with the first assertion.
\end{proof}

\begin{proposition}
\label{prop:reedy_quasicategorical_ambient}
For $\mathbf M=\sset_{(\infty,1)}$, the model category
$\genreedyslice$ is left proper and combinatorial. Its cofibrations
are the levelwise monomorphisms, its weak equivalences are the levelwise
weak equivalences in $\sset_{(\infty,1)}$, and its fibrations are
the Reedy fibrations over $\mathbb B$. In particular, every object
is cofibrant.
\end{proposition}
\begin{proof}
This follows directly from
\cref{lem:presheaf_reedy_cofibrancy,prop:reedy_localization_existence}
applied to $\mathbf M=\sset_{(\infty,1)}$, and from the definition
of the slice model structure.
\end{proof}

\begin{notation}
\label{not:joyal_homotopy_generators}
We fix a set $\HGen$ of homotopy generators of
$\sset_{(\infty,1)}$, as in \cref{not:homotopy_generators}. The set
$\HGen$ will be used to define the localizing maps on both sides of
the comparison.
\end{notation}

\section{Symmetric monoidal \texorpdfstring{$(\infty,n)$}{(infinity,n)}-categories as left double fibrations}
\label{sec:internal_description}

In this section, we define symmetric monoidal $(\infty,n)$-categories
as special left double $(\infty,n)$-fibrations over $\bisimFin$, in
a model $\Minfn$ of $(\infty,n)$-categories. We introduce left
double $(\infty,n)$-fibrations in
\cref{subsec:internal_horizontal_segal} and the special ones in
\cref{subsec:internal_specialness}. In
\cref{subsec:general_reedy_arguments,subsec:general_fiber_comparison},
we give equivalent formulations of these two notions and show that
the specialness condition is invariant under weak equivalences. Assuming that $\Minfn$ is
left proper and combinatorial, we then construct in
\cref{subsec:general_left_double_model_structure,subsec:general_internal_special_model_structure}
model structures whose fibrant objects are the left double
$(\infty,n)$-fibrations and the special left double
$(\infty,n)$-fibrations, respectively. Finally, we specialize to
the case of quasi-categories in \cref{subsec:general_n_specialization}.

\begin{hypothesis}
\label{hyp:internal_standing_assumptions}
Throughout this section, we work with a model category $\Minfn$
satisfying the following assumptions.
\begin{enumerate}[label=\textup{(A\arabic*)}]
\item The model category $\Minfn$ is a model for the homotopy theory
of $(\infty,n)$-categories.
\item The simplicial object $\bisimFin$ of $\Minfn$ from
\cref{subsec:discrete_base} is Reedy fibrant and its Segal maps are
isomorphisms.
\end{enumerate}
\end{hypothesis}

Assumption \textup{(A1)} only fixes the terminology, and the proofs in
\crefrange{subsec:internal_horizontal_segal}{subsec:general_fiber_comparison}
only use \textup{(A2)}. In particular, we do not assume that $\Minfn$
is cartesian or simplicial, or that its objects are cofibrant, and we
assume that $\Minfn$ is left proper and combinatorial only in
\cref{subsec:general_left_double_model_structure,subsec:general_internal_special_model_structure}.
Fixed-shape fibers over $\bisimFin$ are always taken in the sense of
\cref{not:fixed_simplex_fibers}.

Assumption \textup{(A2)} holds under the hypotheses of
\cref{prop:discrete_base}, and it is a condition on the model
structure by \cref{rem:discrete_base_hypotheses}. For left Bousfield
localizations of injective model structures on simplicial presheaves,
it reduces to the fibrancy of the constant presheaves.

\begin{proposition}
\label{prop:presheaf_discrete_criterion}
Let $\mathrm C$ be a small category and let $\mathbf M$ be a left
Bousfield localization, at a small set $\mathsf L$ of maps, of the
injective model structure on simplicial presheaves
$\sset^{\mathrm C^{\mathrm{op}}}$. If the constant presheaf
$S^{\mathrm{disc}}$ is fibrant in $\mathbf M$ for every set $S$,
then $\mathbf M$ satisfies
\cref{hyp:internal_standing_assumptions}\textup{(A2)}.
\end{proposition}

\begin{proof}
The discrete realization functor has a left adjoint
$\pi_0:\sset^{\mathrm C^{\mathrm{op}}}\to\Set$, given by
\[
\pi_0\mathcal X=
\operatorname*{colim}_{c\in\mathrm C^{\mathrm{op}}}
\pi_0(\mathcal X(c)),
\qquad
\sset^{\mathrm C^{\mathrm{op}}}(\mathcal X,S^{\mathrm{disc}})
\cong\Set(\pi_0\mathcal X,S).
\]
Let $\mathcal A\hookrightarrow\mathcal B$ be a trivial cofibration in
the injective model structure. Since the map
$\mathcal A(c)\to\mathcal B(c)$ induces a bijection on connected
components for every object $c$ of $\mathrm C$, the function
$\pi_0\mathcal A\to\pi_0\mathcal B$ is a bijection. By adjunction, a
lifting problem of $\mathcal A\hookrightarrow\mathcal B$ against a
map $S^{\mathrm{disc}}\to T^{\mathrm{disc}}$ induced by a function
of sets $S\to T$ corresponds to a lifting problem of the bijection
$\pi_0\mathcal A\to\pi_0\mathcal B$ against the function $S\to T$,
which admits a solution. Hence every map
$S^{\mathrm{disc}}\to T^{\mathrm{disc}}$ induced by a function of
sets is an injective fibration. Since $S^{\mathrm{disc}}$ and
$T^{\mathrm{disc}}$ are fibrant in $\mathbf M$ by hypothesis, such
a map is also a fibration in $\mathbf M$ by
\cref{prop:localization_recollections}(1). As coproducts in a presheaf
category are disjoint and universal, the result follows from
\cref{prop:discrete_base}.
\end{proof}

\begin{rmk}
\label{rem:presheaf_discrete_components}
The hypothesis of \cref{prop:presheaf_discrete_criterion} can be
checked on connected components. Indeed, since
$\pi_0(\mathcal X\times\Delta[r])\cong\pi_0\mathcal X$, the
simplicial mapping space from a simplicial presheaf $\mathcal X$ to
$S^{\mathrm{disc}}$ is the discrete simplicial set
$\Set(\pi_0\mathcal X,S)^{\mathrm{disc}}$. As $S^{\mathrm{disc}}$
is injectively fibrant by the proof of
\cref{prop:presheaf_discrete_criterion}, it follows that the constant
presheaf $S^{\mathrm{disc}}$ is fibrant in $\mathbf M$ for every
set $S$ if and only if the functor $\pi_0$ sends every map in
$\mathsf L$ to a bijection.
\end{rmk}

\Cref{hyp:internal_standing_assumptions} holds for the following
models of $(\infty,n)$-categories, where the presheaf categories in
\textup{(2)} and \textup{(3)} are endowed with localizations of the
injective model structure.
\begin{enumerate}[label=\textup{(\arabic*)}]
\item The Joyal model structure $\sset_{(\infty,1)}$, for $n=1$;
see \cref{prop:quasicategorical_standing_assumptions}.
\item Rezk's model structure for complete Segal $\Theta_n$-spaces
on $\sset^{\Theta_n^{\mathrm{op}}}$ \cite{Rezk2010}. Given a set
$S$, the constant presheaf $S^{\mathrm{disc}}$ is the
$\Theta_n$-nerve of the discrete strict $n$-category on $S$,
which is rigid, as all its equivalences are identities. Hence
$S^{\mathrm{disc}}$ is fibrant by
\cite[Proposition~11.24]{Rezk2010}, and
\cref{prop:presheaf_discrete_criterion} applies.
\item The model structure for complete $n$-fold Segal spaces on
$\sset^{(\Delta^{\mathrm{op}})^n}$, obtained by localizing at the
Segal, essential constancy, and completeness maps described in
\cite[\S2.3.1]{calaque2019cobordisms}. Given a set $S$, the constant
presheaf $S^{\mathrm{disc}}$ is injectively fibrant by the proof of
\cref{prop:presheaf_discrete_criterion}, its Segal and essential
constancy maps are identities, and it is complete, as its morphisms in
every direction are identities. Hence $S^{\mathrm{disc}}$ is fibrant
in the localization, and \cref{prop:presheaf_discrete_criterion}
applies.
\end{enumerate}
These three model categories are moreover left proper and
combinatorial: for $\sset_{(\infty,1)}$, this fact is recalled in
\cref{subsec:prelim_simplicial_models}, and the model structures in
\textup{(2)} and \textup{(3)} are left Bousfield localizations of
injective model structures, which are left proper and combinatorial,
so that \cref{thm:left_bousfield_localization} applies.

\subsection{Double categories and left double fibrations}
\label{subsec:internal_horizontal_segal}

In this section, we introduce double $(\infty,n)$-categories and
left double $(\infty,n)$-fibrations in $\Minfn$. Given a simplicial
object $\mathbb X$ of $\Minfn$, we refer to $\mathbb X_0$ as the
$(\infty,n)$-category of objects and \emph{vertical} morphisms, and
to $\mathbb X_1$ as the $(\infty,n)$-category of \emph{horizontal}
morphisms and squares, with source and target maps
$d_1^h,d_0^h:\mathbb X_1\to\mathbb X_0$.

\begin{definition}[{cf. \cite[Definition~1.1.3]{Moser2026InftyN}}]
\label{def:horizontal_segal_object}
A \emph{double $(\infty,n)$-category in $\Minfn$} is a Reedy
fibrant simplicial object $\mathbb X$ of $\Minfn$ such that the
Segal map $\varphi_m:\mathbb X_m\to\mathbb X(I[m])$ from
\eqref{eq:horizontal_segal_map} is a weak equivalence in $\Minfn$,
for every $m\geq2$. A morphism of double $(\infty,n)$-categories
is a map of simplicial objects of $\Minfn$.
\end{definition}

The above definition is \cite[Definition~1.1.3]{Moser2026InftyN} with
$n$ replaced by $n+1$. Note that, for a Reedy fibrant simplicial
object $\mathbb X$, the Segal map $\varphi_m$ is a fibration by
\cref{lem:cellular_restriction}, and its target is a homotopy limit,
as it is an iterated pullback of fibrations between fibrant objects
\cite[Theorems~15.10.8\textup{(1)} and~19.9.1\textup{(2)}]{hirschhorn2003model}.

\begin{rmk}[{cf. \cite[Definitions~1.1.3 and~1.1.5]{Moser2026InftyN}}]
\label{rem:horizontal_completeness}
We do not impose any completeness condition in the horizontal
direction, while the completeness conditions of the model $\Minfn$
are part of the fibrancy in $\Minfn$. In particular, the horizontal
completeness condition and the condition on the space of objects used
in the cited definition of an $(\infty,n+1)$-category are not
imposed here.
\end{rmk}

Given a map $P:\mathbb E\to\bisimFin$ of simplicial objects of
$\Minfn$ and $m\geq0$, we denote by
\[
\lambda_m=(u_m^*,P_m):\mathbb E_m\longrightarrow
\mathbb E_0\times_{(\bisimFin)_0}(\bisimFin)_m
\]
the map induced by the inclusion
$u_m:\Delta[0]\hookrightarrow\Delta[m]$ of the initial vertex, where
the pullback is taken along the map
$u_m^*:(\bisimFin)_m\to(\bisimFin)_0$. Note that $u_1^*=d_1^h$.

The following definition is the analogue in $\Minfn$ of the left
fibrations of simplicial presheaves studied in
\cite[Definition~2.1 and Lemma~2.6]{Rasekh2021Yoneda}; in the case of
$\Theta$-spaces, the corresponding condition appears in
\cite[\S2.2 and Remark~2.2.5]{MoserRasekhRovelli2024LimitsII}.

\begin{definition}
\label{def:left-double-infinity-n-fibration}
A morphism $Q:\mathbb X\to\mathbb Y$ of double
$(\infty,n)$-categories is a \emph{left double
$(\infty,n)$-fibration} if it is a Reedy fibration and the map
\[
(d_1^h,Q_1):\mathbb X_1\longrightarrow
\mathbb X_0\times_{\mathbb Y_0}\mathbb Y_1
\]
is a weak equivalence in $\Minfn$, where the pullback is taken along
the source map $d_1^h:\mathbb Y_1\to\mathbb Y_0$. In particular,
under \cref{hyp:internal_standing_assumptions}, a map
$P:\mathbb E\to\bisimFin$ is a left double $(\infty,n)$-fibration
if $\mathbb E$ is a double $(\infty,n)$-category, $P$ is a Reedy
fibration, and the map $\lambda_1$ is a weak equivalence in
$\Minfn$.
\end{definition}

With this convention, left double $(\infty,n)$-fibrations are tested
at the initial vertex; see \cite{Rasekh2023Yoneda} for a discussion in
the case of Segal spaces. Note that, given a Reedy fibration
$P:\mathbb E\to\bisimFin$, the simplicial object $\mathbb E$ is
Reedy fibrant, as $\bisimFin$ is, and the map $\lambda_1$ is a
fibration by \cref{lem:cellular_restriction}. Hence the last condition
of \cref{def:left-double-infinity-n-fibration} says that $\lambda_1$
is a trivial fibration. Informally, every object $x$ of $\mathbb E$
over $a$ and every morphism $f:a\to b$ of $\Fin$ admit a lift to
a horizontal morphism of $\mathbb E$ with source $x$, and such a
lift is unique up to a contractible space of choices.

\subsection{Inert spans and symmetric monoidality}
\label{subsec:internal_specialness}

In this section, we define special left double
$(\infty,n)$-fibrations over $\bisimFin$. Let
$P:\mathbb E\to\bisimFin$ be a left double $(\infty,n)$-fibration.
The horizontal morphisms of $\mathbb E$ over the inert maps
$\rho_1,\ldots,\rho_k$ with a common source form a span from the
fiber $\mathbb E_0(\langle k\rangle)$ to the product
$\mathbb E_0(\langle1\rangle)^k$, and the specialness condition
requires the right leg of this span to be a weak equivalence.

Given a morphism $f:a\to b$ in $\Fin$, the face maps $d_1^h$ and
$d_0^h$ restrict to maps between fixed-shape fibers
\[
s_f^P:\mathbb E_1(f)\longrightarrow\mathbb E_0(a),
\qquad
t_f^P:\mathbb E_1(f)\longrightarrow\mathbb E_0(b),
\]
which form a span in $\Minfn$ from $\mathbb E_0(a)$ to
$\mathbb E_0(b)$. For $k\geq1$, we consider the \emph{inert span}
\begin{equation}
\label{eq:segality_double_span}
\mathbb E_0(\langle k\rangle)
\xleftarrow{\ \pi_k^P\ }
\mathbb E_1(\rho_1)
\times_{\mathbb E_0(\langle k\rangle)}\cdots
\times_{\mathbb E_0(\langle k\rangle)}\mathbb E_1(\rho_k)
\xrightarrow{\ t_k^P\ }
\mathbb E_0(\langle1\rangle)^k,
\end{equation}
where the pullbacks are taken along the maps $s_{\rho_i}^P$, the map
$\pi_k^P$ is the canonical projection, and
$t_k^P=(t_{\rho_i}^P\pr_i)_{i=1}^k$. For $k=0$, we set
\begin{equation}
\label{eq:internal_nullary_span}
\pi_0^P=\id_{\mathbb E_0(\langle0\rangle)},
\qquad
t_0^P:\mathbb E_0(\langle0\rangle)\longrightarrow[0].
\end{equation}

\begin{definition}[{cf. \cite{Seg74};
\cite[Definition~2.0.0.7 and \S2.4.2]{LurieHA}}]
\label{def:special_left_double}
A left double $(\infty,n)$-fibration
$P:\mathbb E\to\bisimFin$ is \emph{special} if the map
\[
t_k^P=(t_{\rho_i}^P\pr_i)_{i=1}^k:
\mathbb E_1(\rho_1)
\times_{\mathbb E_0(\langle k\rangle)}\cdots
\times_{\mathbb E_0(\langle k\rangle)}\mathbb E_1(\rho_k)
\longrightarrow\mathbb E_0(\langle1\rangle)^k
\]
is a weak equivalence in $\Minfn$ for every $k\geq1$, and so is
the map $t_0^P$ from \eqref{eq:internal_nullary_span}.
\end{definition}

The terminology comes from the special $\Gamma$-spaces of
\cite{Seg74}. The maps $t_k^P$ play the role of the maps induced by
the inert maps in \cite[Definition~2.0.0.7]{LurieHA}; compare also the
monoid objects of \cite[\S2.4.2]{LurieHA}. The source and target of
$t_k^P$ can be expressed in terms of the inert fan of
\cref{def:inert_fan} as follows.

\begin{lemma}
\label{lem:internal_fan_identification}
Let $P:\mathbb E\to\bisimFin$ be a left double
$(\infty,n)$-fibration and let $k\geq1$. There are isomorphisms
\[
\mathbb E(\mathcal J_k)_{\beta_k}\cong
\mathbb E_1(\rho_1)
\times_{\mathbb E_0(\langle k\rangle)}\cdots
\times_{\mathbb E_0(\langle k\rangle)}\mathbb E_1(\rho_k),
\qquad
\mathbb E(\mathcal V_k)_{\beta_k|\mathcal V_k}
\cong\mathbb E_0(\langle1\rangle)^k,
\]
natural in $P$, under which the restriction maps along the
inclusions $\Delta\{c\}\hookrightarrow\mathcal J_k$ and
$\mathcal V_k\hookrightarrow\mathcal J_k$ correspond to the maps
$\pi_k^P$ and $t_k^P$, respectively. For $k=0$, the fixed-shape
fibers at $\mathcal J_0$ and $\mathcal V_0$ are
$\mathbb E_0(\langle0\rangle)$ and $[0]$, and the restriction maps
are the maps of \eqref{eq:internal_nullary_span}.
\end{lemma}

\begin{proof}
The inert fan $\mathcal J_k$ is the pushout of $k$ copies of
$\Delta[1]$ along their initial vertices, and the functor
$\mathbb E(-)$ sends this pushout to a pullback. Passing to the
fibers at $\beta_k$, we get the first isomorphism, under which the
restriction along $\Delta\{c\}\hookrightarrow\mathcal J_k$ is the
projection $\pi_k^P$. Similarly, the simplicial set $\mathcal V_k$
is the coproduct of the vertices $v_1,\ldots,v_k$, which gives the
second isomorphism, and the restriction along
$\mathcal V_k\hookrightarrow\mathcal J_k$ is then the map $t_k^P$.
Both isomorphisms are natural in $P$ by construction. For $k=0$,
we have
$\mathbb E(\Delta[0])_{\langle0\rangle}=\mathbb E_0(\langle0\rangle)$
and $\mathbb E(\varnothing)\cong[0]$.
\end{proof}

\begin{definition}
\label{def:symmetric_monoidal_internal}
A \emph{symmetric monoidal $(\infty,n)$-category in $\Minfn$}
is a special left double $(\infty,n)$-fibration
$P:\mathbb E\to\bisimFin$. Its \emph{underlying
$(\infty,n)$-category} is the fibrant object
$\mathbb E_0(\langle1\rangle)$ of $\Minfn$. A \emph{functor of
symmetric monoidal $(\infty,n)$-categories} is a map of simplicial
objects over $\bisimFin$.
\end{definition}

\begin{rmk}
\label{rem:internal_monoidal_functors}
For coCartesian fibrations over $\simFin$, a symmetric monoidal
functor is required to preserve coCartesian edges; see
\cite[Definition~2.1.3.7]{LurieHA}. No such condition is needed in
\cref{def:symmetric_monoidal_internal}, as every horizontal morphism
of a left double $(\infty,n)$-fibration plays the role of a
coCartesian edge. Note that, if not every object of $\Minfn$ is
cofibrant, then the morphisms in the homotopy category between two
symmetric monoidal $(\infty,n)$-categories are represented by maps
out of a cofibrant replacement in the model structure of
\cref{thm:general_internal_special_model_structure}. In the models
listed after \cref{prop:presheaf_discrete_criterion}, every object is
cofibrant. For $n=1$, the comparison with Lurie's definition is
given in
\cref{thm:staircase_preserves_reflects_specialness,cor:strings_specialness}.
\end{rmk}

\begin{rmk}
\label{rem:higher_n_comparison}
For $n=1$, we prove in \cref{sec:comparisons} that
\cref{def:symmetric_monoidal_internal} is equivalent to Lurie's
definition. For $n>1$, we do not prove a comparison in this paper,
but we expect the following. Left fibrations over a base that is
discrete in the $(\infty,n)$-direction should straighten to diagrams
of $(\infty,n)$-categories, as is the case for Segal copresheaves
over an $(n+1)$-category by \cite[Corollary~4.33]{Nuiten2023}, and
for the Grothendieck constructions of \cite{Rasekh2021Yoneda} and
\cite{moser2023inftyncategoricalstraighteningunstraighteningconstruction}
in their respective models. Under such an equivalence, a left double
$(\infty,n)$-fibration over $\bisimFin$ should correspond to a
diagram of $(\infty,n)$-categories indexed by $\Fin$, the span
$(\pi_k^P,t_k^P)$ to the map into the $k$-fold product induced by
the inert maps, and the specialness condition to the $\Gamma$-object
condition of \cite[Definition~3.3]{calaque2019cobordisms}. We have not
verified that the cited equivalences apply to the model structures of
\cref{subsec:general_left_double_model_structure,subsec:general_internal_special_model_structure}.
\end{rmk}

\begin{ex}
\label{ex:discrete_commutative_monoid}
Let $n=1$ and $\Minfn=\sset_{(\infty,1)}$, and let
$F:\Fin\to\mathbf{Set}$ be a functor with category of elements
$\int F$. The projection $\int F\to\Fin$ induces a map
\[
P:\mathbb E=N(\textstyle\int F)^h\longrightarrow\bisimFin.
\]
Since the levels and the matching objects of $\mathbb E$ and
$\bisimFin$ are discrete simplicial sets, the map $P$ is a Reedy
fibration, and the Segal maps of $\mathbb E$ are isomorphisms, as
$N(\int F)$ is a nerve. Moreover, a string of composable maps in
$\Fin$ together with an element of $F$ at its source lifts
uniquely to a string of composable morphisms in $\int F$, so that
each $\lambda_m$ is an isomorphism, and $P$ is a left double
$(\infty,1)$-fibration. For $k\geq0$, we have
\[
\mathbb E_0(\langle k\rangle)=F\langle k\rangle,
\qquad
\mathbb E(\mathcal J_k)_{\beta_k}\cong F\langle k\rangle,
\qquad
t_k^P=(F(\rho_1),\ldots,F(\rho_k)),
\]
where sets are regarded as discrete simplicial sets. Hence $P$ is
special if and only if the map
\[
\theta_k=(F(\rho_i))_{i=1}^k:
F\langle k\rangle\longrightarrow F\langle1\rangle^k
\]
is a bijection for every $k\geq0$, i.e., if and only if $F$ is a
special $\Gamma$-set in the sense of \cite{Seg74}. In this case,
the set $M=F\langle1\rangle$ is a commutative monoid, whose unit is
the image under $F(\mu_0)$ of the unique element of
$F\langle0\rangle$ and whose multiplication is given by
\[
x\cdot y=F(\mu_2)\theta_2^{-1}(x,y).
\]
More generally, for a map
$f:\langle k\rangle\to\langle\ell\rangle$ in $\Fin$, we have
\[
\bigl(\theta_{\ell}F(f)\theta_k^{-1}(x_1,\ldots,x_k)\bigr)_j
=\prod_{i\in f^{-1}(j)}x_i,
\qquad 1\leq j\leq\ell,
\]
and conversely this formula defines a special functor for every
commutative monoid \cite{Seg74}. On the other hand, a constant functor
at a set with at least two elements gives a left double
$(\infty,1)$-fibration that is not special.
\end{ex}

We now show that the left leg of the span
\eqref{eq:segality_double_span} is a trivial fibration and that its
right leg is a fibration.

\begin{proposition}
\label{prop:cleavage_free_segal}
Let $P:\mathbb E\to\bisimFin$ be a left double
$(\infty,n)$-fibration. For every morphism $f:a\to b$ in $\Fin$,
the map $s_f^P$ is a trivial fibration and the map
\[
(s_f^P,t_f^P):\mathbb E_1(f)\longrightarrow
\mathbb E_0(a)\times\mathbb E_0(b)
\]
is a fibration. Moreover, for every $k\geq0$, the three objects of
the span
\[
\mathbb E_0(\langle k\rangle)
\xleftarrow{\ \pi_k^P\ }
\mathbb E(\mathcal J_k)_{\beta_k}
\xrightarrow{\ t_k^P\ }
\mathbb E_0(\langle1\rangle)^{\,k}
\]
are fibrant in $\Minfn$, the map $\pi_k^P$ is a trivial fibration,
and the map $t_k^P$ is a fibration.
\end{proposition}

\begin{proof}
Since $P$ is a Reedy fibration, the map
$P_0:\mathbb E_0\to(\bisimFin)_0$ is a fibration, and so its fiber
$\mathbb E_0(a)$ is fibrant for every object $a$ of $\Fin$. Let
$f:a\to b$ be a morphism in $\Fin$, and set
\[
A_{\mathrm{src}}=\mathbb E_0\times_{(\bisimFin)_0}(\bisimFin)_1,
\qquad
A_{\mathrm{end}}=(\mathbb E_0\times\mathbb E_0)
\times_{(\bisimFin)_0\times(\bisimFin)_0}(\bisimFin)_1,
\]
where the pullbacks are taken along the source map $d_1^h$ and the
map $(d_1^h,d_0^h)$, respectively. The relative matching map
$r_1:\mathbb E_1\to A_{\mathrm{end}}$ is a fibration, and the map
$\lambda_1:\mathbb E_1\to A_{\mathrm{src}}$ is a trivial fibration,
as observed after \cref{def:left-double-infinity-n-fibration}. Writing
$j_a,j_b,j_f$ for the inclusions of the fibers
$\mathbb E_0(a),\mathbb E_0(b),\mathbb E_1(f)$ into
$\mathbb E_0,\mathbb E_0,\mathbb E_1$, respectively, we have the
following pullback squares in $\Minfn$.
\[
\begin{tikzcd}[ampersand replacement=\&, sep=large]
\mathbb E_1(f)\arrow[r,"j_f"]
  \arrow[d,two heads,"s_f^P"',"\sim"]
  \arrow[dr,phantom,"\lrcorner",very near start]
\&\mathbb E_1\arrow[d,two heads,"\lambda_1","\sim"']\\
\mathbb E_0(a)\arrow[r,"{(j_a,f!)}"']\&A_{\mathrm{src}}
\end{tikzcd}
\qquad \qquad
\begin{tikzcd}[ampersand replacement=\&, sep=large]
\mathbb E_1(f)\arrow[r,"j_f"]
  \arrow[d,two heads,"{(s_f^P,t_f^P)}"']
  \arrow[dr,phantom,"\lrcorner",very near start]
\&\mathbb E_1\arrow[d,two heads,"r_1"]\\
\mathbb E_0(a)\times\mathbb E_0(b)
  \arrow[r,"{(j_a\times j_b,f!)}"']\&A_{\mathrm{end}}
\end{tikzcd}
\]
Indeed, a map from an object $C$ of $\Minfn$ to the pullback of
the left-hand cospan consists of maps $x:C\to\mathbb E_1$ and
$a':C\to\mathbb E_0(a)$ such that $P_1x=f\,!_C$ and
$j_aa'=d_1^hx$. Since $f$ has source $a$, the equation
$P_1x=f\,!_C$ implies that $d_1^hx$ factors uniquely through
$\mathbb E_0(a)$, and so such a pair $(x,a')$ amounts to a map
$C\to\mathbb E_1(f)$. The same argument applies to the right-hand
square. Hence $s_f^P$ is a trivial fibration and $(s_f^P,t_f^P)$
is a fibration. For $k\geq1$, by \cref{lem:common_source_fibrations}
applied to $A=\mathbb E_0(\langle k\rangle)$,
$H_i=\mathbb E_1(\rho_i)$, and $D_i=\mathbb E_0(\langle1\rangle)$,
and by \cref{lem:internal_fan_identification}, the map $\pi_k^P$ is
a trivial fibration, the map $t_k^P$ is a fibration, and the object
$\mathbb E(\mathcal J_k)_{\beta_k}$ is fibrant. Moreover, the object
$\mathbb E_0(\langle1\rangle)^{\,k}$ is fibrant as a product of
fibrant objects. Finally, the case $k=0$ follows from the case
$k=0$ of \cref{lem:common_source_fibrations} and
\eqref{eq:internal_nullary_span}.
\end{proof}

As a consequence, a left double $(\infty,n)$-fibration $P$ is
special if and only if the map $t_k^P$ is a trivial fibration for
every $k\geq0$.

\begin{rmk}
\label{rem:segality_small_arities}
For $k=0$, the specialness condition requires the map
$\mathbb E_0(\langle0\rangle)\to[0]$ to be a weak equivalence, but
not necessarily an isomorphism. For $k=1$, the specialness condition
is automatic. Indeed, we have $\rho_1=\id_{\langle1\rangle}$, and the
degeneracy map restricts to a map
$s_0^h:\mathbb E_0(\langle1\rangle)\to
\mathbb E_1(\id_{\langle1\rangle})$, which is a common section of the
maps $s_{\rho_1}^P$ and $t_{\rho_1}^P$. Since $s_{\rho_1}^P$ is
a trivial fibration by \cref{prop:cleavage_free_segal}, the map
$s_0^h$ is a weak equivalence, and hence so is the map
$t_1^P=t_{\rho_1}^P$, by the two-out-of-three property.
\end{rmk}

\begin{ex}
\label{ex:empty_left_double}
The specialness condition for $k=0$ does not follow from the
conditions for $k\geq1$. Indeed, for
$\Minfn=\sset_{(\infty,1)}$, the unique map
$P:\varnothing\to\bisimFin$ from the empty simplicial object is a
left double $(\infty,1)$-fibration: it is a Reedy fibration, as its
relative matching maps are maps between discrete simplicial sets, and
its Segal maps and the map $\lambda_1$ are isomorphisms between
empty simplicial sets. For every $k\geq1$, the map $t_k^P$ is an
isomorphism between empty simplicial sets, but the map $t_0^P$ is
the map $\varnothing\to\Delta[0]$, which is not a weak equivalence
in $\sset_{(\infty,1)}$.
\end{ex}

\subsection{Source conditions and induced operations}

\label{subsec:general_reedy_arguments}

In this section, we show that a Reedy fibration
$P:\mathbb E\to\bisimFin$ is a left double $(\infty,n)$-fibration
if and only if all the maps $\lambda_m$ are weak equivalences, and
we describe the operations induced by a left double
$(\infty,n)$-fibration on its fibers. In the case of Segal spaces,
the first result is \cite[Lemma~3.29]{Rasekh2023Yoneda}.

\begin{proposition}
\label{prop:minimal_axioms}
Let $P:\mathbb E\to\bisimFin$ be a Reedy fibration. Then
$\mathbb E$ is Reedy fibrant, each $\mathbb E_m$ is fibrant in
$\Minfn$, and each map $\lambda_m$ is a fibration. Moreover, the
following conditions are equivalent.
\begin{enumerate}[label=\textup{(\roman*)}]
\item The Segal map $\varphi_m$ of $\mathbb E$ is a weak
      equivalence for every $m\geq2$ and $\lambda_1$ is a weak
      equivalence, i.e., the map $P$ is a left double
      $(\infty,n)$-fibration.
\item The map $\lambda_m$ is a weak equivalence for every
      $m\geq 0$.
\item The map $\lambda_m$ is a trivial fibration for every
      $m\geq 0$.
\end{enumerate}
If these conditions hold, then the Segal maps of $\mathbb E$ are
trivial fibrations.
\end{proposition}

\begin{proof}
Since $\bisimFin$ is Reedy fibrant, so is $\mathbb E$. By
\cref{lem:cellular_restriction} applied to the inclusions
$\varnothing\hookrightarrow\Delta[m]$ and
$I[m]\hookrightarrow\Delta[m]$, the objects $\mathbb E_m$ are
fibrant and the Segal maps $\varphi_m$ of $\mathbb E$ are
fibrations, and by the same lemma applied to $P$ and the inclusion
$\Delta\{0\}\hookrightarrow\Delta[m]$, the maps $\lambda_m$ are
fibrations. In particular, conditions \textup{(ii)} and \textup{(iii)}
are equivalent. To compare \textup{(i)} and \textup{(ii)}, suppose
that $\lambda_1$ is a trivial fibration, which is the case under
either condition, and let $m\geq1$. For $0\leq j\leq m$, let
$I[j]\subseteq\Delta[m]$ denote the spine of the face
$\Delta\{0,\ldots,j\}$, and set
\[
C_j=
\mathbb E(I[j])
\times_{\bisimFin(I[j])}(\bisimFin)_m,
\]
so that $C_0=\mathbb E_0\times_{(\bisimFin)_0}(\bisimFin)_m$. For
$1\leq j\leq m$, let
$\mathrm{edge}_j:(\bisimFin)_m\to(\bisimFin)_1$ denote the
restriction to the edge $\Delta\{j-1,j\}$, and let
$a_{j-1}:C_{j-1}\to\mathbb E_0$ denote the evaluation at the vertex
$j-1$. Applying the functors $\mathbb E(-)$ and $\bisimFin(-)$ to
the pushout
$I[j]=I[j-1]\amalg_{\Delta\{j-1\}}\Delta\{j-1,j\}$ of simplicial
sets, we get a pullback square in $\Minfn$
\[
\begin{tikzcd}[ampersand replacement=\&, row sep=large, column sep=huge]
C_j\arrow[r]\arrow[d,two heads,"\sim"']
  \arrow[dr,phantom,"\lrcorner",very near start]
\&\mathbb E_1\arrow[d,two heads,"\lambda_1","\sim"']\\
C_{j-1}\arrow[r,"{(a_{j-1},\mathrm{edge}_j\pr_2)}"']
\&\mathbb E_0\times_{(\bisimFin)_0}(\bisimFin)_1,
\end{tikzcd}
\]
and so the projection $C_j\to C_{j-1}$ is a trivial fibration. Since
the Segal map $(\bisimFin)_m\to\bisimFin(I[m])$ is an isomorphism, we
have $C_m\cong\mathbb E(I[m])$, and the composite
$C_m\to C_{m-1}\to\cdots\to C_0$ gives a trivial fibration
\[
\gamma_m:\mathbb E(I[m])\longrightarrow
\mathbb E_0\times_{(\bisimFin)_0}(\bisimFin)_m
\]
such that $\lambda_m=\gamma_m\varphi_m$. By the two-out-of-three
property, the map $\lambda_m$ is a weak equivalence if and only if
the Segal map $\varphi_m$ is one. As $\lambda_0$ is an
isomorphism, this shows that \textup{(i)} and \textup{(ii)} are
equivalent. Finally, if these conditions hold, then the Segal maps of
$\mathbb E$ are trivial fibrations, as they are fibrations.
\end{proof}

\begin{corollary}
\label{cor:fixed_simplex_source}
Let $P:\mathbb E\to\bisimFin$ be a left double
$(\infty,n)$-fibration. For every $m\geq0$ and every $m$-simplex
$\sigma$ of $\simFin$, the restriction map along the inclusion of
the initial vertex
\[
\lambda_m^\sigma=u_m^*:
\mathbb E_m(\sigma)\longrightarrow\mathbb E_0(\sigma(0))
\]
is a trivial fibration between fibrant objects of $\Minfn$.
\end{corollary}

\begin{proof}
The fibers $\mathbb E_m(\sigma)$ and $\mathbb E_0(\sigma(0))$ are
fibrant by \cref{lem:fixed_shape_fiber_fibrant}. Moreover, the map
$\lambda_m$ is a trivial fibration by \cref{prop:minimal_axioms}, and
hence so is its pullback $\lambda_m^\sigma$ displayed in the
following diagram.
\[
\begin{tikzcd}[ampersand replacement=\&, sep=large]
\mathbb E_m(\sigma)\arrow[r,"\pr_1"]
  \arrow[d,two heads,"\lambda_m^\sigma"',"\sim"]
  \arrow[dr,phantom,"\lrcorner",very near start]
\&\mathbb E_m\arrow[d,two heads,"\lambda_m","\sim"']\\
\mathbb E_0(\sigma(0))
  \arrow[r,"{(\pr_1,\sigma !)}"']
\&\mathbb E_0\times_{(\bisimFin)_0}(\bisimFin)_m.
\end{tikzcd}
\]
\end{proof}

\begin{rmk}[Induced operations]
\label{rem:internal_transport_operations}
A left double $(\infty,n)$-fibration $P:\mathbb E\to\bisimFin$
induces a functor
\[
\Fin\longrightarrow\Ho(\Minfn),
\qquad
a\longmapsto\mathbb E_0(a),
\qquad
f\longmapsto T_f=[t_f^P][s_f^P]^{-1},
\]
where $[-]$ denotes the image in the homotopy category, and
$[s_f^P]$ is invertible by \cref{prop:cleavage_free_segal}. Indeed,
the assignment $f\mapsto T_f$ preserves identities, as the
degeneracy map $s_0^h$ is a common section of $s_{\id_a}^P$ and
$t_{\id_a}^P$. Given composable maps
$a\xrightarrow{f}b\xrightarrow{g}c$ in $\Fin$ and the corresponding
$2$-simplex $\sigma=(f,g)$ of $\simFin$, the restriction
$v:\mathbb E_2(\sigma)\to\mathbb E_0(a)$ to the initial vertex is a
trivial fibration by \cref{cor:fixed_simplex_source}, and the face
maps $d_2^h:\mathbb E_2(\sigma)\to\mathbb E_1(f)$,
$d_1^h:\mathbb E_2(\sigma)\to\mathbb E_1(gf)$, and
$d_0^h:\mathbb E_2(\sigma)\to\mathbb E_1(g)$ satisfy
\[
s_f^Pd_2^h=v=s_{gf}^Pd_1^h,
\qquad
t_f^Pd_2^h=s_g^Pd_0^h,
\qquad
t_g^Pd_0^h=t_{gf}^Pd_1^h.
\]
Hence
\[
T_gT_f[v]
=[t_g^P][s_g^P]^{-1}[t_f^Pd_2^h]
=[t_g^Pd_0^h]
=[t_{gf}^Pd_1^h]
=T_{gf}[v],
\]
and so $T_gT_f=T_{gf}$, as $[v]$ is an isomorphism. Moreover,
since finite products of fibrant objects are products in
$\Ho(\Minfn)$ by \cref{prop:model_category_recollections}(5) and
$s_{\rho_i}^P\pr_i=\pi_k^P$, the map
\[
S_k=[t_k^P][\pi_k^P]^{-1}:
\mathbb E_0(\langle k\rangle)
\longrightarrow\mathbb E_0(\langle1\rangle)^k
\]
is the map $(T_{\rho_i})_{i=1}^k$, for every $k\geq1$. Hence $P$
is special if and only if $S_k$ is an isomorphism for every
$k\geq0$, where $S_0$ denotes the map to the terminal object. In
this case, the active maps $\mu_k$ induce the $k$-ary operations
\[
T_{\mu_k}S_k^{-1}:\mathbb E_0(\langle1\rangle)^k
\longrightarrow\mathbb E_0(\langle1\rangle)
\]
in $\Ho(\Minfn)$, with the unit given by the case $k=0$. The
higher coherences of these operations are encoded by the simplicial
object $\mathbb E$.
\end{rmk}

\subsection{Homotopy invariance of fibers and specialness}

\label{subsec:general_fiber_comparison}

In this section, we show that the specialness condition is invariant
under weak equivalences between left double $(\infty,n)$-fibrations.
For this, we first show that these weak equivalences are detected in
degree zero; see \cite[Theorem~4.62]{Rasekh2021Yoneda} for a similar
result for localized left fibrations of simplicial presheaves.

\begin{lemma}
\label{lem:left_double_degree_zero}
Let $F:\mathbb E\to\mathbb E'$ be a map over $\bisimFin$
between left double $(\infty,n)$-fibrations. Then $F$ is a
levelwise weak equivalence if and only if $F_0$ is a weak
equivalence in $\Minfn$.
\end{lemma}
\begin{proof}
Suppose that $F_0$ is a weak equivalence, and let
$P:\mathbb E\to\bisimFin$ and $P':\mathbb E'\to\bisimFin$ denote
the structure maps. We write
$b_m:(\bisimFin)_m\to(\bisimFin)_0$ for the restriction to the
initial vertex. Since $P_0$ and $P'_0$ are fibrations, the map
$F_0$ is a weak equivalence between fibrant objects of the slice
$\Minfn/(\bisimFin)_0$, and so its pullback $b_m^*F_0$ along
$b_m$ is a weak equivalence by
\cref{prop:model_category_recollections}(6). In the commutative square
\[
\begin{tikzcd}[ampersand replacement=\&, sep=large]
\mathbb E_m\arrow[r,"F_m"]
  \arrow[d,two heads,"\lambda_m"',"\sim"]
\&\mathbb E'_m\arrow[d,two heads,"\lambda'_m","\sim"']\\
\mathbb E_0\times_{(\bisimFin)_0}(\bisimFin)_m
  \arrow[r,"b_m^*F_0"',"\sim"]
\&\mathbb E'_0\times_{(\bisimFin)_0}(\bisimFin)_m
\end{tikzcd}
\]
the vertical maps are trivial fibrations by
\cref{prop:minimal_axioms}. Hence $F_m$ is a weak equivalence for
every $m\geq0$, by the two-out-of-three property.
\end{proof}

As a consequence, such a map $F$ induces weak equivalences between
the fibers considered in \cref{subsec:internal_specialness}.
\begin{lemma}
\label{lem:general_fiber_comparison}
Let $P:\mathbb E\to\bisimFin$ and
$P':\mathbb E'\to\bisimFin$ be left double
$(\infty,n)$-fibrations, and let
$F:\mathbb E\to\mathbb E'$ be a map over $\bisimFin$ such that
$F_0$ is a weak equivalence. Then, for every $m\geq0$ and every
$m$-simplex $\sigma$ of $\simFin$, the induced map
\[
F_m(\sigma):\mathbb E_m(\sigma)\longrightarrow
\mathbb E'_m(\sigma)
\]
is a weak equivalence in $\Minfn$, and so is the induced map
$\mathbb E(\mathcal J_k)_{\beta_k}\to
\mathbb E'(\mathcal J_k)_{\beta_k}$, for every $k\geq0$. In
particular, the left double $(\infty,n)$-fibration $P$ is special
if and only if $P'$ is special.
\end{lemma}

\begin{proof}
By \cref{lem:left_double_degree_zero}, the map $F_m$ is a weak
equivalence for every $m\geq0$. Since $P_m$ and $P'_m$ are
fibrations by \cref{lem:cellular_restriction} applied to the inclusion
$\varnothing\hookrightarrow\Delta[m]$, the map $F_m$ is a weak
equivalence between fibrant objects of the slice
$\Minfn/(\bisimFin)_m$. Hence its pullback $F_m(\sigma)$ along the
point $\sigma:[0]\to(\bisimFin)_m$ is a weak equivalence by
\cref{prop:model_category_recollections}(6). Next, for $k\geq0$, let
$r_k:\mathbb E(\mathcal J_k)_{\beta_k}\to
\mathbb E'(\mathcal J_k)_{\beta_k}$ denote the map induced by $F$.
In the commutative square
\[
\begin{tikzcd}[ampersand replacement=\&, sep=large]
\mathbb E(\mathcal J_k)_{\beta_k}\arrow[r,"r_k"]
  \arrow[d,two heads,"\pi_k^P"',"\sim"]
\&\mathbb E'(\mathcal J_k)_{\beta_k}\arrow[d,two heads,"\pi_k^{P'}","\sim"']\\
\mathbb E_0(\langle k\rangle)
  \arrow[r,"F_0(\langle k\rangle)"',"\sim"]
\&\mathbb E'_0(\langle k\rangle)
\end{tikzcd}
\]
the vertical maps are trivial fibrations by
\cref{prop:cleavage_free_segal}, and the map $F_0(\langle k\rangle)$
is a weak equivalence by the first part of the proof. Hence $r_k$ is
a weak equivalence, by the two-out-of-three property. Moreover, in the
commutative square
\[
\begin{tikzcd}[ampersand replacement=\&, sep=large]
\mathbb E(\mathcal J_k)_{\beta_k}\arrow[r,"r_k","\sim"']\arrow[d,two heads,"t_k^P"']
\&\mathbb E'(\mathcal J_k)_{\beta_k}\arrow[d,two heads,"t_k^{P'}"]\\
\mathbb E_0(\langle1\rangle)^k
  \arrow[r,"F_0(\langle1\rangle)^{\,k}"',"\sim"]
\&\mathbb E'_0(\langle1\rangle)^k
\end{tikzcd}
\]
the map $F_0(\langle1\rangle)^{\,k}$ is a weak equivalence by
\cref{prop:model_category_recollections}(5), as a finite product of
weak equivalences between fibrant objects. Hence the map $t_k^P$ is
a weak equivalence if and only if $t_k^{P'}$ is one, by the
two-out-of-three property, and so $P$ is special if and only if
$P'$ is special.
\end{proof}

\subsection{The general model structure for left double fibrations}
\label{subsec:general_left_double_model_structure}

In this section, we construct a model structure on the category of
simplicial objects of $\Minfn$ over $\bisimFin$ whose fibrant
objects are the left double $(\infty,n)$-fibrations. We obtain this
model structure as the left Bousfield localization of the slice Reedy
model structure at a set of maps that detect when the maps
$\lambda_m$ are weak equivalences.

\begin{hypothesis}
\label{hyp:internal_localization_assumptions}
In this section and in
\cref{subsec:general_internal_special_model_structure}, we assume,
in addition to \cref{hyp:internal_standing_assumptions}, that
$\Minfn$ is left proper and combinatorial.
\end{hypothesis}

By \cref{prop:reedy_localization_existence}, the slice Reedy model
structure
$\bigl((\Minfn)^{\Delta^{\mathrm{op}}}_{/\bisimFin}\bigr)$
is then left proper and combinatorial. Model structures of the same
kind are constructed in
\cite[Theorems~4.17 and~5.15]{Rasekh2021Yoneda} for simplicial
presheaves, and the model structure for double left fibrations of
$\Theta$-spaces is recalled in
\cite[\S2.2 and Remark~2.2.5]{MoserRasekhRovelli2024LimitsII}; we do
not compare these model structures with the one constructed below.

As before, we write $b_m:(\bisimFin)_m\to(\bisimFin)_0$ for the
restriction to the initial vertex. To define the localizing maps, we
use the left adjoint of the evaluation functor in degree $m$.

\begin{lemma}
\label{lem:slice_evaluation_source_adjunction}
Let $m\geq0$. Then the following statements hold.
\begin{enumerate}[label=\textup{(\arabic*)}]
\item There is a Quillen adjunction
\[
\adjunction{\Minfn/(\bisimFin)_m}
{\bigl((\Minfn)^{\Delta^{\mathrm{op}}}_{/\bisimFin}\bigr)}
{\mathbb F_m}{\operatorname{ev}_m}
\]
where, for an object $a:A\to(\bisimFin)_m$, the simplicial object
$\mathbb F_mA$ is given by
\[
(\mathbb F_mA)_\ell=
\coprod_{\alpha:[\ell]\to[m]}A,
\]
with structure map $\alpha_h^*a$ on the summand indexed by
$\alpha$.
\item The functor
$(b_m)_!:\Minfn/(\bisimFin)_m\to\Minfn/(\bisimFin)_0$ given by
postcomposition with $b_m$ is left Quillen, and the inclusion of the
initial vertex induces a natural map
\[
j^{\mathrm{glob}}_{m,A}:
\mathbb F_0((b_m)_!A)\longrightarrow\mathbb F_mA,
\]
which is a cofibration between cofibrant objects if $A$ is cofibrant
in $\Minfn/(\bisimFin)_m$.
\item Given a cofibrant object $A$ of $\Minfn/(\bisimFin)_m$ and a
Reedy fibration $P:\mathbb E\to\bisimFin$, the map induced by
$j^{\mathrm{glob}}_{m,A}$ on derived mapping spaces into
$\mathbb E$ is identified with the map
\[
\RMap(A,\mathbb E_m)
\xrightarrow{\RMap(A,\lambda_m)}
\RMap\bigl(A,b_m^*\mathbb E_0\bigr),
\]
where the derived mapping spaces are taken in
$\Minfn/(\bisimFin)_m$ and
$b_m^*\mathbb E_0=
\mathbb E_0\times_{(\bisimFin)_0}(\bisimFin)_m$.
\end{enumerate}
\end{lemma}

\begin{proof}
A map $\mathbb F_mA\to\mathbb E$ over $\bisimFin$ is determined by
its component at the summand indexed by $\id_{[m]}$, which is a map
$A\to\mathbb E_m$ over $(\bisimFin)_m$, and this gives the
adjunction in \textup{(1)}. The underlying simplicial object of
$\mathbb F_mA$ is $\Delta[m]\boxtimes A$, and so, by the
computation of the latching objects in the proof of
\cref{lem:general_shape_quillen}, the relative latching map in degree
$\ell$ of the map $\mathbb F_mA\to\mathbb F_mA'$ induced by a map
$A\to A'$ in $\Minfn/(\bisimFin)_m$ is the map
\[
\left(\coprod_{\mathsf D_\ell\Delta[m]}A'\right)
\amalg\left(\coprod_{\mathsf N_\ell\Delta[m]}A\right)
\longrightarrow
\left(\coprod_{\mathsf D_\ell\Delta[m]}A'\right)
\amalg\left(\coprod_{\mathsf N_\ell\Delta[m]}A'\right),
\]
which is the identity on the summands indexed by the degenerate
simplices and the map $A\to A'$ on the summands indexed by the
nondegenerate simplices. As in the proof of
\cref{lem:general_shape_quillen}, it follows that $\mathbb F_m$
preserves cofibrations and trivial cofibrations, since the
cofibrations and weak equivalences of both slices are created in the
underlying model categories. For the same reason, the functor
$(b_m)_!$ is left Quillen. The underlying map of
$j^{\mathrm{glob}}_{m,A}$ is the map $u_m\boxtimes A$, and the
underlying object of a cofibrant object $A$ of
$\Minfn/(\bisimFin)_m$ is cofibrant in $\Minfn$. Hence
$j^{\mathrm{glob}}_{m,A}$ is a cofibration between cofibrant objects
by \cref{lem:shape_inclusion_cofibrant}, as the cofibrations of the
slice are created in the Reedy model structure, which proves
\textup{(2)}.
For \textup{(3)}, let $A^\bullet\to\operatorname{const}(A)$ be a
cosimplicial resolution of $A$ in $\Minfn/(\bisimFin)_m$. Since
$A$ is cofibrant, the left Quillen functors $\mathbb F_m$ and
$\mathbb F_0(b_m)_!$ send $A^\bullet$ to cosimplicial resolutions
of $\mathbb F_mA$ and $\mathbb F_0((b_m)_!A)$. In each cosimplicial
degree, the adjunctions identify the maps from these resolutions to
$\mathbb E$ with the maps from $A^\bullet$ to $\mathbb E_m$ and
to $b_m^*\mathbb E_0$, respectively, in such a way that
precomposition with the initial vertex map corresponds to
postcomposition with $\lambda_m$. Finally, the objects
$\mathbb E_m$ and $b_m^*\mathbb E_0$ are fibrant in
$\Minfn/(\bisimFin)_m$, as $P_m$ is a fibration by
\cref{lem:cellular_restriction} and the structure map of
$b_m^*\mathbb E_0$ is a pullback of $P_0$. Hence these
identifications compute the derived mapping spaces in \textup{(3)}
and the map between them.
\end{proof}

For each $m\geq1$, we choose a set $\mathsf H_m$ of cofibrant
homotopy generators of the combinatorial model category
$\Minfn/(\bisimFin)_m$, which exists by
\cref{not:homotopy_generators}, and we define
\begin{equation}
\label{eq:general_source_localizers}
\mathsf J_{\mathrm{glob}}=
\left\{j^{\mathrm{glob}}_{m,A}\ \middle|\
 m\geq1,\ A\in\mathsf H_m\right\},
\qquad
\genldblfib=
L_{\mathsf J_{\mathrm{glob}}}
\bigl((\Minfn)^{\Delta^{\mathrm{op}}}_{/\bisimFin}\bigr).
\end{equation}
Note that the maps in $\mathsf J_{\mathrm{glob}}$ are cofibrations
between cofibrant objects by
\cref{lem:slice_evaluation_source_adjunction}, and that the case
$m=0$ is not needed, as $\lambda_0$ is an isomorphism.

\begin{theorem}
\label{thm:general_left_double_model_structure}
Under \cref{hyp:internal_localization_assumptions}, the model
structure $\genldblfib$ exists, and it is left proper and
combinatorial. Its cofibrations are the Reedy cofibrations and its
fibrant objects are the left double $(\infty,n)$-fibrations
$P:\mathbb E\to\bisimFin$. Moreover, a map between fibrant objects
is a weak equivalence if and only if it is a levelwise weak
equivalence, if and only if it is a weak equivalence in degree zero,
and it is a fibration if and only if it is a Reedy fibration.
\end{theorem}

\begin{proof}
The existence of $\genldblfib$, its left properness and
combinatoriality, and the description of its cofibrations follow from
\cref{prop:reedy_localization_existence,thm:left_bousfield_localization}.
Let $P:\mathbb E\to\bisimFin$ be a Reedy fibration and let
$m\geq1$. By \cref{lem:slice_evaluation_source_adjunction}, the map
$P$ is local with respect to the maps $j^{\mathrm{glob}}_{m,A}$,
for $A\in\mathsf H_m$, if and only if the map
$\RMap(A,\lambda_m)$ of derived mapping spaces in
$\Minfn/(\bisimFin)_m$ is a weak equivalence for every
$A\in\mathsf H_m$. Since the source and target of $\lambda_m$ are
fibrant in $\Minfn/(\bisimFin)_m$, this is the case if and only if
$\lambda_m$ is a weak equivalence, by
\cref{lem:homotopy_generators_detect}. Hence, by
\cref{prop:minimal_axioms}, the fibrant objects of $\genldblfib$ are
the left double $(\infty,n)$-fibrations. The descriptions of the
weak equivalences and of the fibrations between fibrant objects follow
from \cref{prop:localization_recollections}(1) and
\cref{lem:left_double_degree_zero}.
\end{proof}

\subsection{The general model structure for special left double fibrations}
\label{subsec:general_internal_special_model_structure}

In this section, we localize the model structure $\genldblfib$ to
obtain a model structure whose fibrant objects are the special left
double $(\infty,n)$-fibrations. We assume
\cref{hyp:internal_localization_assumptions}, and we choose a set
$\mathsf H$ of cofibrant homotopy generators of $\Minfn$. Recall
from \cref{lem:internal_fan_identification} that the map $t_k^P$ is
the restriction along the inclusion
$\mathcal V_k\hookrightarrow\mathcal J_k$. We define
\begin{equation}
\label{eq:general_internal_special_localizers}
\mathsf S_{\mathrm{sp}}^{\Minfn}=
\left\{
 s^{\Minfn}_{k,A}:
 (\mathcal V_k\boxtimes A)_{\beta_k|\mathcal V_k}
 \longrightarrow
 (\mathcal J_k\boxtimes A)_{\beta_k}
 \ \middle|\ k\geq0,\ A\in\mathsf H
\right\},
\end{equation}
where the map $s^{\Minfn}_{k,A}$ is induced by the inclusion
$\mathcal V_k\hookrightarrow\mathcal J_k$, and
\[
\genldblfibsp=
L_{\mathsf S_{\mathrm{sp}}^{\Minfn}}\genldblfib.
\]
Note that the maps in $\mathsf S_{\mathrm{sp}}^{\Minfn}$ are
cofibrations between cofibrant objects by
\cref{lem:shape_inclusion_cofibrant}, and that, for $k=0$, they are
the maps
$\varnothing\to(\Delta[0]\boxtimes A)_{\langle0\rangle}$.

\begin{theorem}
\label{thm:general_internal_special_model_structure}
Under \cref{hyp:internal_localization_assumptions}, the model
structure $\genldblfibsp$ exists, and it is left proper and
combinatorial. Its cofibrations are the Reedy cofibrations and its
fibrant objects are the special left double $(\infty,n)$-fibrations
over $\bisimFin$. Moreover, a map between fibrant objects is a weak
equivalence if and only if it is a levelwise weak equivalence, if and
only if it is a weak equivalence in degree zero, and it is a fibration
if and only if it is a Reedy fibration.
\end{theorem}

\begin{proof}
The existence of $\genldblfibsp$, its left properness and
combinatoriality, and the description of its cofibrations follow from
\cref{thm:left_bousfield_localization} applied to the model structure
of \cref{thm:general_left_double_model_structure}. Let
$P:\mathbb E\to\bisimFin$ be a left double $(\infty,n)$-fibration
and let $k\geq0$. By
\cref{lem:shape_locality,lem:internal_fan_identification}, the map
$P$ is local with respect to the maps $s^{\Minfn}_{k,A}$, for
$A\in\mathsf H$, if and only if the map $t_k^P$ is a weak
equivalence. Hence the fibrant objects of $\genldblfibsp$ are the
special left double $(\infty,n)$-fibrations. The descriptions of the
weak equivalences and of the fibrations between fibrant objects follow
from \cref{prop:localization_recollections}(1) and
\cref{thm:general_left_double_model_structure}.
\end{proof}

\begin{rmk}
\label{rem:independence_generators}
The model structures $\genldblfib$ and $\genldblfibsp$ do not
depend on the choice of the sets of homotopy generators. Indeed, a
left Bousfield localization is determined by its cofibrations and its
fibrant objects, and
\cref{thm:general_left_double_model_structure,thm:general_internal_special_model_structure}
describe the fibrant objects without reference to these sets. For the
same reason, the model structure $\genldblfibsp$ is the left
Bousfield localization of the slice Reedy model structure at the set
$\mathsf J_{\mathrm{glob}}\cup\mathsf S_{\mathrm{sp}}^{\Minfn}$.
\end{rmk}

\subsection{The quasi-categorical presentation}
\label{subsec:general_n_specialization}

In this section, we specialize to the case
$\Minfn=\sset_{(\infty,1)}$, with the conventions of
\cref{subsec:prelim_simplicial_models}, and we give presentations of
the model structures of
\cref{thm:general_left_double_model_structure,thm:general_internal_special_model_structure}
that are adapted to the comparison of \cref{sec:comparisons}.

\begin{proposition}
\label{prop:quasicategorical_standing_assumptions}
The model category $\sset_{(\infty,1)}$ satisfies
\cref{hyp:internal_standing_assumptions,hyp:internal_localization_assumptions}.
\end{proposition}

\begin{proof}
The model category $\sset_{(\infty,1)}$ is a model for
$(\infty,1)$-categories, and it is left proper and combinatorial, as
recalled in \cref{subsec:prelim_simplicial_models}. Moreover,
coproducts of simplicial sets are disjoint and universal, and every
map between discrete simplicial sets is a fibration in
$\sset_{(\infty,1)}$. Hence assumption \textup{(A2)} follows from
\cref{prop:discrete_base}.
\end{proof}

We denote by
\[
\ldblfib
\qquad\text{and}\qquad
\ldblfibsp
\]
the model structures of
\cref{thm:general_left_double_model_structure,thm:general_internal_special_model_structure}
for $\Minfn=\sset_{(\infty,1)}$. We now show that these model
structures can be obtained by localizing at maps defined over the
simplices of $\simFin$ and over the inert fans.

\subsubsection{The source localization over individual simplices}
\label{subsubsec:left_double_model_structure}

Recall the set $\HGen$ of homotopy generators of
$\sset_{(\infty,1)}$ from \cref{not:joyal_homotopy_generators}. For
$m\geq1$, an $m$-simplex $\sigma$ of $\simFin$, and
$\mathcal A\in\HGen$, we denote by
\[
j_{m,\sigma,\mathcal A}:
(\Delta[0]\boxtimes\mathcal A)_{\sigma(0)}
\longrightarrow(\Delta[m]\boxtimes\mathcal A)_\sigma
\]
the map induced by the inclusion
$u_m:\Delta[0]\hookrightarrow\Delta[m]$ of the initial vertex, and
we denote by $\mathsf J$ the set of all such maps.

\begin{theorem}
\label{thm:left_double_model_structure}
The model structure $\ldblfib$ is the left Bousfield localization
\[
\ldblfib=L_{\mathsf J}\reedyslice.
\]
It is left proper and combinatorial, its cofibrations are the
levelwise monomorphisms, and its fibrant objects are the left double
$(\infty,1)$-fibrations $P:\mathbb E\to\bisimFin$. In particular,
every object is cofibrant. Moreover, for a map $F$ between fibrant
objects, the following conditions are equivalent.
\begin{enumerate}[label=\textup{(\roman*)}]
\item The map $F$ is a weak equivalence in $\ldblfib$.
\item The map $F$ is a levelwise weak equivalence in
      $\sset_{(\infty,1)}$.
\item The map $F_0$ is a weak equivalence in
      $\sset_{(\infty,1)}$.
\end{enumerate}
Finally, a map between fibrant objects is a fibration in $\ldblfib$
if and only if it is a Reedy fibration.
\end{theorem}

\begin{proof}
The left Bousfield localization $L_{\mathsf J}\reedyslice$ exists by
\cref{prop:reedy_quasicategorical_ambient,thm:left_bousfield_localization}.
Let $P:\mathbb E\to\bisimFin$ be a Reedy fibration. By
\cref{lem:shape_locality}, the map $P$ is local with respect to the
maps $j_{m,\sigma,\mathcal A}$, for $\mathcal A\in\HGen$, if and
only if the map
\[
\lambda_m^\sigma:
\mathbb E_m(\sigma)\longrightarrow\mathbb E_0(\sigma(0))
\]
is a weak equivalence in $\sset_{(\infty,1)}$. Since
$\lambda_m=\coprod_{\sigma\in(\simFin)_m}\lambda_m^\sigma$, the maps
$\lambda_m^\sigma$ are weak equivalences for all
$\sigma\in(\simFin)_m$ if and only if $\lambda_m$ is a weak
equivalence, by \cref{lem:fixed_index_coproducts}. Hence the model
structures $L_{\mathsf J}\reedyslice$ and $\ldblfib$ have the same
fibrant objects by \cref{thm:general_left_double_model_structure}. As
both are left Bousfield localizations of $\reedyslice$, they
coincide. The remaining assertions follow from
\cref{thm:general_left_double_model_structure,prop:reedy_quasicategorical_ambient}.
\end{proof}

\subsubsection{The special localization over inert fans}
\label{subsubsec:internal_special_model_structure}

We denote by
\begin{equation}
\label{eq:internal_special_localizers}
\mathsf S_{\mathrm{sp}}=
\left\{
 s_{k,\mathcal A}:
 (\mathcal V_k\boxtimes\mathcal A)_{\beta_k|\mathcal V_k}
 \longrightarrow
 (\mathcal J_k\boxtimes\mathcal A)_{\beta_k}
 \ \middle|\ k\geq0,\ \mathcal A\in\HGen
\right\}
\end{equation}
the set \eqref{eq:general_internal_special_localizers} for
$\Minfn=\sset_{(\infty,1)}$ and $\mathsf H=\HGen$. A localization
of the same kind is carried out in \cite[\S5]{Sharma2021} for diagrams
of quasi-categories indexed by finite pointed sets.

\begin{theorem}
\label{thm:internal_special_model_structure}
The model structure $\ldblfibsp$ is the left Bousfield localization
\[
\ldblfibsp=L_{\mathsf S_{\mathrm{sp}}}\ldblfib.
\]
It is left proper and combinatorial, its cofibrations are the
levelwise monomorphisms, and its fibrant objects are the special left
double $(\infty,1)$-fibrations. In particular, every object is
cofibrant. Moreover, a map between fibrant objects is a weak
equivalence if and only if it is a levelwise weak equivalence in
$\sset_{(\infty,1)}$, if and only if it is a weak equivalence in
$\sset_{(\infty,1)}$ in degree zero, and it is a fibration if and
only if it is a Reedy fibration.
\end{theorem}

\begin{proof}
The displayed equality holds by definition of $\genldblfibsp$ for
the choice $\mathsf H=\HGen$, and the model structure
$\genldblfibsp$ does not depend on this choice by
\cref{rem:independence_generators}. The remaining assertions follow
from
\cref{thm:general_internal_special_model_structure,prop:reedy_quasicategorical_ambient}.
\end{proof}

\section{Symmetric monoidal \texorpdfstring{$(\infty,1)$}{(infinity,1)}-categories as coCartesian fibrations}
\label{sec:enriched_description}

In this section, we recall Lurie's definition of symmetric monoidal
$(\infty,1)$-categories as coCartesian fibrations over $\simFin$,
in a form parallel to \cref{sec:internal_description}. We recall
coCartesian fibrations and the coCartesian model structure in
\cref{subsec:cocartesian_definitions,subsec:marked_mapping_objects},
and we study the restriction maps between marked mapping objects in
\cref{subsec:enriched_mapping_fibrations,subsec:enriched_marked_simplices}.
In \cref{subsec:enriched_specialness}, we express the Segal condition
on a coCartesian fibration over $\simFin$ in terms of spans analogous
to those of \cref{subsec:internal_specialness}, and we construct in
\cref{subsec:enriched_special_model_structure} a model structure whose
fibrant objects are the coCartesian fibrations satisfying this
condition. Finally, in \cref{subsec:enriched_boundary_tests}, we
characterize these coCartesian fibrations by a lifting property.

Throughout this section, simplicial sets are endowed with the model
structure $\sset_{(\infty,1)}$ recalled in
\cref{subsec:prelim_simplicial_models}; see also
\cite{Joyal2008Theory}.

\subsection{Marked simplicial sets and coCartesian fibrations}
\label{subsec:cocartesian_definitions}

In this section, we recall the main definitions about marked
simplicial sets and coCartesian fibrations from
\cite[\S2.4 and \S3.1]{lurie2009higher}, dualized from Cartesian to
coCartesian fibrations. A \emph{marked simplicial set} is a pair
$(\mathcal K,\mathsf E)$ of a simplicial set $\mathcal K$ and a set
$\mathsf E\subseteq\mathcal K_1$ of \emph{marked} edges containing
the degenerate edges, and a map of marked simplicial sets is a map of
simplicial sets that preserves marked edges
\cite[Definition~3.1.0.1]{lurie2009higher}. We denote by $\sset^+$
the category of marked simplicial sets, in which an edge of a product
is marked if and only if both of its components are marked. Given a
simplicial set $\mathcal K$, we write $\mathcal K^\flat$ for the
marked simplicial set in which only the degenerate edges are marked,
and $\mathcal K^\sharp$ for the one in which all edges are marked.

Let $p:\mathcal E\to\mathcal K$ be an inner fibration. An edge
$e:x\to y$ of $\mathcal E$ is \emph{$p$-coCartesian} if the map
\[
\mathcal E_{e/}\longrightarrow
\mathcal E_{x/}\times_{\mathcal K_{p(x)/}}\mathcal K_{p(e)/}
\]
is a trivial fibration of simplicial sets, where the slices are those
of \cite[Remark~1.2.9.5]{lurie2009higher}; equivalently, if
$e^{\mathrm{op}}$ is $p^{\mathrm{op}}$-Cartesian in the sense of
\cite[Definition~2.4.1.1]{lurie2009higher}. The inner fibration $p$
is a \emph{coCartesian fibration} if, for every vertex $x$ of
$\mathcal E$ and every edge $f:p(x)\to b$ of $\mathcal K$, there
is a $p$-coCartesian edge $e:x\to y$ of $\mathcal E$ such that
$p(e)=f$ \cite[Definition~2.4.2.1]{lurie2009higher}. We will use the
following criterion to recognize coCartesian edges.

\begin{lemma}[{dual of \cite[Remark~2.4.1.4]{lurie2009higher}}]
\label{lem:cocartesian_edge_horn_criterion}
Let $p:\mathcal E\to\mathcal K$ be an inner fibration. An edge $e$
of $\mathcal E$ is $p$-coCartesian if and only if, for every
$q\geq2$, every commutative square of simplicial sets
\[
\begin{tikzcd}[ampersand replacement=\&, sep=large]
\Lambda^0[q]\arrow[r,"\eta"]\arrow[d,hook]
\&\mathcal E\arrow[d,"p"]\\
\Delta[q]\arrow[r,"\sigma"']\arrow[ur,dashed,"\widetilde\eta"]
\&\mathcal K
\end{tikzcd}
\]
with $\eta|_{\Delta\{0,1\}}=e$ admits a lift $\widetilde\eta$ as
indicated.
\end{lemma}

Given a vertex $a$ of $\mathcal K$, the fiber
$\mathcal E_a\coloneqq\mathcal E\times_{\mathcal K}\Delta[0]$ of an
inner fibration $p:\mathcal E\to\mathcal K$ is a quasi-category. If
$\mathcal K$ is a quasi-category, then so is $\mathcal E$, and the
degenerate edges of $\mathcal E$ are $p$-coCartesian by the dual of
\cite[Corollary~2.4.1.6]{lurie2009higher}.

\begin{notation}[{dual of \cite[Definition~3.1.1.8]{lurie2009higher}}]
\label{not:natural_marking}
Let $p:\mathcal E\to\mathcal K$ be a coCartesian fibration over a
quasi-category $\mathcal K$. We denote by $\mathcal E^\natural$ the
marked simplicial set obtained from $\mathcal E$ by marking the
$p$-coCartesian edges, and we refer to this marking as the
\emph{natural marking}. Given a map $\beta:\mathcal A\to\mathcal K$
of simplicial sets, we write $\mathcal A^\sharp_\beta$ and
$\mathcal A^\flat_\beta$ for the marked simplicial sets
$\mathcal A^\sharp$ and $\mathcal A^\flat$ over
$\mathcal K^\sharp$ with structure map $\beta$.
\end{notation}

In the rest of \cref{sec:enriched_description}, we take
$\mathcal K=\simFin$, and
$p:\mathcal E\to\simFin$ denotes a coCartesian fibration, with fiber
$\mathcal E_a$ over an object $a$ of $\Fin$.

\subsection{Marked mapping objects and the coCartesian model structure}
\label{subsec:marked_mapping_objects}

In this section, we recall the mapping objects of marked simplicial
sets over $\simFin^\sharp$ and the coCartesian model structure from
\cite[\S3.1.3]{lurie2009higher}. Given marked simplicial sets
$\mathcal A$ and $\mathcal X$ over $\simFin^\sharp$, we denote by
$\Map^\flat(\mathcal A,\mathcal X)$ and
$\Map^\sharp(\mathcal A,\mathcal X)$ the simplicial sets given by
\begin{equation}
\label{eq:marked_flat_mapping}
\begin{aligned}
\Map^\flat(\mathcal A,\mathcal X)_r
&=(\sset^+)_{/\simFin^\sharp}
  \bigl(\Delta[r]^\flat\times\mathcal A,\mathcal X\bigr),\\
\Map^\sharp(\mathcal A,\mathcal X)_r
&=(\sset^+)_{/\simFin^\sharp}
  \bigl(\Delta[r]^\sharp\times\mathcal A,\mathcal X\bigr),  
\end{aligned}
\end{equation}
where the products $\Delta[r]^\flat\times\mathcal A$ and
$\Delta[r]^\sharp\times\mathcal A$ are regarded as objects over
$\simFin^\sharp$ through the projection to $\mathcal A$.

Similarly, given maps of simplicial sets
$\beta:\mathcal A\to\simFin$ and $g:\mathcal Y\to\simFin$, we
denote by $\Map(\mathcal A,\mathcal Y)$ the simplicial set of maps
over $\simFin$, given by
\begin{equation}
\label{eq:unmarked_relative_mapping}
\begin{aligned}
\Map(\mathcal A,\mathcal Y)_r
 &=\sset_{/\simFin}(\Delta[r]\times\mathcal A,\mathcal Y),\\
\Map(\mathcal A,\mathcal Y)
 &=\mathcal Y^{\mathcal A}
   \times_{(\simFin)^{\mathcal A}}\Delta[0],
\end{aligned}
\end{equation}
where $\mathcal Y^{\mathcal A}$ denotes the internal hom of
simplicial sets and the map $\Delta[0]\to(\simFin)^{\mathcal A}$ is
the vertex $\beta$. For instance, an $r$-simplex of
$\Map^\flat(\mathcal K^\sharp_\beta,\mathcal E^\natural)$
is a map $F:\Delta[r]\times\mathcal K\to\mathcal E$
such that $pF=\beta\pr_2$ and such that, for each vertex $i$ of
$\Delta[r]$, the restriction of $F$ to $\{i\}\times\mathcal K$
sends every edge to a $p$-coCartesian edge. We refer to the edges of
$\{i\}\times\mathcal K$ as the \emph{horizontal} edges of
$\Delta[r]\times\mathcal K$.

\begin{theorem}[{dual of \cite[Proposition~3.1.3.7, Remark~3.1.3.9,
and Proposition~3.1.4.1]{lurie2009higher}}]
\label{thm:cocartesian_model_structure}
The category $(\sset^+)_{/\simFin^\sharp}$ admits a left proper
combinatorial model structure, denoted by $\cocart$, whose
cofibrations are the monomorphisms and whose fibrant objects are the
naturally marked coCartesian fibrations over $\simFin$. In
particular, every object is cofibrant.
\end{theorem}

\begin{proposition}[{dual of \cite[Proposition~3.1.3.3 and the definition
following it, and Proposition~3.3.1.5]{lurie2009higher}}]
\label{prop:cocartesian_weak_equivalences}
A map $g:\mathcal A\to\mathcal D$ of marked simplicial sets over
$\simFin^\sharp$ is a weak equivalence in $\cocart$ if and only
if, for every coCartesian fibration $\mathcal Z\to\simFin$, the
induced map
\[
g^*:\Map^\flat(\mathcal D,\mathcal Z^\natural)
\longrightarrow\Map^\flat(\mathcal A,\mathcal Z^\natural)
\]
is a weak equivalence in $\sset_{(\infty,1)}$. Moreover, a map
between naturally marked coCartesian fibrations over $\simFin$ is a
weak equivalence in $\cocart$ if and only if its underlying map of
simplicial sets is a weak equivalence in $\sset_{(\infty,1)}$, if
and only if it induces a weak equivalence in $\sset_{(\infty,1)}$ on
every fiber.
\end{proposition}

The model category $\cocart$ is enriched over
$\sset_{(\infty,1)}$, with tensor $\mathcal U^\flat\times\mathcal X$
and mapping object $\Map^\flat(\mathcal X,\mathcal Y)$, by the dual
of \cite[Remark~3.1.4.5]{lurie2009higher}, and over
$\sset_{(\infty,0)}$, with tensor
$\mathcal U^\sharp\times\mathcal X$ and mapping object
$\Map^\sharp(\mathcal X,\mathcal Y)$, by the dual of
\cite[Corollary~3.1.4.4]{lurie2009higher}. In particular, the
simplicial set $\Map^\flat(\mathcal X,\mathcal Y)$ is a
quasi-category if $\mathcal Y$ is fibrant. Recall that the
\emph{core} $\mathcal U^\simeq$ of a quasi-category $\mathcal U$
is its maximal Kan subcomplex, consisting of the simplices of
$\mathcal U$ all of whose edges are equivalences. Given an object
$\mathcal A$ and a fibrant object $\mathcal X$ of $\cocart$, we
have
\begin{equation}
\label{eq:marked_mapping_core}
\Map^\sharp(\mathcal A,\mathcal X)
=\bigl(\Map^\flat(\mathcal A,\mathcal X)\bigr)^\simeq
\simeq\RMap(\mathcal A,\mathcal X),
\end{equation}
where the equality is the dual of
\cite[Remark~3.1.3.1]{lurie2009higher}, and the weak equivalence
follows from the enrichment over $\sset_{(\infty,0)}$, as
$\mathcal A$ is cofibrant and $\mathcal X$ is fibrant.

\subsection{Restriction maps and weak equivalences}
\label{subsec:enriched_mapping_fibrations}

In this section, we show that the restriction maps between the mapping
objects $\Map^\flat$ along monomorphisms are fibrations in
$\sset_{(\infty,1)}$, and we deduce a criterion for locality in
$\cocart$, analogous to \cref{lem:shape_locality}.

\begin{lemma}[{dual of \cite[Corollary~3.1.4.3]{lurie2009higher}}]
\label{lem:cocartesian_joyal_pushout_product}
Let $i:\mathcal A\hookrightarrow\mathcal D$ be a cofibration
in $\cocart$ and let
$j:\mathcal K\hookrightarrow\mathcal L$ be a cofibration
in $\sset_{(\infty,1)}$. Then the pushout-product
\[
j^\flat\mathbin{\square}i:
(\mathcal K^\flat\times\mathcal D)
\amalg_{\mathcal K^\flat\times\mathcal A}
(\mathcal L^\flat\times\mathcal A)
\longrightarrow\mathcal L^\flat\times\mathcal D
\]
is a cofibration in $\cocart$, which is trivial if $i$ is a weak
equivalence in $\cocart$ or $j$ is a weak equivalence in
$\sset_{(\infty,1)}$.
\end{lemma}

\begin{proof}
Since taking opposites preserves products, pushouts, and
monomorphisms, the dual of
\cite[Corollary~3.1.4.3]{lurie2009higher} says that the
pushout-product of a cofibration in the coCartesian model structure
over $\Delta[0]^\sharp$ with a cofibration in $\cocart$ is a
cofibration in the coCartesian model structure over
$\Delta[0]^\sharp\times\simFin^\sharp\cong\simFin^\sharp$, which is
trivial if either of the two cofibrations is trivial. The map
$j^\flat$ is a monomorphism of marked simplicial sets over
$\Delta[0]^\sharp$. Moreover, if $j$ is a weak equivalence in
$\sset_{(\infty,1)}$, then $j^\flat$ is a weak equivalence in the
Cartesian model structure over $\Delta[0]^\sharp$ by
\cite[Proposition~3.1.5.3]{lurie2009higher}, and hence in the
coCartesian model structure over $\Delta[0]^\sharp$ by
\cite[Remark~3.1.4.6]{lurie2009higher}. The result follows by applying
the above statement to $j^\flat$ and $i$.
\end{proof}

\begin{proposition}
\label{prop:marked_mapping_fibration}
Let $\mathcal X$ be a fibrant object of $\cocart$ and let
$i:\mathcal A\hookrightarrow\mathcal D$ be a monomorphism of marked
simplicial sets over $\simFin^\sharp$. Then the restriction map
\[
i^*:\Map^\flat(\mathcal D,\mathcal X)
\longrightarrow\Map^\flat(\mathcal A,\mathcal X)
\]
is a fibration in $\sset_{(\infty,1)}$, which is trivial if $i$ is
a weak equivalence in $\cocart$. In particular, the simplicial set
$\Map^\flat(\mathcal A,\mathcal X)$ is a quasi-category.
\end{proposition}

\begin{proof}
Let $j:\mathcal K\hookrightarrow\mathcal L$ be a trivial
cofibration in $\sset_{(\infty,1)}$. By adjunction, a lifting
problem in simplicial sets
\[
\begin{tikzcd}[ampersand replacement=\&, sep=large]
\mathcal K\arrow[r]\arrow[d,hook,"j"',"\sim"]
\&\Map^\flat(\mathcal D,\mathcal X)\arrow[d,"i^*"]\\
\mathcal L\arrow[r]\arrow[ur,dashed]
\&\Map^\flat(\mathcal A,\mathcal X)
\end{tikzcd}
\]
corresponds to a lifting problem in marked simplicial sets
\[
\begin{tikzcd}[ampersand replacement=\&, sep=large]
(\mathcal K^\flat\times\mathcal D)
\amalg_{\mathcal K^\flat\times\mathcal A}
(\mathcal L^\flat\times\mathcal A)
\arrow[r]\arrow[d,hook,"j^\flat\mathbin{\square}i"',"\sim"]
\&\mathcal X\arrow[d,two heads]\\
\mathcal L^\flat\times\mathcal D
\arrow[r]\arrow[ur,dashed]\&\simFin^\sharp
\end{tikzcd}
\]
which admits a lift, as the map $j^\flat\mathbin{\square}i$ is a
trivial cofibration in $\cocart$ by
\cref{lem:cocartesian_joyal_pushout_product} and $\mathcal X$ is
fibrant. Hence $i^*$ is a fibration in $\sset_{(\infty,1)}$. If
$i$ is a weak equivalence in $\cocart$, then the map
$j^\flat\mathbin{\square}i$ is a trivial cofibration for every
monomorphism $j$, and so $i^*$ is a trivial fibration by the same
argument. The last assertion follows by taking $i$ to be the
inclusion $\varnothing\hookrightarrow\mathcal A$.
\end{proof}

\begin{lemma}
\label{lem:marked_shape_locality}
Let $i:\mathcal A\hookrightarrow\mathcal D$ be a monomorphism of
marked simplicial sets over $\simFin^\sharp$ and let $\mathcal X$
be a fibrant object of $\cocart$. For every simplicial set
$\mathcal U$, there is a weak equivalence in $\sset_{(\infty,0)}$
\[
\RMap(\mathcal D\times\mathcal U^\flat,\mathcal X)
\simeq
\RMap\bigl(\mathcal U,\Map^\flat(\mathcal D,\mathcal X)\bigr),
\]
natural in $\mathcal D$. In particular, the object $\mathcal X$
is local with respect to the set of maps
\[
\left\{i\times\mathcal U^\flat
  \ \middle|\ \mathcal U\in\HGen\right\}
\]
if and only if the restriction map
$i^*:\Map^\flat(\mathcal D,\mathcal X)
\to\Map^\flat(\mathcal A,\mathcal X)$
is a weak equivalence in $\sset_{(\infty,1)}$, where $\HGen$ is
the set of homotopy generators from
\cref{not:joyal_homotopy_generators}.
\end{lemma}

\begin{proof}
By \cref{lem:cocartesian_joyal_pushout_product} applied to the
inclusion $\varnothing\hookrightarrow\mathcal D$, we have a Quillen
adjunction
\[
\begin{tikzcd}[ampersand replacement=\&,column sep=huge]
\sset_{(\infty,1)}
\arrow[r,shift left=1.2ex,"\mathcal D\times(-)^\flat"]
\&\cocart
\arrow[l,shift left=1.2ex,"{\Map^\flat(\mathcal D,-)}"]
\end{tikzcd}
\]
and similarly for $\mathcal A$. Let
$\mathcal U^\bullet\to\operatorname{const}(\mathcal U)$ be a
cosimplicial resolution of $\mathcal U$ in $\sset_{(\infty,1)}$.
As in the proof of \cref{lem:derived_shape_adjunction}, the left
adjoint sends $\mathcal U^\bullet$ to a cosimplicial resolution of
$\mathcal D\times\mathcal U^\flat$, and the adjunction in each
cosimplicial degree gives the desired weak equivalence; see also
\cite[Proposition~17.4.16\textup{(1)}]{hirschhorn2003model}. Using the
same resolution for $\mathcal A$ and $\mathcal D$, we get a
commutative square
\[
\begin{tikzcd}[ampersand replacement=\&, sep=large]
\RMap(\mathcal D\times\mathcal U^\flat,\mathcal X)
\arrow[r,"\simeq"]\arrow[d,"(i\times\mathcal U^\flat)^*"']
\&\RMap\bigl(\mathcal U,\Map^\flat(\mathcal D,\mathcal X)\bigr)
\arrow[d,"{\RMap(\mathcal U,i^*)}"]\\
\RMap(\mathcal A\times\mathcal U^\flat,\mathcal X)
\arrow[r,"\simeq"']
\&\RMap\bigl(\mathcal U,\Map^\flat(\mathcal A,\mathcal X)\bigr).
\end{tikzcd}
\]
Since $\Map^\flat(\mathcal D,\mathcal X)$ and
$\Map^\flat(\mathcal A,\mathcal X)$ are quasi-categories by
\cref{prop:marked_mapping_fibration}, the right-hand vertical map is a
weak equivalence for every $\mathcal U\in\HGen$ if and only if
$i^*$ is a weak equivalence in $\sset_{(\infty,1)}$, by
\cref{lem:homotopy_generators_detect}. This proves the last assertion,
as the left-hand vertical maps are the maps defining locality with
respect to the maps $i\times\mathcal U^\flat$.
\end{proof}

\subsection{Initial-vertex restriction}
\label{subsec:enriched_marked_simplices}

In this section, we show that the inclusion of the initial vertex of a
marked edge over $\simFin^\sharp$ is a trivial cofibration in
$\cocart$.

\begin{convention}[{dual of \cite[Definition~3.1.1.1]{lurie2009higher}}]
\label{conv:cocart-anodyne}
A map of marked simplicial sets is \emph{coCartesian marked anodyne}
if its opposite is marked anodyne. By the dual of
\cite[Remark~3.1.3.4]{lurie2009higher}, every coCartesian marked
anodyne map over $\simFin^\sharp$ is a trivial cofibration in
$\cocart$.
\end{convention}

Given an edge $f$ of $\simFin$, we denote by
\begin{equation}
\label{eq:initial-vertex-inclusion}
j_f:(\Delta\{0\})^\sharp_{f(0)}
\hookrightarrow\Delta[1]^\sharp_f
\end{equation}
the inclusion of the initial vertex.

\begin{lemma}[{dual of \cite[Definition~3.1.1.1\textup{(2)}]{lurie2009higher},
in dimension one}]
\label{lem:initial-vertex-cocart-anodyne}
The map $j_f$ is coCartesian marked anodyne. In particular, it is a
trivial cofibration in $\cocart$.
\end{lemma}

\begin{proof}
The opposite of the map $j_f$ is the inclusion
\[
(\Delta\{1\})^\flat=(\Delta\{1\})^\sharp
\hookrightarrow\Delta[1]^\sharp
\]
over $((\simFin)^{\mathrm{op}})^\sharp$, which is the generating
marked anodyne map of
\cite[Definition~3.1.1.1\textup{(2)}]{lurie2009higher} in dimension
one, as $\Lambda^1[1]=\Delta\{1\}$ and the edge $\Delta\{0,1\}$ is
the only nondegenerate edge of $\Delta[1]$. The last assertion
follows from \cref{conv:cocart-anodyne}.
\end{proof}

\subsection{Inert spans and special coCartesian fibrations}
\label{subsec:enriched_specialness}

In this section, we define special coCartesian fibrations over
$\simFin$ in terms of spans analogous to those of
\cref{subsec:internal_specialness}, and we show that they coincide
with the symmetric monoidal $\infty$-categories of
\cite[Definition~2.0.0.7]{LurieHA}. Given a coCartesian fibration
$p:\mathcal E\to\simFin$ and an edge $f:a\to b$ of $\simFin$, we
define
\[
\Cocart_p(f)=\Map^\flat
\bigl(\Delta[1]^\sharp_f,\mathcal E^\natural\bigr),
\]
whose vertices are the $p$-coCartesian edges over $f$. The
restrictions to the vertices of $\Delta[1]$ give a span of
simplicial sets
\begin{equation}
\label{eq:single_cocartesian_span}
\mathcal E_a\xleftarrow{\ s_f^p\ }
\Cocart_p(f)\xrightarrow{\ t_f^p\ }\mathcal E_b.
\end{equation}

For $k\geq1$, we denote by $\pi_k^p$ the canonical projection to
$\mathcal E_{\langle k\rangle}$ from the iterated pullback
\[
\Cocart_p(\rho_1)
\times_{\mathcal E_{\langle k\rangle}}\cdots
\times_{\mathcal E_{\langle k\rangle}}\Cocart_p(\rho_k),
\]
of the maps $s_{\rho_i}^p$, and we consider the map
\begin{equation}
\label{eq:marked_inert_target}
t_k^p=(t_{\rho_i}^p\pr_i)_{i=1}^k:
\Cocart_p(\rho_1)
\times_{\mathcal E_{\langle k\rangle}}\cdots
\times_{\mathcal E_{\langle k\rangle}}\Cocart_p(\rho_k)
\longrightarrow\mathcal E_{\langle1\rangle}^{\,k}.
\end{equation}
For $k=0$, we consider the span
\begin{equation}
\label{eq:marked_nullary_span}
\mathcal E_{\langle0\rangle}
\xleftarrow{\ \pi_0^p=\mathrm{id}\ }
\mathcal E_{\langle0\rangle}
\xrightarrow{\ t_0^p\ }\Delta[0].
\end{equation}

\begin{definition}[{cf. \cite{Seg74} and \cite[Definition~2.0.0.7]{LurieHA}}]
\label{def:special_cocartesian_fibration}
A coCartesian fibration $p:\mathcal E\to\simFin$ is
\emph{special} if the map $t_k^p$ from
\eqref{eq:marked_inert_target} and \eqref{eq:marked_nullary_span} is a
weak equivalence in $\sset_{(\infty,1)}$ for every $k\geq0$.
\end{definition}

As in \cref{lem:internal_fan_identification}, the source and target of
$t_k^p$ can be expressed in terms of the inert fan of
\cref{def:inert_fan}.

\begin{lemma}
\label{lem:inert_fan_mapping_objects}
Let $p:\mathcal E\to\simFin$ be a coCartesian fibration. There are
isomorphisms
\[
\Map^\flat\bigl((\mathcal J_k)^\sharp_{\beta_k},
\mathcal E^\natural\bigr)
\cong
\Cocart_p(\rho_1)
\times_{\mathcal E_{\langle k\rangle}}\cdots
\times_{\mathcal E_{\langle k\rangle}}\Cocart_p(\rho_k)
\qquad(k\geq1)
\]
and
\[
\Map^\flat\bigl((\mathcal V_k)^\sharp_{\beta_k|\mathcal V_k},
\mathcal E^\natural\bigr)
\cong\mathcal E_{\langle1\rangle}^{\,k}
\qquad(k\geq0),
\]
under which the restriction maps along the inclusions
$\Delta\{c\}\hookrightarrow\mathcal J_k$ and
$\mathcal V_k\hookrightarrow\mathcal J_k$ correspond to the maps
$\pi_k^p$ and $t_k^p$, respectively. For $k=0$, we have
$\Map^\flat\bigl((\mathcal J_0)^\sharp_{\beta_0},
\mathcal E^\natural\bigr)=\mathcal E_{\langle0\rangle}$, and the
restriction maps are the maps of \eqref{eq:marked_nullary_span}.
\end{lemma}

\begin{proof}
Let $k\geq1$. An $r$-simplex of
$\Map^\flat\bigl((\mathcal J_k)^\sharp_{\beta_k},
\mathcal E^\natural\bigr)$ is a map of marked simplicial sets
\[
\Delta[r]^\flat\times\mathcal J_k^\sharp
\longrightarrow\mathcal E^\natural
\]
over $\beta_k\pr_2$. Since $\mathcal J_k$ is the union of $k$
copies of $\Delta[1]$ along the vertex $c$, such a map amounts to
a family of maps $F_i:\Delta[r]\times\Delta[1]\to\mathcal E$ over
$\rho_i\pr_2$, for $1\leq i\leq k$, that have the same restriction
to $\Delta[r]\times\{0\}$ and send the horizontal edges to
$p$-coCartesian edges, i.e., to a family of $r$-simplices $F_i$
of $\Cocart_p(\rho_i)$ with the same image under the maps
$s_{\rho_i}^p$. This correspondence is compatible with the
simplicial structure in $r$, which gives the first isomorphism, and
the restriction along $\Delta\{c\}\hookrightarrow\mathcal J_k$ then
corresponds to $\pi_k^p$. Since
$\mathcal V_k=\coprod_{i=1}^k\Delta[0]$ and $\beta_k$ sends every
vertex of $\mathcal V_k$ to $\langle1\rangle$, we get the second
isomorphism, under which the restriction along
$\mathcal V_k\hookrightarrow\mathcal J_k$ sends a family
$(F_i)_i$ to the family $(F_i|_{\Delta[r]\times\{1\}})_i$, and
hence corresponds to $t_k^p$. For $k=0$, we have
$\mathcal J_0=\Delta[0]$ over $\langle0\rangle$ and
$\mathcal V_0=\varnothing$, so that the two mapping objects are
$\mathcal E_{\langle0\rangle}$ and $\Delta[0]$, respectively.
\end{proof}

In what follows, we use the first isomorphism of
\cref{lem:inert_fan_mapping_objects} to write the \emph{inert span}
formed by the maps $\pi_k^p$ and $t_k^p$ as
\begin{equation}
\label{eq:marked_inert_span}
\mathcal E_{\langle k\rangle}
\xleftarrow{\ \pi_k^p\ }
\Map^\flat\bigl((\mathcal J_k)^\sharp_{\beta_k},
                 \mathcal E^\natural\bigr)
\xrightarrow{\ t_k^p\ }\mathcal E_{\langle1\rangle}^{\,k}.
\end{equation}

\begin{ex}
\label{ex:enriched_discrete_functor}
Let $F:\Fin\to\Set$ be a functor with category of elements
$\int F$. The map $p:N(\int F)\to\simFin$ is a left fibration, and
so every edge of $N(\int F)$ is $p$-coCartesian
\cite[\S2.1.1]{lurie2009higher}. Its fibers are the discrete simplicial
sets $F\langle k\rangle$, and the map $t_k^p$ is identified with
the map $(F(\rho_1),\ldots,F(\rho_k))$. Hence, as in
\cref{ex:discrete_commutative_monoid}, the coCartesian fibration $p$
is special if and only if $F$ is a special $\Gamma$-set.
\end{ex}

\begin{ex}
\label{ex:ordinary_symmetric_monoidal}
Let $\mathrm C$ be a small symmetric monoidal category, and let
$\mathrm C^\otimes\to\Fin$ be the opfibration of
\cite[Construction~2.0.0.1]{LurieHA}, whose objects are the finite
sequences $(x_1,\ldots,x_k)$ of objects of $\mathrm C$ and whose
morphisms $(x_1,\ldots,x_k)\to(y_1,\ldots,y_{k'})$ over a map
$f:\langle k\rangle\to\langle k'\rangle$ in $\Fin$ are the
families of morphisms
\[
\bigotimes_{i\in f^{-1}(j)}x_i\longrightarrow y_j
\qquad(1\leq j\leq k')
\]
in $\mathrm C$. Its nerve is a coCartesian fibration
$p:\mathcal E=N(\mathrm C^\otimes)\to\simFin$ with fibers
$\mathcal E_{\langle k\rangle}\cong N(\mathrm C)^k$, and an edge of
$\mathcal E$ is $p$-coCartesian if and only if all its components
are isomorphisms. Hence the mapping object
$\Map^\flat\bigl((\mathcal J_k)^\sharp_{\beta_k},
\mathcal E^\natural\bigr)$ is the nerve of the category of families
$(\alpha_i:x_i\xrightarrow{\cong}y_i)_{i=1}^k$ of isomorphisms in
$\mathrm C$, and the map $t_k^p$ is the nerve of the functor $T$
sending such a family to $(y_i)_{i=1}^k$. The functor $T$ is an
equivalence of categories, with inverse given by
$(y_i)_i\mapsto(\id_{y_i})_i$, and so $t_k^p$ is a weak
equivalence in $\sset_{(\infty,1)}$ for every $k\geq1$. As
$\mathcal E_{\langle0\rangle}=\Delta[0]$, the coCartesian fibration
$p$ is special.
\end{ex}

\begin{lemma}
\label{lem:cocart_transport_span}
Let $p:\mathcal E\to\simFin$ be a coCartesian fibration and let
$f:a\to b$ be an edge of $\simFin$. Then the simplicial set
$\Cocart_p(f)$ is a quasi-category, the map $s_f^p$ is a trivial
fibration of simplicial sets, and the map
\[
(s_f^p,t_f^p):\Cocart_p(f)\longrightarrow
\mathcal E_a\times\mathcal E_b
\]
is a fibration in $\sset_{(\infty,1)}$.
\end{lemma}

\begin{proof}
The simplicial set $\Cocart_p(f)$ is a quasi-category by
\cref{prop:marked_mapping_fibration}. By
\eqref{eq:marked_flat_mapping}, we have
$\Map^\flat\bigl((\Delta\{0\})^\sharp_a,\mathcal E^\natural\bigr)
=\mathcal E_a$ and
$\Map^\flat\bigl((\partial\Delta[1])^\sharp_{f|\partial\Delta[1]},
\mathcal E^\natural\bigr)=\mathcal E_a\times\mathcal E_b$, so that
the maps $s_f^p$ and $(s_f^p,t_f^p)$ are the restriction maps along
the inclusions $j_f$ and
$(\partial\Delta[1])^\sharp\hookrightarrow\Delta[1]^\sharp_f$,
respectively. Hence $(s_f^p,t_f^p)$ is a fibration in
$\sset_{(\infty,1)}$ by \cref{prop:marked_mapping_fibration}, and
$s_f^p$ is a trivial fibration by the same proposition, as $j_f$
is a trivial cofibration in $\cocart$ by
\cref{lem:initial-vertex-cocart-anodyne}.
\end{proof}

\begin{proposition}
\label{prop:marked_inert_span_fibrations}
Let $p:\mathcal E\to\simFin$ be a coCartesian fibration and let
$k\geq0$. Then the three objects of the span
\eqref{eq:marked_inert_span} are quasi-categories, the map
$\pi_k^p$ is a trivial fibration, and the map $t_k^p$ is a
fibration in $\sset_{(\infty,1)}$. In particular, the coCartesian
fibration $p$ is special if and only if the map $t_k^p$ is a
trivial fibration for every $k\geq0$.
\end{proposition}

\begin{proof}
For $k\geq1$, by \cref{lem:cocart_transport_span}, the spans
\[
\mathcal E_{\langle k\rangle}
\xleftarrow{\ s_{\rho_i}^p\ }\Cocart_p(\rho_i)
\xrightarrow{\ t_{\rho_i}^p\ }\mathcal E_{\langle1\rangle},
\]
for $1\leq i\leq k$, satisfy the hypotheses of
\cref{lem:common_source_fibrations}. Hence the map $\pi_k^p$ is a
trivial fibration, its source is a quasi-category, and the map
$t_k^p$ is a fibration in $\sset_{(\infty,1)}$. The case $k=0$
follows from the case $k=0$ of \cref{lem:common_source_fibrations}
applied to \eqref{eq:marked_nullary_span}. The last assertion holds
as the maps $t_k^p$ are fibrations.
\end{proof}

\begin{rmk}
\label{rem:marked_arity_one}
As in \cref{rem:segality_small_arities}, the specialness condition is
automatic for $k=1$. Indeed, the projection
$\mathcal E_{\langle1\rangle}\times\Delta[1]
\to\mathcal E_{\langle1\rangle}$ sends the horizontal edges to
degenerate edges, which are $p$-coCartesian, and so it defines a
common section $s$ of the maps $s_{\rho_1}^p$ and
$t_{\rho_1}^p$. Since $s_{\rho_1}^p$ is a trivial fibration by
\cref{lem:cocart_transport_span}, the maps $s$ and
$t_1^p=t_{\rho_1}^p$ are weak equivalences, by the two-out-of-three
property.
\end{rmk}

\begin{proposition}
\label{prop:special_cocartesian_HA}
A coCartesian fibration $p:\mathcal E\to\simFin$ is special if and
only if it is a symmetric monoidal $\infty$-category in the sense of
\cite[Definition~2.0.0.7]{LurieHA}.
\end{proposition}

\begin{proof}
Recall that, given an edge $f:a\to b$ of $\simFin$, a functor
$f_!:\mathcal E_a\to\mathcal E_b$ induced by $f$ is obtained as
the restriction $f_!=h|_{\mathcal E_a\times\{1\}}$ of a map
$h:\mathcal E_a\times\Delta[1]\to\mathcal E$ over $f\pr_2$ that
restricts to the inclusion on $\mathcal E_a\times\{0\}$ and sends
each edge $\{x\}\times\Delta[1]$ to a $p$-coCartesian edge, and
that $f_!$ is well defined up to equivalence
\cite[\S5.2.1]{lurie2009higher}. By definition of $\Cocart_p(f)$,
such a map $h$ amounts to a section $\varsigma$ of the trivial
fibration $s_f^p$ of \cref{lem:cocart_transport_span}, with
$f_!=t_f^p\varsigma$. Hence, in $\Ho(\sset_{(\infty,1)})$, we have
\[
[f_!]=[t_f^p][\varsigma]=[t_f^p][s_f^p]^{-1}.
\]
Since $s_{\rho_i}^p\pr_i=\pi_k^p$, the $i$-th component of the map
$[t_k^p][\pi_k^p]^{-1}$ is
\[
[t_{\rho_i}^p][\pr_i][\pi_k^p]^{-1}
=[t_{\rho_i}^p][s_{\rho_i}^p]^{-1}=[(\rho_i)_!],
\]
and so, as products of quasi-categories are products in the homotopy
category by \cref{prop:model_category_recollections}(5), the map
$[t_k^p][\pi_k^p]^{-1}$ is the map
$\mathcal E_{\langle k\rangle}\to\mathcal E_{\langle1\rangle}^{\,k}$
with components $(\rho_i)_!$, which is required to be an equivalence
in \cite[Definition~2.0.0.7]{LurieHA}. Since $\pi_k^p$ is a weak
equivalence, this map is an isomorphism in
$\Ho(\sset_{(\infty,1)})$ if and only if $t_k^p$ is a weak
equivalence in $\sset_{(\infty,1)}$, by
\cref{prop:model_category_recollections}(3). For $k=0$, both
definitions require the map
$\mathcal E_{\langle0\rangle}\to\Delta[0]$ to be a weak equivalence.
\end{proof}

\subsection{The special coCartesian model structure}
\label{subsec:enriched_special_model_structure}

In this section, we localize the model structure $\cocart$ to obtain
a model structure whose fibrant objects are the naturally marked
special coCartesian fibrations; a related localization of the category
of marked simplicial sets over finite pointed sets appears in
\cite[\S8]{Sharma2021}. We consider the set of maps over
$\simFin^\sharp$
\begin{equation}
\label{eq:enriched_special_localizers}
\mathsf S_{\mathrm{sp}}^+
=\left\{
\mathcal V_k^\sharp\times\mathcal U^\flat
\longrightarrow\mathcal J_k^\sharp\times\mathcal U^\flat
\ \middle|\ k\geq0,\ \mathcal U\in\HGen
\right\},
\end{equation}
where the structure maps are induced by $\beta_k$, and we write
\[
\cocartsp=L_{\mathsf S_{\mathrm{sp}}^+}\cocart.
\]

\begin{theorem}
\label{thm:enriched_special_model_structure}
The model structure $\cocartsp$ exists, and it is left proper and
combinatorial. Its cofibrations are the monomorphisms and its fibrant
objects are the naturally marked special coCartesian fibrations over
$\simFin$. In particular, every object is cofibrant. Moreover, a map
between fibrant objects is a weak equivalence in $\cocartsp$ if and
only if it is a weak equivalence in $\cocart$, if and only if it
induces a weak equivalence in $\sset_{(\infty,1)}$ on every fiber.
\end{theorem}

\begin{proof}
The existence of $\cocartsp$, its left properness and
combinatoriality, and the description of its cofibrations follow from
\cref{thm:left_bousfield_localization} applied to the model structure
of \cref{thm:cocartesian_model_structure}. Let
$p:\mathcal E\to\simFin$ be a coCartesian fibration and let
$k\geq0$. By \cref{lem:marked_shape_locality} applied to the
inclusion
$(\mathcal V_k)^\sharp_{\beta_k|\mathcal V_k}\hookrightarrow
(\mathcal J_k)^\sharp_{\beta_k}$ and
\cref{lem:inert_fan_mapping_objects}, the object
$\mathcal E^\natural$ is local with respect to the maps
$\mathcal V_k^\sharp\times\mathcal U^\flat\to
\mathcal J_k^\sharp\times\mathcal U^\flat$, for
$\mathcal U\in\HGen$, if and only if the map $t_k^p$ is a weak
equivalence in $\sset_{(\infty,1)}$. Hence the fibrant objects of
$\cocartsp$ are the naturally marked special coCartesian fibrations.
The description of the weak equivalences between fibrant objects
follows from \cref{prop:localization_recollections}(1) and
\cref{prop:cocartesian_weak_equivalences}.
\end{proof}

\begin{rmk}
\label{rem:compat_generator_tensors}
The products with the objects $\mathcal U^\flat$ in
\eqref{eq:enriched_special_localizers} are needed. Indeed, by
\eqref{eq:marked_mapping_core}, a naturally marked coCartesian
fibration $\mathcal E^\natural$ is local with respect to the
inclusion $\mathcal V_k^\sharp\hookrightarrow\mathcal J_k^\sharp$
over $\beta_k$ if and only if the map induced by $t_k^p$ on cores
is a weak equivalence in $\sset_{(\infty,0)}$, and this condition
does not imply that $t_k^p$ is a weak equivalence in
$\sset_{(\infty,1)}$. For instance, let $\mathcal C$ be the nerve
of the monoid $(\mathbb Z_{\geq0},+)$, regarded as a category with
one object, and consider the coCartesian fibration
$p:\mathcal C\times\simFin\to\simFin$ given by the projection. Since
$0$ is the only invertible element of this monoid, the
$\mathcal C$-component of a $p$-coCartesian edge is an identity,
and so
\[
\Map^\flat\bigl((\mathcal J_k)^\sharp_{\beta_k},
                 (\mathcal C\times\simFin)^\natural\bigr)
\cong\mathcal C,
\]
with $t_k^p$ given by the diagonal $\mathcal C\to\mathcal C^k$.
The diagonal induces an isomorphism on cores, as the core of
$\mathcal C$ is $\Delta[0]$, but it is not a weak equivalence in
$\sset_{(\infty,1)}$ for $k=0$ and for $k\geq2$, as the induced
map on endomorphism monoids is not an isomorphism. Hence
$(\mathcal C\times\simFin)^\natural$ is local with respect to the
inclusions $\mathcal V_k^\sharp\hookrightarrow\mathcal J_k^\sharp$,
but $p$ is not special.
\end{rmk}

\subsection{Boundary lifting tests for specialness}
\label{subsec:enriched_boundary_tests}

In this section, we characterize the special coCartesian fibrations by
a right lifting property against the pushout-products of the boundary
inclusions with the inclusions
$\mathcal V_k\hookrightarrow\mathcal J_k$, and we show that
$\cocartsp$ is also the left Bousfield localization of $\cocart$
at these pushout-products; compare
\cite[Theorem~4.46 and its proof]{barwick2010bousfield}.

\begin{proposition}
\label{prop:special_localization_explicit_lifts}
For $k,r\geq0$, let $w_{k,r}$ denote the monomorphism of marked
simplicial sets over $\simFin^\sharp$
\begin{equation}
\label{eq:special_boundary_localizers}
w_{k,r}:
\bigl(\mathcal J_k^\sharp\times(\partial\Delta[r])^\flat\bigr)
\amalg_{\mathcal V_k^\sharp\times(\partial\Delta[r])^\flat}
\bigl(\mathcal V_k^\sharp\times\Delta[r]^\flat\bigr) \longrightarrow
\mathcal J_k^\sharp\times\Delta[r]^\flat,
\end{equation}
with structure map induced by $\beta_k$, and let
$\mathsf W=\{w_{k,r}\mid k,r\geq0\}$. For a coCartesian fibration
$p:\mathcal E\to\simFin$, the following conditions are equivalent.
\begin{enumerate}[label=\textup{(\roman*)}]
\item The object $\mathcal E^\natural$ is fibrant in $\cocartsp$.
\item The map $t_k^p$ is a weak equivalence in
$\sset_{(\infty,1)}$ for every $k\geq0$.
\item The map $t_k^p$ is a trivial fibration of simplicial sets for
every $k\geq0$.
\item The map $\mathcal E^\natural\to\simFin^\sharp$ has the right
lifting property against the map $w_{k,r}$ for all $k,r\geq0$.
\end{enumerate}
Moreover, we have
\[
L_{\mathsf W}\cocart=\cocartsp.
\]
\end{proposition}

\begin{proof}
Conditions \textup{(i)} and \textup{(ii)} are equivalent by
\cref{thm:enriched_special_model_structure}, and conditions
\textup{(ii)} and \textup{(iii)} are equivalent by
\cref{prop:marked_inert_span_fibrations}. To show that \textup{(iii)}
and \textup{(iv)} are equivalent, recall that the trivial fibrations
of simplicial sets are the maps with the right lifting property
against the boundary inclusions, and consider a lifting problem in
simplicial sets
\[
\begin{tikzcd}[ampersand replacement=\&, sep=large]
\partial\Delta[r]\arrow[r,"a"]\arrow[d,hook]
\&\Map^\flat\bigl((\mathcal J_k)^\sharp_{\beta_k},
\mathcal E^\natural\bigr)\arrow[d,two heads,"t_k^p"]\\
\Delta[r]\arrow[r,"b"']\arrow[ur,dashed]
\&\mathcal E_{\langle1\rangle}^{\,k}.
\end{tikzcd}
\]
By \cref{lem:inert_fan_mapping_objects} and adjunction, the maps $a$
and $b$ correspond to maps of marked simplicial sets over
$\simFin^\sharp$
\[
\mathcal J_k^\sharp\times(\partial\Delta[r])^\flat
\longrightarrow\mathcal E^\natural,
\qquad
\mathcal V_k^\sharp\times\Delta[r]^\flat
\longrightarrow\mathcal E^\natural,
\]
which agree on $\mathcal V_k^\sharp\times(\partial\Delta[r])^\flat$,
i.e., to a map from the source of $w_{k,r}$ to $\mathcal E^\natural$
over $\simFin^\sharp$, and the lifts in the above square correspond
to the extensions of this map along $w_{k,r}$. Hence \textup{(iii)}
and \textup{(iv)} are equivalent. It remains to show that the model
structures $L_{\mathsf W}\cocart$, which exists by
\cref{thm:left_bousfield_localization}, and $\cocartsp$ coincide. As
both are left Bousfield localizations of $\cocart$, it is enough to
show that they have the same fibrant objects. Let
$\mathcal X=\mathcal E^\natural$ be a fibrant object of $\cocart$.
By the enrichment of $\cocart$ over $\sset_{(\infty,0)}$, the
restriction map
\[
\Map^\sharp(\operatorname{cod}(w_{k,r}),\mathcal X)
\longrightarrow
\Map^\sharp(\operatorname{dom}(w_{k,r}),\mathcal X)
\]
is a fibration in $\sset_{(\infty,0)}$ between models of the
derived mapping spaces, by \eqref{eq:marked_mapping_core}. Hence, if
$\mathcal X$ is local with respect to the maps in $\mathsf W$,
then these restriction maps are trivial fibrations, and in particular
they are surjective on vertices, so that \textup{(iv)} holds and
$p$ is special. Conversely, suppose that $t_k^p$ is a trivial fibration for
every $k\geq0$. By \cref{lem:inert_fan_mapping_objects} and
adjunction, the restriction map along $w_{k,r}$ between the mapping
objects $\Map^\flat$ into $\mathcal X$ is identified with the map
\[
\Map^\flat\bigl((\mathcal J_k)^\sharp_{\beta_k},
                   \mathcal E^\natural\bigr)^{\Delta[r]}
\longrightarrow
\Map^\flat\bigl((\mathcal J_k)^\sharp_{\beta_k},
                   \mathcal E^\natural\bigr)^{\partial\Delta[r]}
\mathop{\times}_{(\mathcal E_{\langle1\rangle}^{\,k})^{\partial\Delta[r]}}
(\mathcal E_{\langle1\rangle}^{\,k})^{\Delta[r]},
\]
which is a trivial fibration between quasi-categories, as
$\sset_{(\infty,1)}$ is cartesian and by
\cref{prop:marked_mapping_fibration}. Hence the induced map on cores
is a weak equivalence in $\sset_{(\infty,0)}$ by
\cite[Lemma~3.1.3.2]{lurie2009higher}, and, by
\eqref{eq:marked_mapping_core}, this map on cores is the above
restriction map between the mapping objects $\Map^\sharp$. This
shows that $\mathcal X$ is local with respect to $w_{k,r}$.
\end{proof}

\section{Comparison of the two descriptions}
\label{sec:comparisons}

In this section, we take $n=1$ and $\Minfn=\sset_{(\infty,1)}$,
and we compare the left double $(\infty,1)$-fibrations over
$\bisimFin$ of \cref{sec:internal_description} with the coCartesian
fibrations over $\simFin$ of \cref{sec:enriched_description}. In
\cref{subsec:comparison_internal_to_enriched}, we construct the
staircase $\Stair(P):\Stair(\mathbb E)\to\simFin$ of a left double
$(\infty,1)$-fibration $P:\mathbb E\to\bisimFin$ and show that it
is a coCartesian fibration, which is special if and only if $P$ is
special. In \cref{subsec:comparison_enriched_to_internal}, we
construct the left double $(\infty,1)$-fibration
$\Int(\mathcal E^\natural)\to\bisimFin$ of coCartesian strings of a
coCartesian fibration $p:\mathcal E\to\simFin$, and we prove the
analogous statement. We then show in
\cref{subsec:marked_diagonal_adjunction} that the functor $\Int$
admits a left adjoint $\Diag$, the marked diagonal, and in
\cref{subsec:comparison_quillen_equivalences} that the adjunction
$\Diag\dashv\Int$ induces the horizontal Quillen equivalences in the
following diagram, in which the vertical identity functors are left
Quillen; see
\cref{thm:left_double_cocart_quillen_equivalence,thm:special_quillen_equivalence}.
\[
\begin{tikzcd}[ampersand replacement=\&, sep=large]
\ldblfib
 \arrow[r,shift left=1.2ex,"\Diag"]
 \arrow[r,phantom,"\bot" description]
 \arrow[d,"\id"']
\&\cocart
 \arrow[l,shift left=1.2ex,"\Int"]
 \arrow[d,"\id"]\\
\ldblfibsp
 \arrow[r,shift left=1.2ex,"\Diag"]
 \arrow[r,phantom,"\bot" description]
\&\cocartsp
 \arrow[l,shift left=1.2ex,"\Int"].
\end{tikzcd}
\]
The staircase is related to this adjunction as follows. For a left
double $(\infty,1)$-fibration $P:\mathbb E\to\bisimFin$, we
construct in \cref{const:diagonal_to_staircase} a map
$\psi_{\mathbb E}:\Diag(\mathbb E)\to\Stair(\mathbb E)^\natural$,
which is a weak equivalence in $\cocart$ with fibrant target by
\cref{prop:diagonal_to_staircase_cocart_equivalence}. Since every
object of $\ldblfib$ is cofibrant, the staircase is therefore a
fibrant model for the value at $\mathbb E$ of the left derived
functor of $\Diag$. A comparison between marked simplicial sets and
bisimplicial sets as models of Cartesian fibrations, in the spirit of
\cite{JoyalTierney2007}, is given in
\cite[\S4.1]{Rasekh2021QuasiCategories}, in terms of complete Segal
objects in right fibrations.

\begin{rmk}
\label{rem:nuiten_dictionary}
The equivalence of homotopy theories underlying
\cref{thm:left_double_cocart_quillen_equivalence} can also be obtained
from Nuiten's work on straightening for Segal spaces
\cite{Nuiten2023}, and we do not claim it as new. Nuiten considers two
kinds of fibrations of double $\infty$-categories, namely those that
are left fibrations horizontally and Cartesian fibrations vertically,
and those that are coCartesian fibrations horizontally and right
fibrations vertically, and shows in \cite[Theorem~3.1]{Nuiten2023}
that they are equivalent through a \emph{reflection}, in which a
horizontal arrow is a horizontal arrow followed by a vertical arrow in
a fiber. The edges of the staircase of
\cref{const:staircase_codiagonal} are of this form, and the staircase
plays the role of the reflection in our setting. For a base given by
a category $\mathrm C$, \cite[Corollary~4.35]{Nuiten2023} identifies
Segal copresheaves of $(\infty,1)$-categories on $\mathrm C$ with
coCartesian fibrations over $\mathrm C$, which for
$\mathrm C=\Fin$ corresponds to
\cref{thm:left_double_cocart_quillen_equivalence}; compare
\cref{rem:general_base}. Nuiten also remarks that the reflection
functors arise as right Quillen functors on marked bisimplicial
spaces, and that the homotopy theory in question is a left Bousfield
localization of an injective model structure
\cite[Remarks~3.2 and~3.9]{Nuiten2023}. We compare instead a
localization of the Reedy model structure on bisimplicial sets with
Lurie's coCartesian model structure, through the adjunction
$\Diag\dashv\Int$. The localizations for symmetric monoidal
structures, including the condition for $k=0$, and the comparison of
the spans of
\cref{subsec:internal_specialness,subsec:enriched_specialness} are not
treated in \cite{Nuiten2023}.
\end{rmk}

\subsection{From left double fibrations to coCartesian fibrations}
\label{subsec:comparison_internal_to_enriched}
\subsubsection{The staircase construction}

In this section, we recall the total simplicial set of a bisimplicial
set of Artin and Mazur \cite{ArtinMazur1966VanKampen}, in the form of
\cite[Definition~3.1]{Stevenson2025TotalSimplicialSet} with the two
variables interchanged, which we call the \emph{staircase}. The shapes
indexing its simplices are those of Waldhausen's
$S_\bullet$-construction
\cite{Waldhausen1985Spaces,BergnerEtAl2021TwoSegalObjects}. Since the
horn-filling arguments of the next sections are carried out
componentwise, we describe the simplicial structure of the staircase
in components.

Let $\mathbb X$ be a bisimplicial set, with the conventions of
\cref{subsec:prelim_simplicial_models}. Given a map $\alpha$ in
$\Delta$, we write $\alpha_h^*$ and $\alpha_v^*$ for the induced
maps in the horizontal and vertical directions, respectively, and we
write $\delta_i:[r-1]\hookrightarrow[r]$ for the injection omitting
$i$.

\begin{construction}[{\cite[Definition~3.1]{Stevenson2025TotalSimplicialSet}, with the variables interchanged}]
\label{const:staircase_codiagonal}
For $q\geq 0$, we define $\Stair(\mathbb X)_q$ to be the set of
tuples
$
\xi=(\xi_0,\ldots,\xi_q)
$
with
$
\xi_j\in\mathbb X_{q-j,j}
$
such that, for $1\leq j\leq q$,
\begin{equation}
\label{eq:staircase_relation}
d_0^h\xi_{j-1}=d_j^v\xi_j.
\end{equation}
Given a map $\alpha:[m]\to[q]$ in $\Delta$ and $0\leq i\leq m$,
let $\alpha_i^-:[i]\to[\alpha(i)]$ and
$\alpha_i^+:[m-i]\to[q-\alpha(i)]$ be the maps given by
\[
\alpha_i^-(t)=\alpha(t)
\qquad\text{and}\qquad
\alpha_i^+(t)=\alpha(i+t)-\alpha(i).
\]
We define the map
\[
\Stair(\mathbb X)(\alpha):
\Stair(\mathbb X)_q\longrightarrow\Stair(\mathbb X)_m
\]
by
\begin{equation}
\label{eq:staircase_alpha}
\bigl(\Stair(\mathbb X)(\alpha)(\xi)\bigr)_i
=
(\alpha_i^+)_{h}^{*}(\alpha_i^-)_{v}^{*}
\xi_{\alpha(i)}.
\end{equation}
\end{construction}

\begin{ex}
\label{ex:low_dimensional_staircases}
We draw the horizontal direction from left to right and the vertical
direction downward. An edge $\xi=(\xi_0,\xi_1)$ of
$\Stair(\mathbb X)$ is depicted as
\[
\begin{tikzcd}[ampersand replacement=\&, sep=large]
\boxed{x_{00}}\arrow[r,"\xi_0"]
\&x_{10}\arrow[d,"\xi_1"]\\
\&\boxed{x_{11}}
\end{tikzcd}
\]
where $\xi_0\in\mathbb X_{1,0}$ is the horizontal $1$-simplex, $\xi_1\in\mathbb X_{0,1}$ is the vertical $1$-simplex, and the boxed vertices are the vertices of $\xi$. The relations \eqref{eq:staircase_relation} read
\[d_0^h\xi_0=d_1^v\xi_1=x_{10}.\]
A $2$-simplex
$\xi=(\xi_0,\xi_1,\xi_2)$ of $\Stair(\mathbb X)$ is depicted as
\[
\begin{tikzcd}[ampersand replacement=\&, sep=large]
\boxed{x_{00}}\arrow[r,"h_{01}"]
 \arrow[rr,bend left=30,"h_{02}\ (\xi_0)"]
\&x_{10}\arrow[r,"h_{12}"]\arrow[d,"v^1_{01}"']
 \arrow[dr,phantom,"\xi_1" description]
\&x_{20}\arrow[d,"v^2_{01}"]
 \arrow[dd,bend left=50,"v^2_{02}\ (\xi_2)"]\\
\&\boxed{x_{11}}\arrow[r,"h'_{12}"']
\&x_{21}\arrow[d,"v^2_{12}"]\\
\&\&\boxed{x_{22}},
\end{tikzcd}
\]
where $\xi_0\in\mathbb X_{2,0}$ is the horizontal $2$-simplex with
long edge $h_{02}$, $\xi_1\in\mathbb X_{1,1}$ is the square, and
$\xi_2\in\mathbb X_{0,2}$ is the vertical $2$-simplex with long
edge $v^2_{02}$. The relations \eqref{eq:staircase_relation} read
\[
 d_0^h\xi_0=d_1^v\xi_1=h_{12},\qquad
 d_0^h\xi_1=d_2^v\xi_2=v^2_{01}.
\]
In general, the component $\xi_j$ of a $q$-simplex occupies the
box $[j,q]\times[0,j]$ of the grid $[q]\times[q]$, and the union of
these boxes is the triangular region
$\{(h,v)\mid 0\leq v\leq h\leq q\}$.
\end{ex}

\begin{proposition}
\label{prop:staircase_is_simplicial}
\Cref{const:staircase_codiagonal} defines a simplicial set
$\Stair(\mathbb X)$, naturally in the bisimplicial set $\mathbb X$.
\end{proposition}

\begin{proof}
This is classical, and we only indicate the argument. Let
$\xi\in\Stair(\mathbb X)_q$. For $0\leq a\leq b\leq q$, write
$\mathrm{init}_{a,b}:[a]\hookrightarrow[b]$ for the inclusion of the
initial segment and
$\mathrm{term}_{a,b}:[q-b]\hookrightarrow[q-a]$, $t\mapsto b-a+t$,
for the inclusion of the terminal segment. Iterating the relation
\eqref{eq:staircase_relation} gives
\[
(\mathrm{term}_{a,b})_h^*\xi_a=(\mathrm{init}_{a,b})_v^*\xi_b
\qquad\text{in }\mathbb X_{q-b,a}.
\]
Given a map $\alpha:[m]\to[q]$ in $\Delta$ and $1\leq i\leq m$,
put $a=\alpha(i-1)$ and $b=\alpha(i)$. Then
$\alpha_{i-1}^+\delta_0=\mathrm{term}_{a,b}\alpha_i^+$ and
$\alpha_i^-\delta_i=\mathrm{init}_{a,b}\alpha_{i-1}^-$, and so, by
the displayed identity, the tuple \eqref{eq:staircase_alpha} satisfies
the relations \eqref{eq:staircase_relation}. Functoriality follows
from the identities
$(\alpha\beta)_i^-=\alpha_{\beta(i)}^-\beta_i^-$ and
$(\alpha\beta)_i^+=\alpha_{\beta(i)}^+\beta_i^+$, and from the fact
that the horizontal and vertical structure maps commute. Finally, a
map of bisimplicial sets acts componentwise, which gives naturality.
\end{proof}

We will use the following description of the face and degeneracy maps
of $\Stair(\mathbb X)$.

\begin{rmk}
\label{rem:staircase_faces_degeneracies}
We denote by $D_r:\Stair(\mathbb X)_q\to\Stair(\mathbb X)_{q-1}$
and $S_r:\Stair(\mathbb X)_q\to\Stair(\mathbb X)_{q+1}$, for
$0\leq r\leq q$, the face and degeneracy maps of
$\Stair(\mathbb X)$. By \eqref{eq:staircase_alpha}, we have
\begin{equation}\label{eq:staircase_faces}
\bigl(D_r\xi\bigr)_i
=
\begin{cases}
d_{r-i}^h \xi_i, & 0\leq i<r,\\[1mm]
d_r^v \xi_{i+1}, & r\leq i\leq q-1,
\end{cases}
\end{equation}
and
\begin{equation}\label{eq:staircase_degeneracies}
\bigl(S_r\xi\bigr)_i
=
\begin{cases}
s_{r-i}^h \xi_i, & 0\leq i\leq r,\\[1mm]
s_r^v \xi_{i-1}, & r+1\leq i\leq q+1.
\end{cases}
\end{equation}
\end{rmk}

Writing $\mathbb X^{\mathrm{tr}}$ for the bisimplicial set given by
$(\mathbb X^{\mathrm{tr}})_{i,j}=\mathbb X_{j,i}$, the above
description shows that
\begin{equation}
\label{eq:staircase_total_identification}
\Stair(\mathbb X)=T(\mathbb X^{\mathrm{tr}}),
\end{equation}
where $T$ denotes the total simplicial set functor of
\cite[Definition~3.1]{Stevenson2025TotalSimplicialSet}. We now compute
the staircase of an external product.

\begin{lemma}[{\cite[Lemma~3.8]{Stevenson2025TotalSimplicialSet}, with the variables interchanged}]
\label{lem:staircase_external_product}
Given simplicial sets $\mathcal K$ and $\mathcal A$, there is an
isomorphism of simplicial sets
\begin{equation}
\label{eq:staircase_external_product}
\Stair(\mathcal K\boxtimes \mathcal A)\cong \mathcal K\times \mathcal A,
\end{equation}
natural in $\mathcal K$ and $\mathcal A$, which sends a
$q$-simplex $(z_j,x_j)_{j=0}^q$ of
$\Stair(\mathcal K\boxtimes\mathcal A)$ to the pair $(z_0,x_q)$.
\end{lemma}

\begin{proof}
This is \cite[Lemma~3.8]{Stevenson2025TotalSimplicialSet} under the
identification \eqref{eq:staircase_total_identification}; we recall
the argument, as the explicit formulas are used below. Write the
$j$-th component of a $q$-simplex of
$\Stair(\mathcal K\boxtimes\mathcal A)$ as
$(z_j,x_j)\in\mathcal K_{q-j}\times\mathcal A_j$. The relations
\eqref{eq:staircase_relation} read $z_j=d_0z_{j-1}$ and
$x_{j-1}=d_jx_j$, and so $z_j=z_0|_{[j,q]}$ and
$x_j=x_q|_{[0,j]}$, where the intervals $[j,q]$ and $[0,j]$ are
identified with the corresponding faces of $\Delta[q]$. Hence the
assignment $(z_j,x_j)_{j=0}^q\mapsto(z_0,x_q)$ is a bijection
$\Stair(\mathcal K\boxtimes\mathcal A)_q\cong
\mathcal K_q\times\mathcal A_q$, with inverse
\[
(z,x)\longmapsto\bigl(z|_{[j,q]},\,x|_{[0,j]}\bigr)_{j=0}^q.
\]
By \eqref{eq:staircase_alpha}, these bijections are compatible with
the simplicial structures, and they are natural in $\mathcal K$ and
$\mathcal A$ by construction.
\end{proof}

\begin{ex}
\label{ex:staircase_external_product}
For $\mathcal K=\mathcal A=\Delta[1]$, we get
$\Stair(\Delta[1]\boxtimes\Delta[1])\cong\Delta[1]\times\Delta[1]$,
where the diagonal edge of $\Delta[1]\times\Delta[1]$ corresponds to
the edge of the staircase formed by the horizontal edge
$(0,0)\to(1,0)$ and the vertical edge $(1,0)\to(1,1)$. For the
left double $(\infty,1)$-fibration $N(\int F)^h\to\bisimFin$ of
\cref{ex:discrete_commutative_monoid}, we get
\[
 \Stair\bigl(N(\textstyle\int F)^h\bigr)
 \cong N(\textstyle\int F),
\]
and the staircase is the coCartesian fibration
$N(\int F)\to\simFin$ of \cref{ex:enriched_discrete_functor}.
\end{ex}

\subsubsection{The staircase over the fixed base}

We now apply the staircase construction to bisimplicial sets over
$\bisimFin$.

\begin{lemma}
\label{lem:staircase_of_horizontal_base}
There is an isomorphism of simplicial sets
\[
\Stair(\bisimFin)\cong \simFin,
\]
under which a $q$-simplex $\sigma$ of $\simFin$ corresponds to
the tuple
$
\bigl(
\sigma,\sigma_{\geq1},\ldots,\sigma_{\geq q}
\bigr)
$,
where $\sigma_{\geq j}\coloneqq(d_0)^j\sigma$.
\end{lemma}

\begin{proof}
This follows directly from \cref{lem:staircase_external_product}
applied to $\bisimFin=\simFin\boxtimes\Delta[0]$, as
$\sigma|_{[j,q]}=(d_0)^j\sigma$.
\end{proof}

Given a map $P:\mathbb E\to\bisimFin$ of bisimplicial sets, we
identify $\Stair(\bisimFin)$ with $\simFin$ through the isomorphism
of \cref{lem:staircase_of_horizontal_base}, and we write
\[
\Stair(P):\Stair(\mathbb E)\longrightarrow \simFin
\]
for the induced map. An edge of $\Stair(\mathbb E)$ over a morphism
$f:a\to b$ of $\Fin$ is a pair $(h,v)$ of a horizontal edge
$h\in\mathbb E_{1,0}$ with $P_1(h)=f$ and an edge
$v\in\mathbb E_{0,1}$ of the simplicial set $\mathbb E_0$ such
that $d_0^hh=d_1^vv$; its source is $d_1^hh$ and its target is
$d_0^vv$. Note that $P_0(v)=s_0b$, as $(\bisimFin)_0$ is
discrete.

\subsubsection{Horizontal extension and detection of equivalences}
\label{subsubsec:horizontal_left_extension_converse}

In the rest of \cref{subsec:comparison_internal_to_enriched}, we fix a
left double $(\infty,1)$-fibration $P:\mathbb E\to\bisimFin$. Given
a monomorphism
$i:\mathcal K\hookrightarrow \mathcal L$ of simplicial sets, recall
from \cref{lem:cellular_restriction} that the relative restriction map
\[
R_i:
\mathbb E(\mathcal L)
\longrightarrow
\mathbb E(\mathcal K)\times_{\bisimFin(\mathcal K)}\bisimFin(\mathcal L)
\]
is a fibration in $\sset_{(\infty,1)}$. In this section, we study
the monomorphisms $i$ for which $R_i$ is a trivial fibration. For
$q\geq1$ and $0\leq s\leq q$, we write
$d_s\Delta[q-1]\subseteq\Delta[q]$ for the image of the map
$\delta_s:\Delta[q-1]\hookrightarrow\Delta[q]$, i.e., for the face
opposite to the vertex $s$.

\begin{definition}
\label{def:class_T}
We denote by $\mathsf T$ the class of monomorphisms
$i:\mathcal K\hookrightarrow \mathcal L$ of simplicial sets such
that the map $R_i$ is a trivial fibration of simplicial sets.
\end{definition}

We first study the closure properties of the class $\mathsf T$.

\begin{lemma}
\label{lem:T_closure}
The class $\mathsf T$ satisfies the following properties.
\begin{enumerate}[label=(\arabic*)]
\item It is closed under pushouts along maps of simplicial sets.
\item It is closed under composition.
\item Given monomorphisms
$\mathcal K\xrightarrow{i}\mathcal L\xrightarrow{j}\mathcal M$ of
simplicial sets such that $i$ and $ji$ belong to $\mathsf T$,
the monomorphism $j$ belongs to $\mathsf T$.
\end{enumerate}
\end{lemma}

\begin{proof}
A pushout square of simplicial sets
\[
\begin{tikzcd}[ampersand replacement=\&, sep=large]
\mathcal K\arrow[r,"f"]\arrow[d,hook,"i"']
\&\mathcal K'\arrow[d,hook,"i'"]\\
\mathcal L\arrow[r,"g"']
\&\mathcal L'\arrow[ul,phantom,"\ulcorner",very near start]
\end{tikzcd}
\]
with $i$ a monomorphism induces a pullback square
\[
\begin{tikzcd}[ampersand replacement=\&, sep=large]
\mathbb E(\mathcal L')\arrow[r,"g^*"]\arrow[d,two heads,"R_{i'}"']
  \arrow[dr,phantom,"\lrcorner",very near start]
\&\mathbb E(\mathcal L)\arrow[d,two heads,"R_i"]\\
\mathbb E(\mathcal K')\times_{\bisimFin(\mathcal K')}\bisimFin(\mathcal L')
  \arrow[r,"{(f^*,g^*)}"']
\&\mathbb E(\mathcal K)\times_{\bisimFin(\mathcal K)}\bisimFin(\mathcal L),
\end{tikzcd}
\]
so that $R_{i'}$ is a pullback of $R_i$, which proves \textup{(1)}.
Next, given monomorphisms
$\mathcal K\xrightarrow{i}\mathcal L\xrightarrow{j}\mathcal M$ of
simplicial sets, the map $R_{ji}$ factors as
\[
\mathbb E(\mathcal M)\xrightarrow{R_j}
\mathbb E(\mathcal L)\times_{\bisimFin(\mathcal L)}\bisimFin(\mathcal M)
\xrightarrow{Q_i}
\mathbb E(\mathcal K)\times_{\bisimFin(\mathcal K)}\bisimFin(\mathcal M),
\]
where $Q_i$ is the pullback of $R_i$ along the map
\[
(\pr_1,j^*\pr_2):
\mathbb E(\mathcal K)\times_{\bisimFin(\mathcal K)}\bisimFin(\mathcal M)
\longrightarrow
\mathbb E(\mathcal K)\times_{\bisimFin(\mathcal K)}\bisimFin(\mathcal L).
\]
Hence, if $i$ and $j$ belong to $\mathsf T$, then so does
$ji$, which proves \textup{(2)}. Finally, if $i$ and $ji$ belong
to $\mathsf T$, then $Q_i$ and $R_{ji}=Q_iR_j$ are weak
equivalences in $\sset_{(\infty,1)}$, and hence so is $R_j$, by
the two-out-of-three property. As $R_j$ is a fibration in
$\sset_{(\infty,1)}$ by \cref{lem:cellular_restriction}, it is a
trivial fibration, which proves \textup{(3)}.
\end{proof}

By a \emph{left horn inclusion}, we mean an inclusion
\[
\Lambda^i[m]=\bigcup_{s\ne i}d_s\Delta[m-1]
\hookrightarrow\Delta[m],
\qquad m\geq1,\quad 0\leq i<m,
\]
so that the inner horn inclusions are left horn inclusions. We will
use the following filtration of the left horns.

\begin{lemma}[{\cite[Tag~05HW, Proposition~4.2.4.10]{kerodon}}]
\label{lem:initial_vertex_to_left_horn}
For every $m\geq1$ and $0\leq i<m$, the inclusion
$\Delta\{0\}\hookrightarrow\Lambda^i[m]$
is a finite composite of pushouts of left horn inclusions
$\Lambda^a[q]\hookrightarrow\Delta[q]$ with
$0\leq a<q<m$.
\end{lemma}

With these, we can show that the left horn inclusions belong to
$\mathsf T$.

\begin{proposition}
\label{prop:left_horns_in_T}
For every $m\geq1$ and every $0\leq i<m$, the left horn inclusion
$\Lambda^i[m]\hookrightarrow\Delta[m]$
belongs to $\mathsf T$.
\end{proposition}

\begin{proof}
We proceed by induction on $m$. For $m=1$, the inclusion
$\Lambda^0[1]=\Delta\{0\}\hookrightarrow\Delta[1]$ belongs to
$\mathsf T$, as its relative restriction map is $\lambda_1$, which
is a trivial fibration by \cref{prop:minimal_axioms}. Suppose now that
the statement holds in all dimensions less than $m$. Then the
inclusion $\Delta\{0\}\hookrightarrow\Lambda^i[m]$ belongs to
$\mathsf T$ by \cref{lem:initial_vertex_to_left_horn} and
\cref{lem:T_closure}(1) and~(2). Moreover, the composite
$\Delta\{0\}\hookrightarrow\Lambda^i[m]\hookrightarrow\Delta[m]$
belongs to $\mathsf T$, as its relative restriction map is
$\lambda_m$, which is a trivial fibration by
\cref{prop:minimal_axioms}. Hence the inclusion
$\Lambda^i[m]\hookrightarrow\Delta[m]$ belongs to $\mathsf T$ by
\cref{lem:T_closure}(3).
\end{proof}

We will also need the following variant, in which two faces are
removed from the simplex; compare
\cite[Proposition~2.12\textup{(ii)}]{Joyal2008Theory}.

\begin{lemma}
\label{lem:two_missing_horizontal_faces}
Let $1\leq\ell<m$, and put
\[
\mathcal K_{m,\ell}
=
\bigcup_{\substack{1\leq s\leq m\\s\neq\ell}}
d_s\Delta[m-1].
\]
Then the inclusion
$\mathcal K_{m,\ell}\hookrightarrow\Delta[m]$
belongs to $\mathsf T$.
\end{lemma}

\begin{proof}
The faces of $\Delta[m]$ of codimension one that are not contained
in $\mathcal K_{m,\ell}$ are $d_0\Delta[m-1]$ and
$d_\ell\Delta[m-1]$. For $s\geq1$, the intersection
$d_0\Delta[m-1]\cap d_s\Delta[m-1]$ is the face of
$d_0\Delta[m-1]\cong\Delta[m-1]$ opposite to the vertex $s$, i.e.,
the face $d_{s-1}\Delta[m-2]$ of $\Delta[m-1]$. Hence
\[
\mathcal K_{m,\ell}\cap d_0\Delta[m-1]
  =
\Lambda^{\ell-1}[m-1],
\]
and we have a pushout square
\[
\begin{tikzcd}[ampersand replacement=\&, sep=large]
\Lambda^{\ell-1}[m-1]\arrow[r,hook]\arrow[d,hook]
\&\mathcal K_{m,\ell}\arrow[d,hook]\\
\Delta[m-1]\arrow[r,hook,"\delta_0"']
\&\mathcal K_{m,\ell}\cup d_0\Delta[m-1]
  \arrow[ul,phantom,"\ulcorner",very near start].
\end{tikzcd}
\]
As $\ell-1<m-1$, the inclusion
$\mathcal K_{m,\ell}\hookrightarrow
\mathcal K_{m,\ell}\cup d_0\Delta[m-1]=\Lambda^\ell[m]$ belongs to
$\mathsf T$ by \cref{prop:left_horns_in_T} and
\cref{lem:T_closure}(1). Hence the composite
$\mathcal K_{m,\ell}\hookrightarrow\Lambda^\ell[m]
\hookrightarrow\Delta[m]$ belongs to $\mathsf T$ by
\cref{prop:left_horns_in_T} and \cref{lem:T_closure}(2).
\end{proof}

Finally, we show that the equivalences of $\mathbb E_m$ are detected
at the initial vertex.

\begin{lemma}
\label{lem:source_evaluation_reflects_equivalences}
Let $m\geq0$, and let $e$ be an edge of the quasi-category
$\mathbb E_m$. If the edge $u_m^*e$ is an equivalence in
$\mathbb E_0$, then $e$ is an equivalence in $\mathbb E_m$.
\end{lemma}

\begin{proof}
For $m=0$, the map $u_0^*$ is the identity, so suppose that
$m\geq1$. Since $(\bisimFin)_m$ is discrete, we have
$P_m(e)=s_0^v\beta$ for a unique $m$-simplex $\beta$ of
$\simFin$, and so $\lambda_m(e)=(u_m^*e,s_0^v\beta)$. The edge
$\lambda_m(e)$ is an equivalence in
$\mathbb E_0\times_{(\bisimFin)_0}(\bisimFin)_m$. Indeed, an inverse
of $u_m^*e$ in $\mathbb E_0$ and the $2$-simplices exhibiting it
as an inverse are sent to the vertex $u_m^*\beta$ of the discrete
simplicial set $(\bisimFin)_0$, and so, paired with the degenerate
simplices at $\beta$, they give an inverse of $\lambda_m(e)$. By
\cref{prop:minimal_axioms}, the map $\lambda_m$ is a weak
equivalence in $\sset_{(\infty,1)}$ between quasi-categories, and so
it induces an equivalence of homotopy categories by
\cite[Proposition~1.23]{Joyal2008Theory}. As the equivalences of a
quasi-category are the edges that become isomorphisms in its homotopy
category, the map $\lambda_m$ reflects equivalences, and hence $e$
is an equivalence in $\mathbb E_m$.
\end{proof}

\subsubsection{Horn filling}
\label{subsubsec:staircase_inner_horns}

In this section, we show that the horns $\Lambda^a[q]$, with
$0\leq a<q$, in the staircase $\Stair(\mathbb E)$ can be filled
componentwise, which reduces the lifting problems for the map
$\Stair(P)$ to lifting problems in $\mathbb E_0$. Given a
$q$-simplex $\beta$ of $\simFin$, we denote by
$c_\beta^{(j)}$ the degenerate $j$-simplex of the discrete
simplicial set $(\bisimFin)_q$ at $\beta$.

\begin{definition}
\label{def:relative_staircase_horn}
Let $q\geq2$, let $0\leq a<q$, and let $\sigma$ be a
$q$-simplex of $\simFin$. A \emph{relative $(q,a)$-horn over
$\sigma$} is a map $\eta:\Lambda^a[q]\to\Stair(\mathbb E)$ such
that $\Stair(P)\eta=\sigma|_{\Lambda^a[q]}$. Such a map $\eta$
amounts to a family of $(q-1)$-simplices of $\Stair(\mathbb E)$
\[
\eta^r=(\eta^r_0,\ldots,\eta^r_{q-1}),
\qquad
\eta^r_i\in\mathbb E_{q-1-i,i},
\qquad
0\leq r\leq q,\ r\ne a,
\]
such that $P_{q-1-i}(\eta^r_i)=c_{(d_0)^i d_r\sigma}^{(i)}$ and
$D_r\eta^s=D_{s-1}\eta^r$ for all $0\leq r<s\leq q$ with
$r,s\ne a$. By \cref{rem:staircase_faces_degeneracies}, the latter
condition is equivalent to the identities
\begin{equation}
\label{eq:relative_staircase_horn_overlap}
\begin{cases}
d_{r-i}^h\eta^s_i=d_{s-1-i}^h\eta^r_i,& i<r,\\
d_r^v\eta^s_{i+1}=d_{s-1-i}^h\eta^r_i,& r\leq i<s-1,\\
d_r^v\eta^s_{i+1}=d_{s-1}^v\eta^r_{i+1},& s-1\leq i\leq q-2.
\end{cases}
\end{equation}
\end{definition}

The following lemma constructs all the components of a filler of a
relative horn except the last one.

\begin{lemma}
\label{lem:staircase_extension}
Let $q\geq2$ and $0\leq a<q$, and let $\eta$ be a relative
$(q,a)$-horn over a $q$-simplex $\sigma$ of $\simFin$. There
exist simplices $\xi_j\in\mathbb E_{q-j,j}$, for $0\leq j<q$, such
that $P_{q-j}(\xi_j)=c_{\sigma_{\geq j}}^{(j)}$ and
\begin{numcases}{}
 d_s^h\xi_j=\eta^{j+s}_j
 & $(1\leq s\leq q-j,\ j+s\ne a),$
 \label{eq:IH2-horizontal}\\
 d_r^v\xi_j=\eta^r_{j-1}
 & $(0\leq r<j,\ r\ne a),$
 \label{eq:IH3-vertical}\\
 d_j^v\xi_j=d_0^h\xi_{j-1}
 & $(1\leq j<q).$
 \label{eq:IH3-staircase_gluing}
\end{numcases}
Moreover, for any such simplices, the assignments
\[
\Phi_q|_{d_r\Delta[q-1]}=
\begin{cases}
\eta^r_{q-1},&0\leq r<q,\ r\ne a,\\
d_0^h\xi_{q-1},&r=q
\end{cases}
\]
define a horn $\Phi_q:\Lambda^a[q]\to\mathbb E_0$ over the
degenerate simplex at $\sigma_{\geq q}$, and, for every extension
$\xi_q:\Delta[q]\to\mathbb E_0$ of $\Phi_q$, the tuple
$(\xi_0,\ldots,\xi_q)$ is a $q$-simplex of $\Stair(\mathbb E)$
over $\sigma$ that extends $\eta$.
\end{lemma}

The proof of \cref{lem:staircase_extension} is given in
\cref{app:staircase_extension}. It proceeds by induction on $j$ as
follows. For $j<q$, the faces $d_s^h\xi_j$ prescribed by
\eqref{eq:IH2-horizontal} assemble into a map out of the simplicial
set
\[
\mathcal K_j=
\begin{cases}
\mathcal K_{q-j,a-j}, & j<a,\\
\Lambda^0[q-j], & j\geq a,
\end{cases}
\]
with $\mathcal K_{m,\ell}$ as in
\cref{lem:two_missing_horizontal_faces}, and the faces
$d_r^v\xi_j$ prescribed by \eqref{eq:IH3-vertical} and
\eqref{eq:IH3-staircase_gluing} assemble into a map out of the
simplicial set
\[
\mathcal L_j=
\begin{cases}
\varnothing, & j=0,\\
\partial\Delta[j], & 0<j\leq a,\\
\Lambda^a[j], & j>a.
\end{cases}
\]
Hence the simplex $\xi_j$ is obtained as a lift in a commutative
square
\[
\begin{tikzcd}[sep=large]
  \mathcal L_j \ar[r] \ar[d,hook] &
  \mathbb E_{q-j} \ar[d,two heads,"\sim"'] \\
  \Delta[j] \ar[r] \ar[ru,dashed,"\xi_j"] &
  \mathbb E(\mathcal K_j)\times_{\bisimFin(\mathcal K_j)}(\bisimFin)_{q-j}
\end{tikzcd}
\]
whose right-hand vertical map is the relative restriction map along
the inclusion $\mathcal K_j\hookrightarrow\Delta[q-j]$, which is a
trivial fibration by
\cref{prop:left_horns_in_T,lem:two_missing_horizontal_faces}. The
compatibility of the prescribed faces and the commutativity of the
square follow from the identities
\eqref{eq:relative_staircase_horn_overlap} and from the relations
imposed at the previous stages.

\begin{ex}
\label{ex:staircase_horns_degree_two}
With the notation of \cref{ex:low_dimensional_staircases}, a relative
$(2,1)$-horn consists of the edges $\eta^2=(h_{01},v^1_{01})$ and
$\eta^0=(h'_{12},v^2_{12})$, and \cref{lem:staircase_extension}
produces the dashed part of the following diagram.
\[
\begin{tikzcd}[ampersand replacement=\&, sep=large]
\boxed{x_{00}}\arrow[r,"h_{01}"]\arrow[rr,bend left=30,dashed,"h_{02}"]
\&x_{10}\arrow[r,dashed,"h_{12}"]\arrow[d,"v^1_{01}"']
 \arrow[dr,phantom,"\xi_1" description]
\&x_{20}\arrow[d,dashed,"v^2_{01}"]
 \arrow[dd,bend left=50,dashed,"v^2_{02}"]\\
\&\boxed{x_{11}}\arrow[r,"h'_{12}"']\&x_{21}\arrow[d,"v^2_{12}"]\\
\&\&\boxed{x_{22}}
\end{tikzcd}
\]
For $j=0$, we have
$\mathcal K_0=\mathcal K_{2,1}=\Delta\{0,1\}$ and
$\mathcal L_0=\varnothing$, and the lift $\xi_0$ provides the edges
$h_{12}$ and $h_{02}$. For $j=1$, we have
$\mathcal K_1=\Lambda^0[1]$ and $\mathcal L_1=\partial\Delta[1]$,
and the lift $\xi_1$ is a square with right edge $v^2_{01}$.
Finally, the horn $\Phi_2:\Lambda^1[2]\to\mathbb E_0$ consists of the
edges $v^2_{01}$ and $v^2_{12}$, and a filler $\xi_2$ of
$\Phi_2$ provides the edge $v^2_{02}$.
\end{ex}

As a first consequence of \cref{lem:staircase_extension}, we obtain
the following.

\begin{theorem}
\label{thm:staircase_inner_fibration}
The map $\Stair(P):\Stair(\mathbb E)\to \simFin$ is an inner fibration.
\end{theorem}

\begin{proof}
Let $q\geq2$ and $0<a<q$, and let $\eta$ be a relative
$(q,a)$-horn over a $q$-simplex $\sigma$ of $\simFin$. By
\cref{lem:staircase_extension}, it is enough to show that the inner
horn $\Phi_q:\Lambda^a[q]\to\mathbb E_0$ admits a filler, which is
the case as $\mathbb E_0$ is a quasi-category by
\cref{prop:minimal_axioms}.
\end{proof}

\subsubsection{CoCartesian staircase edges}
\label{subsubsec:cocartesian_staircase_edges}

In this section, we show that an edge $(h,v)$ of
$\Stair(\mathbb E)$ is $\Stair(P)$-coCartesian if its vertical
component $v$ is an equivalence in $\mathbb E_0$, and we deduce
that $\Stair(P)$ is a coCartesian fibration.

\begin{proposition}
\label{prop:cocartesian_staircase_edges}
Let $e=(h,v)$ be an edge of $\Stair(\mathbb E)$. If $v$ is an
equivalence in the quasi-category $\mathbb E_0$, then $e$ is
$\Stair(P)$-coCartesian.
\end{proposition}

\begin{proof}
Since $\Stair(P)$ is an inner fibration by
\cref{thm:staircase_inner_fibration}, it is enough to show, by
\cref{lem:cocartesian_edge_horn_criterion}, that every commutative
square
\begin{equation}
\label{eq:cocart_edge_relative_horn}
\begin{tikzcd}[sep=large]
 \Lambda^0[q]\ar[r,"\eta"]\ar[d,hook] &
 \Stair(\mathbb E)\ar[d,"\Stair(P)"]\\
 \Delta[q]\ar[r,"\sigma"'] & \simFin
\end{tikzcd}
\end{equation}
with $q\geq2$ and $\eta|_{\Delta\{0,1\}}=e$ admits a lift. By
\cref{lem:staircase_extension} for $a=0$, we get simplices
$\xi_j\in\mathbb E_{q-j,j}$, for $0\leq j<q$, satisfying
\eqref{eq:IH2-horizontal}, \eqref{eq:IH3-vertical}, and
\eqref{eq:IH3-staircase_gluing}, and a horn
\[
\Phi_q:\Lambda^0[q]\longrightarrow\mathbb E_0,
\]
and it remains to show that $\Phi_q$ admits a filler. For this, we
show that the edge $\Phi_q|_{\Delta\{0,1\}}$ is an equivalence in
$\mathbb E_0$. By \eqref{eq:IH2-horizontal} for $j=1$ and
$s=q-1$, the edge $\xi_1$ of the quasi-category
$\mathbb E_{q-1}$ satisfies
\[
 u_{q-1}^*\xi_1
 =u_{q-2}^*d_{q-1}^h\xi_1
 =u_{q-2}^*\eta^q_1
 =v,
\]
where the last equality holds as the face $\eta^q$ contains the edge
$e$, whose vertical component is $u_{q-2}^*\eta^q_1$ by
\eqref{eq:staircase_alpha}. Hence $\xi_1$ is an equivalence in
$\mathbb E_{q-1}$ by
\cref{lem:source_evaluation_reflects_equivalences}. Moreover, since
the edge $\Delta\{0,1\}$ is contained in the last face of
$\Delta[j]$ for $2\leq j\leq q$, the relations
\eqref{eq:IH3-staircase_gluing} and the definition of $\Phi_q$ give
\begin{equation}
\label{eq:cocart_edge_final_initial_edge}
 \Phi_q|_{\Delta\{0,1\}}
 =d_0^h\bigl(\xi_{q-1}|_{\Delta\{0,1\}}\bigr)
 =(d_0^h)^{q-1}\xi_1.
\end{equation}
As the map $(d_0^h)^{q-1}:\mathbb E_{q-1}\to\mathbb E_0$ induces a
functor between homotopy categories
\cite[\S1.2.3]{lurie2009higher}, it preserves equivalences, and so the
edge $\Phi_q|_{\Delta\{0,1\}}$ is an equivalence in $\mathbb E_0$.
Hence the horn $\Phi_q$ admits a filler $\xi_q$ by
\cite[Proposition~1.2.4.3]{lurie2009higher}, and the tuple
$(\xi_0,\ldots,\xi_q)$ is the desired lift in
\eqref{eq:cocart_edge_relative_horn} by
\cref{lem:staircase_extension}.
\end{proof}

We refer to an edge of $\Stair(\mathbb E)$ of the form
$(h,s_0^v d_0^h h)$ as a \emph{pure horizontal staircase edge}. By
\cref{prop:cocartesian_staircase_edges}, every pure horizontal
staircase edge is $\Stair(P)$-coCartesian.

\begin{corollary}
\label{cor:staircase_is_cocartesian_fibration}
The map $\Stair(P):\Stair(\mathbb E)\rightarrow \simFin$ is a
coCartesian fibration. In particular, the simplicial set
$\Stair(\mathbb E)$ is a quasi-category.
\end{corollary}

\begin{proof}
By \cref{thm:staircase_inner_fibration}, the map $\Stair(P)$ is an
inner fibration. Let $x$ be a vertex of $\Stair(\mathbb E)$, i.e.,
an element of $\mathbb E_{0,0}$, and let $f:P_0(x)\to b$ be an
edge of $\simFin$. Since $\lambda_1$ is a trivial fibration by
\cref{prop:minimal_axioms}, there is an element
$h\in\mathbb E_{1,0}$ such that $d_1^hh=x$ and $P_1(h)=f$. Then
the pure horizontal staircase edge $(h,s_0^vd_0^hh)$ has source
$x$, lies over $f$, and is $\Stair(P)$-coCartesian by
\cref{prop:cocartesian_staircase_edges}. The last assertion holds as
$\simFin$ is a quasi-category.
\end{proof}

\subsubsection{Fibers and inert spans}
\label{subsubsec:staircase_fibers_spans}

In this section, we compare the fibers and the spans
\eqref{eq:marked_inert_span} of the coCartesian fibration
$\Stair(P)$ with the fibers and the spans
\eqref{eq:segality_double_span} of the left double
$(\infty,1)$-fibration $P$, and we deduce that $\Stair(P)$ is
special if and only if $P$ is special. Throughout, we write
$\mathcal X=\Stair(\mathbb E)$, and
$\mathcal X_a=\mathcal X\times_{\simFin}\Delta[0]$ for the fiber of
$\Stair(P)$ at an object $a$ of $\Fin$. Recall that the
simplicial sets $\mathbb E_m(\sigma)$ and $\mathcal X_a$ are
quasi-categories, by \cref{lem:fixed_shape_fiber_fibrant} and since
$\Stair(P)$ is an inner fibration.

\begin{lemma}
\label{lem:staircase_fiber_comparison}
For every object $a$ of $\Fin$, there is a map
\[
 \iota_a:\mathbb E_0(a)\longrightarrow\mathcal X_a,
\]
natural in the left double $(\infty,1)$-fibration $P$, which is a
monomorphism and a weak equivalence in $\sset_{(\infty,1)}$. It is
given on a $q$-simplex $x$ of $\mathbb E_0(a)$ by
\begin{equation}
\label{eq:staircase_fiber_inclusion}
 \bigl(\iota_a(x)\bigr)_j
 =(c_{q-j})_h^*\bigl(x|_{[0,j]}\bigr),
 \qquad 0\leq j\leq q,
\end{equation}
where $c_m:[m]\to[0]$ denotes the unique map in $\Delta$.
\end{lemma}

\begin{proof}
Let $\id_a^m$ denote the degenerate $m$-simplex of $\simFin$ at
$a$, and consider the bisimplicial set
$\mathbb E\times_{\bisimFin}(\Delta[0])^h$, where the map
$(\Delta[0])^h\to\bisimFin$ is induced by $a$, whose $m$-th level
is $\mathbb E_m(\id_a^m)$. The maps
\[
 (c_m)_h^*:\mathbb E_0(a)\longrightarrow\mathbb E_m(\id_a^m)
\]
are compatible with the horizontal structure maps, as
$c_m\alpha=c_\ell$ for every map $\alpha:[\ell]\to[m]$ in
$\Delta$, and so they assemble into a map of bisimplicial sets
\[
 \upsilon_a:\Delta[0]\boxtimes\mathbb E_0(a)
 \longrightarrow
 \mathbb E\times_{\bisimFin}(\Delta[0])^h.
\]
Since $\lambda_m^{\id_a^m}(c_m)_h^*=\id_{\mathbb E_0(a)}$ and the
map $\lambda_m^{\id_a^m}$ is a trivial fibration by
\cref{cor:fixed_simplex_source}, the map $(c_m)_h^*$ is a weak
equivalence in $\sset_{(\infty,1)}$ for every $m\geq0$, i.e., the
map $\upsilon_a$ is a levelwise weak equivalence. We claim that the
map $\Stair(\upsilon_a)=T(\upsilon_a^{\mathrm{tr}})$ is then a weak
equivalence in $\sset_{(\infty,1)}$. Indeed, the diagonal functor
$\mathbb X\mapsto(\mathbb X_{q,q})_{q\geq0}$ from bisimplicial sets
to simplicial sets is left Quillen from the Reedy model structure to
$\sset_{(\infty,1)}$, as the relative matching maps of its right
adjoint $\mathcal Y\mapsto(\mathcal Y^{\Delta[m]})_{m\geq0}$ are the
pullback-exponentials with the inclusions
$\partial\Delta[m]\hookrightarrow\Delta[m]$ and
$\sset_{(\infty,1)}$ is cartesian. Since every bisimplicial set is
Reedy cofibrant by \cref{lem:presheaf_reedy_cofibrancy}, the diagonal
functor preserves levelwise weak equivalences, in either variable by
symmetry, and the claim follows from the natural weak equivalence in
$\sset_{(\infty,1)}$ from the diagonal to the total simplicial set
of \cite[Theorem~1.2]{Stevenson2025TotalSimplicialSet}; see also
\cite[Proposition~7.1]{Stevenson2025TotalSimplicialSet}. Finally, we
have isomorphisms of simplicial sets
\[
\begin{aligned}
 \Stair\bigl(\Delta[0]\boxtimes\mathbb E_0(a)\bigr)
 &\cong\mathbb E_0(a),\\
 \Stair\bigl(\mathbb E\times_{\bisimFin}(\Delta[0])^h\bigr)
 &\cong\mathcal X_a,
\end{aligned}
\]
where the first one is given by \cref{lem:staircase_external_product},
and the second one holds as a $q$-simplex of $\mathcal X$ lies over
$\id_a^q$ if and only if its $j$-th component lies over
$\id_a^{q-j}$ for every $0\leq j\leq q$. We define $\iota_a$ to
be the map $\Stair(\upsilon_a)$ under these isomorphisms. By
construction, the map $\iota_a$ is given by
\eqref{eq:staircase_fiber_inclusion} and is natural in $P$, and it
is a monomorphism, as the last component of $\iota_a(x)$ is $x$.
\end{proof}

\begin{proposition}
\label{prop:staircase_transport_comparison}
For every edge $f:a\to b$ of $\simFin$, there is a weak equivalence
in $\sset_{(\infty,1)}$
\[
 \chi_f:\mathbb E_1(f)\longrightarrow\Cocart_{\Stair(P)}(f)
\]
making the following diagram of simplicial sets commute.
\begin{equation}
\label{eq:staircase_transport_diagram}
\begin{tikzcd}[ampersand replacement=\&, sep=large]
\mathbb E_0(a)\arrow[d,hook,"\iota_a"',"\sim"]
\&\mathbb E_1(f)\arrow[l,two heads,"s_f^P"',"\sim"]
 \arrow[r,"t_f^P"]\arrow[d,"\chi_f","\sim"']
\&\mathbb E_0(b)\arrow[d,hook,"\iota_b","\sim"']\\
\mathcal X_a
\&\Cocart_{\Stair(P)}(f)
 \arrow[l,two heads,"s_f^{\Stair(P)}","\sim"']
 \arrow[r,"t_f^{\Stair(P)}"']
\&\mathcal X_b.
\end{tikzcd}
\end{equation}
\end{proposition}

\begin{proof}
An $r$-simplex $h$ of $\mathbb E_1(f)$ is a map of bisimplicial
sets $h:\Delta[1]\boxtimes\Delta[r]\to\mathbb E$ over $f$. We
denote by $F_h$ the composite
\[
\begin{tikzcd}[ampersand replacement=\&,column sep=small]
\Delta[r]\times\Delta[1]\arrow[r,"\mathrm{tw}"]
\&\Delta[1]\times\Delta[r]\arrow[r,"\cong"]
\&\Stair(\Delta[1]\boxtimes\Delta[r])
 \arrow[r,"\Stair(h)"]\&\mathcal X,
\end{tikzcd}
\]
where $\mathrm{tw}$ interchanges the two factors and the isomorphism
is that of \cref{lem:staircase_external_product}. By naturality of
\eqref{eq:staircase_external_product}, we have
$\Stair(P)F_h=f\pr_2$. Moreover, for every vertex $z$ of
$\Delta[r]$, the restriction of $F_h$ to
$\Delta\{z\}\times\Delta[1]$ is the pure horizontal staircase edge
$(h_z,s_0^v d_0^h h_z)$, where $h_z$ denotes the restriction of
$h$ to $\Delta[1]\boxtimes\Delta\{z\}$, and this edge is
$\Stair(P)$-coCartesian by \cref{prop:cocartesian_staircase_edges}.
Hence $F_h$ is an $r$-simplex of $\Cocart_{\Stair(P)}(f)$ by
\eqref{eq:marked_flat_mapping}, and, as the construction is natural in
$\Delta[r]$, we get a map of simplicial sets $\chi_f$ sending $h$
to $F_h$. By \eqref{eq:staircase_fiber_inclusion}, the restrictions
of $F_h$ to $\Delta[r]\times\{0\}$ and $\Delta[r]\times\{1\}$ are
$\iota_a(s_f^P h)$ and $\iota_b(t_f^P h)$, respectively, so that
\begin{equation}
\label{eq:staircase_transport_identities}
 s_f^{\Stair(P)}\chi_f=\iota_a s_f^P,
 \qquad
 t_f^{\Stair(P)}\chi_f=\iota_b t_f^P.
\end{equation}
i.e., the diagram \eqref{eq:staircase_transport_diagram} commutes.
Finally, since the maps $s_f^P$ and $s_f^{\Stair(P)}$ are trivial
fibrations by \cref{prop:cleavage_free_segal,lem:cocart_transport_span}
and $\iota_a$ is a weak equivalence by
\cref{lem:staircase_fiber_comparison}, the map $\chi_f$ is a weak
equivalence, by the two-out-of-three property.
\end{proof}

\begin{theorem}
\label{thm:staircase_preserves_reflects_specialness}
Let $P:\mathbb E\to\bisimFin$ be a left double
$(\infty,1)$-fibration, and write $\mathcal X=\Stair(\mathbb E)$.
For every $k\geq0$, there is a map
\[
\Xi_k:\mathbb E(\mathcal J_k)_{\beta_k}
 \longrightarrow
 \Map^\flat((\mathcal J_k)^\sharp_{\beta_k},\mathcal X^\natural)
\]
making the following diagrams of simplicial sets commute, in which the
vertical maps are weak equivalences in $\sset_{(\infty,1)}$.
\begin{equation}
\label{eq:staircase_specialness_diagram}
\begin{gathered}
\begin{tikzcd}[ampersand replacement=\&, sep=large]
\mathbb E(\mathcal J_k)_{\beta_k}
 \arrow[r,two heads,"\pi_k^P","\sim"']
 \arrow[d,"\Xi_k"',"\sim"]
\&\mathbb E_0(\langle k\rangle)
 \arrow[d,hook,"\iota_{\langle k\rangle}","\sim"']\\
\Map^\flat((\mathcal J_k)^\sharp_{\beta_k},\mathcal X^\natural)
 \arrow[r,two heads,"\pi_k^{\Stair(P)}"',"\sim"]
\&\mathcal X_{\langle k\rangle},
\end{tikzcd}
\qquad
\begin{tikzcd}[ampersand replacement=\&, sep=large]
\mathbb E(\mathcal J_k)_{\beta_k}
 \arrow[r,two heads,"t_k^P"]\arrow[d,"\Xi_k"',"\sim"]
\&\mathbb E_0(\langle1\rangle)^k
 \arrow[d,hook,"(\iota_{\langle1\rangle})^{\times k}","\sim"']\\
\Map^\flat((\mathcal J_k)^\sharp_{\beta_k},\mathcal X^\natural)
 \arrow[r,two heads,"t_k^{\Stair(P)}"']
\&\mathcal X_{\langle1\rangle}^{\,k}.
\end{tikzcd}
\end{gathered}
\end{equation}
In particular, the left double $(\infty,1)$-fibration $P$ is
special if and only if the coCartesian fibration $\Stair(P)$ is
special.
\end{theorem}

\begin{proof}
For $k\geq1$, by
\cref{lem:internal_fan_identification,lem:inert_fan_mapping_objects},
the sources of the maps $\pi_k^P$ and $\pi_k^{\Stair(P)}$ are the
iterated pullbacks of the maps $s_{\rho_i}^P$ and
$s_{\rho_i}^{\Stair(P)}$, respectively. By
\eqref{eq:staircase_transport_identities}, we have
\[
 s_{\rho_i}^{\Stair(P)}\chi_{\rho_i}\pr_i
 =\iota_{\langle k\rangle}s_{\rho_i}^P\pr_i
 =\iota_{\langle k\rangle}\pi_k^P
\]
for every $1\leq i\leq k$, and so the maps $\chi_{\rho_i}\pr_i$
induce a map $\Xi_k$ such that
\begin{align}
 \pi_k^{\Stair(P)}\Xi_k
 &=\iota_{\langle k\rangle}\pi_k^P,
 \label{eq:staircase_fan_source_identity}\\
 t_k^{\Stair(P)}\Xi_k
 &=(\iota_{\langle1\rangle})^{\times k}t_k^P.
 \label{eq:staircase_fan_target_identity}
\end{align}
i.e., both squares of \eqref{eq:staircase_specialness_diagram}
commute. Since the maps $\pi_k^P$ and $\pi_k^{\Stair(P)}$ are
trivial fibrations by
\cref{prop:cleavage_free_segal,prop:marked_inert_span_fibrations} and
$\iota_{\langle k\rangle}$ is a weak equivalence by
\cref{lem:staircase_fiber_comparison}, the map $\Xi_k$ is a weak
equivalence, by \eqref{eq:staircase_fan_source_identity} and the
two-out-of-three property. Moreover, the map
$(\iota_{\langle1\rangle})^{\times k}$ is a weak equivalence by
\cref{prop:model_category_recollections}(5), as a finite product of
weak equivalences between quasi-categories. Hence, by
\eqref{eq:staircase_fan_target_identity} and the two-out-of-three
property, the map $t_k^P$ is a weak equivalence if and only if
$t_k^{\Stair(P)}$ is one. For $k=0$, we set
$\Xi_0=\iota_{\langle0\rangle}$, and the same argument applies, as
$t_0^P=t_0^{\Stair(P)}\iota_{\langle0\rangle}$.
\end{proof}

\subsection{From coCartesian fibrations to left double fibrations}
\label{subsec:comparison_enriched_to_internal}
\subsubsection{CoCartesian strings}

\label{subsubsec:horizontal_construction}

In this section, we associate to a marked simplicial set
$\mathcal X$ over $\simFin^\sharp$ a bisimplicial set
$\Int(\mathcal X)$ over $\bisimFin$, whose $m$-th level is the
simplicial set of marked $m$-simplices of $\mathcal X$. The functor
$\Int$ will be the right adjoint in
\cref{thm:left_double_cocart_quillen_equivalence}. We first show that
the inclusion of the spine into a maximally marked simplex is marked
anodyne, which will imply the Segal condition for $\Int(\mathcal X)$.

\begin{lemma}
\label{lem:marked-spine-anodyne}
For every $m\geq0$, the inclusion
$I[m]^\sharp\hookrightarrow\Delta[m]^\sharp$ is both marked
anodyne and coCartesian marked anodyne.
\end{lemma}

\begin{proof}
For $m\leq1$, the inclusion is an identity, so suppose that
$m\geq2$. For $1\leq\ell\leq m$, let $\mathcal D^{(\ell)}$ denote
the marked simplicial set with underlying simplicial set $\Delta[m]$
in which the marked edges are the degenerate edges and the edges
$[i,j]$ with $j-i\leq\ell$, so that
$\mathcal D^{(m)}=\Delta[m]^\sharp$. Since the inclusion
$I[m]\hookrightarrow\Delta[m]$ is inner anodyne by
\cite[Proposition~2.13]{Joyal2008Theory}, the inclusion
$I[m]^\flat\hookrightarrow\Delta[m]^\flat$ is marked anodyne by
\cite[Remark~3.1.1.4]{lurie2009higher}, and hence so is its pushout
$I[m]^\sharp\hookrightarrow\mathcal D^{(1)}$ in the following
square.
\[
\begin{tikzcd}[ampersand replacement=\&, sep=large]
I[m]^\flat\arrow[r,hook]\arrow[d,hook]
\&\Delta[m]^\flat\arrow[d,hook]\\
I[m]^\sharp\arrow[r,hook]
\&\mathcal D^{(1)}\arrow[ul,phantom,"\ulcorner",very near start].
\end{tikzcd}
\]
Next, let
$\mathcal M=(\Lambda^1[2])^\sharp
  \amalg_{(\Lambda^1[2])^\flat}\Delta[2]^\flat$, so that the inclusion
$g:\mathcal M\hookrightarrow\Delta[2]^\sharp$ is the generating
marked anodyne map of
\cite[Definition~3.1.1.1\textup{(3)}]{lurie2009higher}. For
$2\leq\ell\leq m$ and $0\leq i<j\leq m$ with $j-i=\ell$, the
$2$-simplex of $\Delta[m]$ with vertices $i,j-1,j$ defines a map
of marked simplicial sets
$\theta_{i,j}:\mathcal M\to \mathcal D^{(\ell-1)}$, as its edges
$[i,j-1]$ and $[j-1,j]$ are marked in $\mathcal D^{(\ell-1)}$,
and we have a pushout square
\[
\begin{tikzcd}[ampersand replacement=\&, sep=large]
\displaystyle\coprod_{j-i=\ell}\mathcal M
  \arrow[r,"{\coprod\theta_{i,j}}"]
  \arrow[d,hook,"\coprod g"']
\&\mathcal D^{(\ell-1)}\arrow[d,hook]\\
\displaystyle\coprod_{j-i=\ell}\Delta[2]^\sharp\arrow[r]
\&\mathcal D^{(\ell)}\arrow[ul,phantom,"\ulcorner",very near start].
\end{tikzcd}
\]
Hence the inclusion
$\mathcal D^{(\ell-1)}\hookrightarrow \mathcal D^{(\ell)}$ is marked
anodyne, and so is the composite
\[
I[m]^\sharp\hookrightarrow \mathcal D^{(1)}\hookrightarrow\cdots
\hookrightarrow \mathcal D^{(m)}=\Delta[m]^\sharp.
\]
Finally, reversing the order of the vertices identifies the opposite
of the inclusion $I[m]^\sharp\hookrightarrow\Delta[m]^\sharp$ with
itself, and so this inclusion is also coCartesian marked anodyne.
\end{proof}

We now define the functor $\Int$.

\begin{construction}
\label{const:cocartesian_strings}
Given a marked simplicial set $\mathcal X$ over $\simFin^\sharp$,
we define the bisimplicial set $\Int(\mathcal X)$ over $\bisimFin$
by
\[
\Int(\mathcal X)_m=
\coprod_{\sigma\in (\simFin)_m}
\Map^\flat
\bigl(\Delta[m]^\sharp_\sigma,\mathcal X\bigr),
\]
where the map $\Int(\mathcal X)_m\to(\bisimFin)_m$ sends the summand
indexed by $\sigma$ to the point $\sigma$ of $(\bisimFin)_m$. A
map $\alpha:[\ell]\to[m]$ in $\Delta$ induces maps of marked
simplicial sets
$\alpha^\sharp_\sigma:\Delta[\ell]^\sharp_{\sigma\alpha}\to
\Delta[m]^\sharp_\sigma$, and the map
$\alpha_h^*:\Int(\mathcal X)_m\to\Int(\mathcal X)_\ell$ is given on
the summand indexed by $\sigma$ by precomposition with
$\alpha^\sharp_\sigma$. The construction is functorial in
$\mathcal X$ by postcomposition, and so it defines a functor
\[
 \Int:(\sset^+)_{/\simFin^\sharp}
   \longrightarrow (\sset^{\Delta^{\mathrm{op}}})_{/\bisimFin}.
\]
For a coCartesian fibration $p:\mathcal E\to \simFin$, we denote the
map $\Int(\mathcal E^\natural)\to\bisimFin$ by $\Int(p)$, and we
refer to $\Int(\mathcal E^\natural)$ as the bisimplicial set of
\emph{coCartesian strings} of $p$.
\end{construction}

Note that an $r$-simplex of $\Int(\mathcal E^\natural)_m$ in the
summand indexed by $\sigma$ is a map
$F:\Delta[r]\times\Delta[m]\to\mathcal E$ over $\sigma\pr_2$ whose
restriction to $\{a\}\times\Delta[m]$ is, for each vertex $a$ of
$\Delta[r]$, a string of $m$ composable $p$-coCartesian edges
together with their composites, which are $p$-coCartesian by the
dual of \cite[Proposition~2.4.1.7]{lurie2009higher}.

\begin{ex}
\label{ex:monoidal_category_strings}
Let $\mathrm C$ be a small symmetric monoidal category, and let
$p:\mathcal E=N(\mathrm C^\otimes)\to\simFin$ be the coCartesian
fibration of \cref{ex:ordinary_symmetric_monoidal}. Then
$\Int(\mathcal E^\natural)_0(\langle k\rangle)=N(\mathrm C^k)$, and,
for a map $f:\langle k\rangle\to\langle\ell\rangle$ in $\Fin$, the
simplicial set $\Int(\mathcal E^\natural)_1(f)$ is the nerve of the
category whose objects are the triples $(x,y,\theta)$ of tuples
$x=(x_i)_{i=1}^k$ and $y=(y_j)_{j=1}^\ell$ of objects of
$\mathrm C$ and of isomorphisms
\[
 \theta_j:\bigotimes_{i\in f^{-1}(j)}x_i\xrightarrow{\cong} y_j
 \qquad(1\leq j\leq\ell),
\]
and whose morphisms $(x,y,\theta)\to(x',y',\theta')$ are the tuples
of morphisms $b_i:x_i\to x'_i$ and $c_j:y_j\to y'_j$ in
$\mathrm C$ such that
$c_j\theta_j=\theta'_j\bigl(\bigotimes_{i\in f^{-1}(j)}b_i\bigr)$
for every $1\leq j\leq\ell$.
\end{ex}

\subsubsection{The shape formula}

\label{subsubsec:strings_auxiliary}

We now compute the evaluation of the bisimplicial set
$\Int(\mathcal X)$ at a simplicial set $\mathcal K$.

\begin{proposition}
\label{prop:matching_identification}
Let $\mathcal K$ be a simplicial set and let $\mathcal X$ be a
marked simplicial set over $\simFin^\sharp$. There is an isomorphism
\begin{equation}
\label{eq:strings_shape_identification}
 \Int(\mathcal X)(\mathcal K)\cong\coprod_{\beta:\mathcal K\to \simFin}
 \Map^\flat(\mathcal K^\sharp_\beta,\mathcal X),
\end{equation}
natural in $\mathcal K$ and $\mathcal X$, under which the map to
$\bisimFin(\mathcal K)=\sset(\mathcal K,\simFin)^{\mathrm{disc}}$ is
the projection onto the indexing set. In particular, for every
monomorphism $i:\mathcal K\hookrightarrow\Delta[m]$, we have
\begin{equation}
\label{eq:strings_relative_shape_identification}
 \Int(\mathcal X)(\mathcal K)\times_{\bisimFin(\mathcal K)}(\bisimFin)_m
 \cong\coprod_{\sigma\in (\simFin)_m}
 \Map^\flat(\mathcal K^\sharp_{\sigma i},\mathcal X),
\end{equation}
and the relative restriction map $R_i$ of $\Int(\mathcal X)$ is the
coproduct over $\sigma\in (\simFin)_m$ of the restriction maps along
the inclusions
$
 \mathcal K^\sharp_{\sigma i}\hookrightarrow\Delta[m]^\sharp_\sigma.
$
\end{proposition}

\begin{proof}
The simplicial set $\mathcal K$ is the colimit of its simplices, and
we have an isomorphism of marked simplicial sets
\[
 \mathcal K^\sharp\cong
 \operatorname*{colim}_{(z:\Delta[\ell]\to\mathcal K)\in\Delta/\mathcal K}
 \Delta[\ell]^\sharp,
\]
as every edge of $\mathcal K$ is the image of an edge of a simplex.
By \eqref{eq:global_shape_mapping}, an $r$-simplex of
$\Int(\mathcal X)(\mathcal K)$ is a family of pairs
$(\sigma_z,F_z)$, indexed by the simplices
$z:\Delta[\ell]\to\mathcal K$ of $\mathcal K$, of a simplex
$\sigma_z:\Delta[\ell]\to\simFin$ and a map of marked simplicial
sets
\[
 F_z:\Delta[r]^\flat\times
 \Delta[\ell]^\sharp_{\sigma_z}\longrightarrow \mathcal X
\]
over $\simFin^\sharp$, such that $\sigma_z=\sigma_w\theta$ and
$F_z=F_w(\id\times\theta)$ whenever $z=w\theta$. The simplices
$\sigma_z$ amount to a map $\beta:\mathcal K\to \simFin$ with
$\beta z=\sigma_z$ for every $z$. Moreover, since the category of
marked simplicial sets is cartesian closed
\cite[\S3.1.3]{lurie2009higher}, the product with $\Delta[r]^\flat$
preserves colimits, and so we have an isomorphism
\[
 \operatorname*{colim}_{z\in\Delta/\mathcal K}
 \bigl(\Delta[r]^\flat\times
 \Delta[\ell]^\sharp_{\beta z}\bigr)
 \cong\Delta[r]^\flat\times \mathcal K^\sharp_\beta
\]
over $\simFin^\sharp$. Hence the maps $F_z$ amount to a map of
marked simplicial sets
$F:\Delta[r]^\flat\times \mathcal K^\sharp_\beta\to \mathcal X$ over
$\simFin^\sharp$, i.e., to an $r$-simplex of
$\Map^\flat(\mathcal K^\sharp_\beta,\mathcal X)$. This correspondence
is compatible with the simplicial structure in $r$ and natural in
$\mathcal K$ and $\mathcal X$, which gives the isomorphism
\eqref{eq:strings_shape_identification}, and the map to
$\bisimFin(\mathcal K)$ sends a pair $(\beta,F)$ to $\beta$. The
isomorphism \eqref{eq:strings_relative_shape_identification} follows by
taking the pullback along the map
$(\bisimFin)_m\to\bisimFin(\mathcal K)$ sending $\sigma$ to
$\sigma i$, and the description of $R_i$ follows from the
naturality in $\mathcal K$.
\end{proof}

\subsubsection{Reedy fibrancy and the horizontal Segal condition}

\label{subsubsec:relative_reedy}

In the rest of \cref{subsec:comparison_enriched_to_internal}, we fix a
coCartesian fibration $p:\mathcal E\to \simFin$. We first show that
the map $\Int(p)$ is
a Reedy fibration and that the bisimplicial set
$\Int(\mathcal E^\natural)$ satisfies the Segal condition.

\begin{lemma}
\label{lem:strings_reedy_fibration}
The map $\Int(p):\Int(\mathcal E^\natural)\to\bisimFin$ is a
Reedy fibration in $(\sset_{(\infty,1)})^{\Delta^{\mathrm{op}}}$.
\end{lemma}

\begin{proof}
By \cref{prop:matching_identification} applied to the inclusion
$\partial\Delta[m]\hookrightarrow\Delta[m]$, the relative matching
map of $\Int(p)$ in degree $m$ is the coproduct over
$\sigma\in(\simFin)_m$ of the restriction maps along the inclusions
\[
 (\partial\Delta[m])^\sharp_{\partial\sigma}
   \hookrightarrow\Delta[m]^\sharp_\sigma.
\]
These restriction maps are fibrations in $\sset_{(\infty,1)}$ by
\cref{prop:marked_mapping_fibration}, as $\mathcal E^\natural$ is
fibrant in $\cocart$, and hence so is their coproduct by
\cref{lem:fixed_index_coproducts}.
\end{proof}

Next, since the map $(\simFin)_m\to\sset(I[m],\simFin)$ given by
restriction along the inclusion $I[m]\hookrightarrow\Delta[m]$ is a
bijection, \cref{prop:matching_identification} gives an isomorphism
\begin{equation}
\label{eq:segal-target-as-spines}
\Int(\mathcal E^\natural)(I[m])
\cong\coprod_{\sigma\in (\simFin)_m}
\Map^\flat
\bigl(I[m]^\sharp_{\sigma|I[m]},\mathcal E^\natural\bigr),
\end{equation}
under which the Segal map $\varphi_m$ of
$\Int(\mathcal E^\natural)$ is the coproduct over
$\sigma\in (\simFin)_m$ of the restriction maps
\[
R_{I[m],\sigma}:
\Map^\flat
\bigl(\Delta[m]^\sharp_\sigma,\mathcal E^\natural\bigr)
\longrightarrow
\Map^\flat
\bigl(I[m]^\sharp_{\sigma|I[m]},\mathcal E^\natural\bigr).
\]

\begin{proposition}
\label{prop:marked-direct}
For every $m\geq2$, the Segal map $\varphi_m$ of
$\Int(\mathcal E^\natural)$ is a trivial fibration of simplicial
sets.
\end{proposition}

\begin{proof}
By \cref{lem:marked-spine-anodyne,conv:cocart-anodyne}, the inclusion
$I[m]^\sharp_{\sigma|I[m]}\hookrightarrow\Delta[m]^\sharp_\sigma$ is
a trivial cofibration in $\cocart$ for every
$\sigma\in(\simFin)_m$, and so the map $R_{I[m],\sigma}$ is a
trivial fibration by \cref{prop:marked_mapping_fibration}. Hence the
Segal map $\varphi_m$ is a trivial fibration by
\eqref{eq:segal-target-as-spines} and
\cref{lem:fixed_index_coproducts}.
\end{proof}

\subsubsection{Initial-vertex evaluation and left double fibrancy}

\label{subsubsec:initial_vertex_condition}

It remains to show that the map $\lambda_1$ of $\Int(p)$ is a weak
equivalence. Note that, by \cref{const:cocartesian_strings}, we have
\[
\Int(\mathcal E^\natural)_1(f)=\Cocart_p(f)
\qquad\text{and}\qquad
\Int(\mathcal E^\natural)_0(a)=\mathcal E_a
\]
for every edge $f:a\to b$ of $\simFin$, and that the map
$s_f^{\Int(p)}$ is the map $s_f^p$ of
\eqref{eq:single_cocartesian_span}.

\begin{lemma}
\label{lem:strings_source_evaluation}
The map
\[
\lambda_1:\Int(\mathcal E^\natural)_1\longrightarrow
\Int(\mathcal E^\natural)_0\times_{(\bisimFin)_0}(\bisimFin)_1
\]
is a trivial fibration of simplicial sets.
\end{lemma}

\begin{proof}
By \cref{prop:matching_identification} applied to the inclusion
$\Delta\{0\}\hookrightarrow\Delta[1]$, the map $\lambda_1$ is the
coproduct over the edges $f$ of $\simFin$ of the maps $s_f^p$,
which are trivial fibrations by \cref{lem:cocart_transport_span}.
Hence $\lambda_1$ is a trivial fibration by
\cref{lem:fixed_index_coproducts}.
\end{proof}

By putting the above results together, we get:

\begin{corollary}
\label{cor:left-double}
For every coCartesian fibration $p:\mathcal E\to \simFin$, the map
\[
\Int(p):\Int(\mathcal E^\natural)\longrightarrow\bisimFin
\]
is a left double $(\infty,1)$-fibration. In particular, the map
$\lambda_m$ of $\Int(p)$ is a trivial fibration of simplicial sets
for every $m\geq0$.
\end{corollary}

\begin{proof}
The first assertion follows directly from
\cref{lem:strings_reedy_fibration,prop:marked-direct,lem:strings_source_evaluation},
and the second one from \cref{prop:minimal_axioms}.
\end{proof}

\subsubsection{Compatibility with specialness}

Finally, we show that the spans \eqref{eq:segality_double_span} of the
left double $(\infty,1)$-fibration $\Int(p)$ are isomorphic to the
spans \eqref{eq:marked_inert_span} of $p$.

\begin{corollary}
\label{cor:strings_specialness}
Let $p:\mathcal E\to \simFin$ be a coCartesian fibration. For every
$k\geq0$, there are isomorphisms
\[
 \Int(\mathcal E^\natural)(\mathcal J_k)_{\beta_k}
 \cong\Map^\flat((\mathcal J_k)^\sharp_{\beta_k},\mathcal E^\natural),
 \qquad
 \Int(\mathcal E^\natural)(\mathcal V_k)_{\beta_k|\mathcal V_k}
 \cong\mathcal E_{\langle1\rangle}^{\,k},
\]
natural in $p$, under which the maps $\pi_k^{\Int(p)}$ and
$t_k^{\Int(p)}$ correspond to the maps $\pi_k^p$ and $t_k^p$,
respectively. In particular, the coCartesian fibration $p$ is
special if and only if the left double $(\infty,1)$-fibration
$\Int(p)$ is special.
\end{corollary}

\begin{proof}
The isomorphisms are given by
\eqref{eq:strings_shape_identification} applied to the simplicial sets
$\mathcal J_k$ and $\mathcal V_k$, and by
\cref{lem:inert_fan_mapping_objects}. Since the isomorphism
\eqref{eq:strings_shape_identification} is natural in $\mathcal K$,
the second assertion follows from
\cref{lem:internal_fan_identification,lem:inert_fan_mapping_objects},
as the maps $\pi_k$ and $t_k$ are, in both cases, the restriction
maps along the inclusions $\Delta\{c\}\hookrightarrow\mathcal J_k$
and $\mathcal V_k\hookrightarrow\mathcal J_k$.
\end{proof}

In particular, by
\cref{thm:staircase_preserves_reflects_specialness,prop:special_cocartesian_HA},
the staircase of a special left double $(\infty,1)$-fibration is a
symmetric monoidal $\infty$-category in the sense of
\cite[Definition~2.0.0.7]{LurieHA}.

\subsection{The marked diagonal and the Quillen adjunction}
\label{subsec:marked_diagonal_adjunction}

In this section, we construct a left adjoint $\Diag$ of the functor
$\Int$ of \cref{const:cocartesian_strings}, and we show that the
adjunction $\Diag\dashv\Int$ is a Quillen adjunction between the
model structures $\ldblfib$ and $\cocart$.

\begin{construction}
\label{const:marked_diagonal}
Given a bisimplicial set $P:\mathbb E\to\bisimFin$ over
$\bisimFin$, we define the \emph{marked diagonal} of $\mathbb E$ to
be the marked simplicial set
\[
 \Diag(\mathbb E)=
 \int^{[m]\in\Delta}\mathbb E_m^\flat\times\Delta[m]^\sharp,
\]
with structure map to $\simFin^\sharp$ induced by the maps
$\sigma\operatorname{pr}_2$ on the summands
$\mathbb E_m(\sigma)^\flat\times\Delta[m]^\sharp$, for
$\sigma\in(\simFin)_m$. These maps are compatible with the
identifications of the coend, since, for every map
$\theta:[\ell]\to[m]$ in $\Delta$, the images of
$
 \mathbb E_m(\sigma)^\flat
 \times\Delta[\ell]^\sharp_{\sigma\theta}
$
under the maps $(\theta_{\mathbb E}^*)^\flat\times\id$ and
$\id\times\theta^\sharp$ both lie over
$\sigma\theta\operatorname{pr}_2$.
\end{construction}

We first describe the marked simplicial set $\Diag(\mathbb E)$
explicitly.

\begin{proposition}
\label{prop:marked_diagonal_description}
Let $P:\mathbb E\to\bisimFin$ be a bisimplicial set over
$\bisimFin$. The underlying simplicial set of $\Diag(\mathbb E)$ is
the diagonal $\diagset\mathbb E$ of $\mathbb E$, given by
$(\diagset\mathbb E)_q=\mathbb E_{q,q}$, and its marked edges are the
edges $s_0^v h$, for $h\in\mathbb E_{1,0}$. In particular, we have
$\Diag(\bisimFin)\cong \simFin^\sharp$. Moreover, given a simplicial
set $\mathcal A$ and a map $\beta:\mathcal K\to \simFin$ of
simplicial sets, there is an isomorphism over $\simFin^\sharp$
\begin{equation}
\label{eq:marked_diagonal_external_product}
 \Diag((\mathcal K\boxtimes \mathcal A)_\beta)
 \cong \mathcal K^\sharp_\beta\times \mathcal A^\flat,
\end{equation}
natural in $\beta$ and $\mathcal A$.
\end{proposition}

\begin{proof}
Since the underlying simplicial set of a colimit of marked simplicial
sets is the colimit of the underlying simplicial sets, the set of
$q$-simplices of $\Diag(\mathbb E)$ is
\[
 \int^{[m]\in\Delta}\mathbb E_{m,q}
 \times\operatorname{Hom}_\Delta([q],[m])
 \cong\mathbb E_{q,q},
\]
by the co-Yoneda lemma
\cite[\S3.10, Equation~(3.71), p.~53]{Kelly2005}, where the isomorphism
sends the class of a pair $(x,\alpha)$ to $\alpha_h^*x$, and its
inverse sends $x\in\mathbb E_{q,q}$ to the class of
$(x,\id_{[q]})$. These isomorphisms are compatible with the
simplicial structures, and the structure map to $\simFin$ sends
$x\in\mathbb E_{q,q}$ to $P_{q,q}(x)\in (\simFin)_q$. Moreover, the
marked edges of a colimit of marked simplicial sets are the images of
the marked edges, and an edge of
$\mathbb E_m^\flat\times\Delta[m]^\sharp$ is marked if and only if
its first component is degenerate. Hence the marked edges of
$\Diag(\mathbb E)$ are the images $s_0^v\alpha_h^*x$ of the pairs
$(s_0^v x,\alpha)$, with $x\in\mathbb E_{m,0}$ and
$\alpha:[1]\to[m]$, i.e., the edges $s_0^v h$ with
$h\in\mathbb E_{1,0}$. In particular, every edge of
$\Diag(\bisimFin)$ is marked, as the simplicial sets
$(\bisimFin)_m$ are discrete. Finally, the diagonal of
$\mathcal K\boxtimes \mathcal A$ is $\mathcal K\times \mathcal A$,
and, by the above description, an edge of
$\Diag((\mathcal K\boxtimes \mathcal A)_\beta)$ is marked if and only
if its $\mathcal A$-component is degenerate, i.e., if and only if it
is marked in $\mathcal K^\sharp\times \mathcal A^\flat$. As both
structure maps to $\simFin$ are given by
$\beta\operatorname{pr}_{\mathcal K}$, we get the isomorphism
\eqref{eq:marked_diagonal_external_product}, which is the identity on
underlying simplicial sets, and hence natural in $\beta$ and
$\mathcal A$.
\end{proof}

The marked diagonal is left adjoint to the functor $\Int$.

\begin{proposition}
\label{prop:marked_diagonal_adjunction}
There is an adjunction
\[
 \Diag:(\sset^{\Delta^{\mathrm{op}}})_{/\bisimFin}
 \rightleftarrows(\sset^+)_{/\simFin^\sharp}:\Int.
\]
\end{proposition}

\begin{proof}
By the universal property of the coend, a map of marked simplicial
sets $\Diag(\mathbb E)\to \mathcal X$ over $\simFin^\sharp$ amounts
to a family of maps of marked simplicial sets over $\simFin^\sharp$
\[
 \Phi_{m,\sigma}:\mathbb E_m(\sigma)^\flat
 \times\Delta[m]^\sharp_\sigma\longrightarrow \mathcal X,
\]
for $m\geq0$ and $\sigma\in(\simFin)_m$, such that
\[
 \Phi_{\ell,\sigma\theta}\circ((\theta_{\mathbb E}^*)^\flat\times\id)
 =\Phi_{m,\sigma}\circ(\id\times\theta^\sharp)
\]
for every map $\theta:[\ell]\to[m]$ in $\Delta$. By the definition
of $\Map^\flat$, the maps $\Phi_{m,\sigma}$ correspond to maps of
simplicial sets
\[
 g_{m,\sigma}:\mathbb E_m(\sigma)\longrightarrow
 \Map^\flat(\Delta[m]^\sharp_\sigma,\mathcal X)
 =\Int(\mathcal X)_m(\sigma),
\]
where $g_{m,\sigma}$ sends an $r$-simplex
$a:\Delta[r]\to\mathbb E_m(\sigma)$ to the composite
$\Phi_{m,\sigma}(a^\flat\times\id)$, and the above compatibility
condition corresponds to the condition
\[
 g_{\ell,\sigma\theta}\theta_{\mathbb E}^*
 =\theta_{\Int(\mathcal X)}^*g_{m,\sigma},
\]
i.e., to the condition that the maps $g_{m,\sigma}$ assemble into a
map of bisimplicial sets $\mathbb E\to\Int(\mathcal X)$ over
$\bisimFin$. The resulting bijection is natural in $\mathbb E$ and
$\mathcal X$.
\end{proof}

Next, we show that the adjunction $\Diag\dashv\Int$ is compatible
with the model structures.

\begin{proposition}
\label{prop:marked_diagonal_quillen_pair}
The adjunction $\Diag\dashv\Int$ of
\cref{prop:marked_diagonal_adjunction} induces Quillen adjunctions
\[
\begin{gathered}
 \Diag:\reedyslice\rightleftarrows\cocart:\Int,\\[3pt]
 \Diag:\ldblfib\rightleftarrows\cocart:\Int.
\end{gathered}
\]
\end{proposition}

\begin{proof}
We first show that the functor $\Int$ preserves fibrations and
trivial fibrations from $\cocart$ to $\reedyslice$. Let
$g:\mathcal X\to \mathcal Y$ be a fibration in $\cocart$. For
$\sigma\in (\simFin)_m$, write
$\partial\sigma=\sigma|_{\partial\Delta[m]}$ and
$i:(\partial\Delta[m])^\sharp_{\partial\sigma}
\hookrightarrow\Delta[m]^\sharp_\sigma$ for the inclusion. By
\cref{prop:matching_identification}, the relative matching map of
$\Int(g)$ in degree $m$ is the coproduct over
$\sigma\in(\simFin)_m$ of the maps
\[
 R_{i,g}:\Map^\flat(\Delta[m]^\sharp_\sigma,\mathcal X)\longrightarrow
 \Map^\flat((\partial\Delta[m])^\sharp_{\partial\sigma},\mathcal X)
 \mathop{\times}_{\Map^\flat((\partial\Delta[m])^\sharp_{\partial\sigma},\mathcal Y)}
 \Map^\flat(\Delta[m]^\sharp_\sigma,\mathcal Y)
\]
induced by $i$ and $g$. Given a cofibration
$j:\mathcal K\hookrightarrow\mathcal L$ in $\sset_{(\infty,1)}$, a
lifting problem of $j$ against $R_{i,g}$ corresponds, by
adjunction, to a lifting problem of the pushout-product
$j^\flat\mathbin{\square}i$ against $g$. By
\cref{lem:cocartesian_joyal_pushout_product}, the map
$j^\flat\mathbin{\square}i$ is a cofibration in $\cocart$, which is
trivial if $j$ is a weak equivalence in $\sset_{(\infty,1)}$. Hence
the map $R_{i,g}$ is a fibration in $\sset_{(\infty,1)}$, and it is
a trivial fibration if $g$ is one. By
\cref{lem:fixed_index_coproducts}, the same holds for the relative
matching maps of $\Int(g)$, which shows that the first adjunction is
a Quillen adjunction. For the second adjunction, recall from
\cref{thm:left_double_model_structure} that
$\ldblfib=L_{\mathsf J}\reedyslice$. Since every object of
$\reedyslice$ is cofibrant, it is enough to show, by
\cref{prop:localization_recollections}(2), that the functor $\Diag$
sends the maps $j_{m,\sigma,\mathcal A}$ in $\mathsf J$ to weak
equivalences in $\cocart$. By
\eqref{eq:marked_diagonal_external_product}, the map
$\Diag(j_{m,\sigma,\mathcal A})$ is the map
\[
 \Delta[0]^\sharp_{\sigma(0)}\times \mathcal A^\flat
 \longrightarrow\Delta[m]^\sharp_\sigma\times \mathcal A^\flat.
\]
Given a fibrant object $\mathcal X$ of $\cocart$, the restriction
map along this map between the mapping objects $\Map^\flat$ into
$\mathcal X$ is identified, by adjunction, with the map
\[
 (\lambda_m^\sigma)^{\mathcal A}:
 \bigl(\Int(\mathcal X)_m(\sigma)\bigr)^{\mathcal A}
 \longrightarrow\bigl(\Int(\mathcal X)_0(\sigma(0))\bigr)^{\mathcal A}.
\]
Since $\Int(\mathcal X)\to\bisimFin$ is a left double
$(\infty,1)$-fibration by \cref{cor:left-double}, the map
$\lambda_m^\sigma$ is a trivial fibration by
\cref{cor:fixed_simplex_source}, and hence so is the map
$(\lambda_m^\sigma)^{\mathcal A}$, as $\sset_{(\infty,1)}$ is
cartesian. Hence $\Diag(j_{m,\sigma,\mathcal A})$ is a weak
equivalence in $\cocart$ by
\cref{prop:cocartesian_weak_equivalences}.
\end{proof}

\subsection{The Quillen equivalences}
\label{subsec:comparison_quillen_equivalences}

In this section, we show that the adjunction $\Diag\dashv\Int$ is a
Quillen equivalence between the model structures $\ldblfib$ and
$\cocart$, and between the model structures $\ldblfibsp$ and
$\cocartsp$.

\subsubsection{The first round-trip comparison}
\label{subsubsec:first_roundtrip}

We first construct, for every left double $(\infty,1)$-fibration
$P:\mathbb E\to\bisimFin$, a natural weak equivalence
$\chi_{\mathbb E}:\mathbb E\to\Int(\Stair(\mathbb E)^\natural)$ in
$\ldblfib$. To show that the assignment
$\mathbb E\mapsto\Stair(\mathbb E)^\natural$ is functorial, we need
the following stability property of coCartesian edges. Recall that the
full sub-quasi-category of a quasi-category $\mathcal C$ spanned by
a set of vertices consists of the simplices of $\mathcal C$ all of
whose vertices lie in this set, and that it is \emph{replete} if it
contains every vertex of $\mathcal C$ that is equivalent to one of
its vertices.

\begin{lemma}
\label{lem:cocartesian_edges_replete}
Let $p:\mathcal Y\to \simFin$ be a coCartesian fibration, and let
$\mathcal A^+$ be a marked simplicial set over $\simFin^\sharp$
with underlying simplicial set $\mathcal A$. Then the simplicial set
$\Map(\mathcal A,\mathcal Y)$ is a quasi-category, and the inclusion
\[
 \Map^\flat(\mathcal A^+,\mathcal Y^\natural)
 \hookrightarrow \Map(\mathcal A,\mathcal Y)
\]
identifies $\Map^\flat(\mathcal A^+,\mathcal Y^\natural)$ with the
replete full sub-quasi-category spanned by the maps
$\mathcal A\to\mathcal Y$ over $\simFin$ that send the marked edges
of $\mathcal A^+$ to $p$-coCartesian edges. In particular, for an
edge $f$ of $\simFin$, the quasi-category $\Cocart_p(f)$ is the
replete full sub-quasi-category of $\Map(\Delta[1]_f,\mathcal Y)$
spanned by the $p$-coCartesian edges over $f$.
\end{lemma}

\begin{proof}
By the dual of \cite[Proposition~3.3.1.7]{lurie2009higher}, the map
$p$ is a fibration in $\sset_{(\infty,1)}$, and hence so is the map
$p^{\mathcal A}:\mathcal Y^{\mathcal A}\to \simFin^{\mathcal A}$, as
$\sset_{(\infty,1)}$ is cartesian. Hence its fiber
$\Map(\mathcal A,\mathcal Y)$ is a quasi-category. An $r$-simplex
of $\Map(\mathcal A,\mathcal Y)$ belongs to
$\Map^\flat(\mathcal A^+,\mathcal Y^\natural)$ if and only if the
corresponding map $\Delta[r]\times\mathcal A\to\mathcal Y$ sends the
marked edges of $\Delta[r]^\flat\times \mathcal A^+$ to
$p$-coCartesian edges. As these marked edges are the edges with
degenerate $\Delta[r]$-component and marked $\mathcal A^+$-component,
this condition holds for an $r$-simplex if and only if it holds for
each of its vertices, which shows that
$\Map^\flat(\mathcal A^+,\mathcal Y^\natural)$ is a full
sub-quasi-category. It remains to show that it is replete. We first
consider the case $\mathcal A^+=\Delta[1]^\sharp_f$. An equivalence
in $\Map(\Delta[1]_f,\mathcal Y)$ between edges $g:y_0\to y_1$ and
$g':y'_0\to y'_1$ over $f$ is represented by a map
$\Delta[1]\times\Delta[1]\to \mathcal Y$ over $f\pr_2$, depicted as
\[
\begin{tikzcd}[ampersand replacement=\&, sep=large]
y_0\arrow[r,"g"]\arrow[d,"b_0"',"\simeq"]\arrow[rd,"d"]\&
y_1\arrow[d,"b_1","\simeq"']\\
y'_0\arrow[r,"g'"']\&y'_1,
\end{tikzcd}
\]
where the edges $b_0$ and $b_1$ are equivalences in $\mathcal Y$
lying over degenerate edges of $\simFin$, and hence
$p$-coCartesian by the dual of
\cite[Proposition~2.4.1.5]{lurie2009higher}. If $g$ is
$p$-coCartesian, then so is the diagonal $d$, as a composite of
$g$ and $b_1$, and hence so is $g'$, as $d$ is a composite of
$b_0$ and $g'$, by the dual of
\cite[Proposition~2.4.1.7]{lurie2009higher}. The converse follows by
applying the same argument to an inverse equivalence. In the general
case, if $F$ and $F'$ are equivalent objects of
$\Map(\mathcal A,\mathcal Y)$, then their restrictions along a marked
edge of $\mathcal A^+$ over an edge $f$ of $\simFin$ are
equivalent in $\Map(\Delta[1]_f,\mathcal Y)$, and so $F$ sends this
marked edge to a $p$-coCartesian edge if and only if $F'$ does.
\end{proof}

We now extend the maps $\iota_a$ of
\cref{lem:staircase_fiber_comparison} and $\chi_f$ of
\cref{prop:staircase_transport_comparison} to a map of bisimplicial
sets.

\begin{construction}
\label{const:chi_roundtrip}
Let $P:\mathbb E\to\bisimFin$ be a left double
$(\infty,1)$-fibration. Given an element $x\in\mathbb E_{m,r}$ over
$\sigma\in(\simFin)_m$, represented by a map of bisimplicial sets
$x:\Delta[m]\boxtimes\Delta[r]\to\mathbb E$, we define $\chi_m(x)$
to be the composite
\begin{equation}
\label{eq:compat_chi_formula}
 \Delta[r]\times\Delta[m]
 \xrightarrow{\mathrm{tw}}\Delta[m]\times\Delta[r]
 \cong\Stair(\Delta[m]\boxtimes\Delta[r])
 \xrightarrow{\Stair(x)}\Stair(\mathbb E),
\end{equation}
where the isomorphism is that of
\eqref{eq:staircase_external_product}. As in the proof of
\cref{prop:staircase_transport_comparison}, we have
$\Stair(P)\chi_m(x)=\sigma\operatorname{pr}_2$, and the map
$\chi_m(x)$ sends the edge $\{a\}\times\alpha$, for a vertex $a$
of $\Delta[r]$ and an edge $\alpha$ of $\Delta[m]$, to the pure
horizontal staircase edge $(h,s_0^v d_0^h h)$, where
$h\in\mathbb E_{1,0}$ is the restriction of $x$ along
$\alpha\boxtimes a$, which is $\Stair(P)$-coCartesian by
\cref{prop:cocartesian_staircase_edges}. Hence $\chi_m(x)$ is an
$r$-simplex of
\[
 \Int(\Stair(\mathbb E)^\natural)_m(\sigma)
 =\Map^\flat
   (\Delta[m]^\sharp_\sigma,\Stair(\mathbb E)^\natural),
\]
and we denote by
\[
\chi_{\mathbb E}:\mathbb E\longrightarrow
\Int(\Stair(\mathbb E)^\natural)
\]
the resulting family of maps over $\bisimFin$. Note that the
restriction of $\chi_1$ to the fibers over an edge $f$ of
$\simFin$ is the map $\chi_f$ of
\cref{prop:staircase_transport_comparison}.
\end{construction}

We now show that $\chi_{\mathbb E}$ is a natural map of bisimplicial
sets.

\begin{proposition}
\label{prop:chi_is_simplicial}
The assignment $\mathbb E\mapsto\Stair(\mathbb E)^\natural$ defines a
functor from the full subcategory of left double
$(\infty,1)$-fibrations over $\bisimFin$ to the category of marked
simplicial sets over $\simFin^\sharp$, and the maps
$\chi_{\mathbb E}:\mathbb E\to\Int(\Stair(\mathbb E)^\natural)$ are
maps of bisimplicial sets over $\bisimFin$ that are natural in
$\mathbb E$.
\end{proposition}

\begin{proof}
Given maps $\alpha:[\ell]\to[m]$ and $\beta:[r']\to[r]$ in
$\Delta$, the elements $\alpha_h^*x$ and $\beta_v^*x$ are
represented by the maps $x(\alpha\boxtimes\id)$ and
$x(\id\boxtimes\beta)$, respectively. Hence, by the functoriality of
$\Stair$ and the naturality of the isomorphism
\eqref{eq:staircase_external_product}, we have
\[
 \chi_\ell(\alpha_h^*x)=(\id\times\alpha)^*\chi_m(x),
 \qquad
 \chi_m(\beta_v^*x)=(\beta\times\id)^*\chi_m(x),
\]
which shows that $\chi_{\mathbb E}$ is a map of bisimplicial sets
over $\bisimFin$. Next, let $F:\mathbb E\to\mathbb E'$ be a map
over $\bisimFin$ between left double $(\infty,1)$-fibrations $P$
and $P'$. We show that the map $\Stair(F)$ sends
$\Stair(P)$-coCartesian edges to $\Stair(P')$-coCartesian edges.
Since $F$ is compatible with the face and degeneracy maps, the map
$\Stair(F)$ sends a pure horizontal staircase edge
$(h,s_0^v d_0^h h)$ to the pure horizontal staircase edge
$(Fh,s_0^v d_0^h Fh)$. Let $e$ be a $\Stair(P)$-coCartesian edge
over an edge $f$ of $\simFin$. Since the map $\chi_f$ is a weak
equivalence between quasi-categories by
\cref{prop:staircase_transport_comparison}, there is a vertex $h$ of
$\mathbb E_1(f)$ such that $e$ is equivalent to the pure horizontal
staircase edge $\chi_f(h)$ in
$\Map(\Delta[1]_f,\Stair(\mathbb E))$. As postcomposition with
$\Stair(F)$ preserves equivalences, the edge $\Stair(F)(e)$ is
equivalent to a pure horizontal staircase edge, and hence it is
$\Stair(P')$-coCartesian by
\cref{prop:cocartesian_staircase_edges,lem:cocartesian_edges_replete}.
Hence $\Stair(F)$ induces a map
$\Stair(F)^\natural:\Stair(\mathbb E)^\natural\to
\Stair(\mathbb E')^\natural$ of marked simplicial sets over
$\simFin^\sharp$, and this assignment is functorial, as $\Stair$
is. Finally, applying $\Stair$ to the map $Fx$ in
\eqref{eq:compat_chi_formula} gives
\[
 \Int(\Stair(F)^\natural)\chi_{\mathbb E}
 =\chi_{\mathbb E'}F,
\]
i.e., the maps $\chi_{\mathbb E}$ are natural in $\mathbb E$.
\end{proof}

We can now prove the first comparison result:

\begin{theorem}
\label{thm:first_roundtrip_comparison}
For every left double $(\infty,1)$-fibration
$P:\mathbb E\to\bisimFin$, the map
\[
\chi_{\mathbb E}:\mathbb E\longrightarrow
\Int(\Stair(\mathbb E)^\natural)
\]
is a levelwise weak equivalence in $\sset_{(\infty,1)}$. In
particular, it is a weak equivalence in $\ldblfib$.
\end{theorem}

\begin{proof}
By \eqref{eq:compat_chi_formula} and
\eqref{eq:staircase_fiber_inclusion}, under the identification
$\Int(\Stair(\mathbb E)^\natural)_0(a)=\Stair(\mathbb E)_a$, the
restriction of $\chi_0$ to the fibers over an object $a$ of
$\Fin$ is the map $\iota_a$, which is a weak equivalence in
$\sset_{(\infty,1)}$ by \cref{lem:staircase_fiber_comparison}. Hence
the map $\chi_0=\coprod_{a\in (\simFin)_0}\iota_a$ is a weak
equivalence in $\sset_{(\infty,1)}$ by
\cref{lem:fixed_index_coproducts}. Since $P$ and
$\Int(\Stair(P))$ are left double $(\infty,1)$-fibrations, by
\cref{cor:staircase_is_cocartesian_fibration,cor:left-double}, the
result follows from \cref{thm:left_double_model_structure}.
\end{proof}

As a consequence, we obtain the following:

\begin{corollary}
\label{cor:first_roundtrip_special}
A left double $(\infty,1)$-fibration $P:\mathbb E\to\bisimFin$ is
special if and only if the left double $(\infty,1)$-fibration
$\Int(\Stair(P))$ is special.
\end{corollary}

\begin{proof}
This follows directly from
\cref{thm:staircase_preserves_reflects_specialness,cor:strings_specialness}.
\end{proof}

\subsubsection{From the marked diagonal to the staircase}
\label{subsubsec:diagonal_to_staircase}

We now compare the marked diagonal $\Diag(\mathbb E)$ of a left
double $(\infty,1)$-fibration $P:\mathbb E\to\bisimFin$ with its
staircase $\Stair(\mathbb E)^\natural$, using the comparison map
from the diagonal to the total simplicial set of a bisimplicial set
\cite{CegarraRemedios2005,Stevenson2025TotalSimplicialSet}.

\begin{construction}[{\cite[Definition~3.6]{Stevenson2025TotalSimplicialSet}, with the variables interchanged}]
\label{const:diagonal_to_staircase}
Given a bisimplicial set $\mathbb X$, we denote by
\[
 \psi_{\mathbb X}:\diagset\mathbb X\longrightarrow\Stair(\mathbb X),
 \qquad
 (\psi_{\mathbb X}(z))_j
 =z|_{[j,q]\times[0,j]}\in\mathbb X_{q-j,j}
 \quad(z\in\mathbb X_{q,q}),
\]
the natural map of simplicial sets which, under the identification
\eqref{eq:staircase_total_identification}, is the comparison map of
\cite[Definition~3.6]{Stevenson2025TotalSimplicialSet}. Here the
intervals $[j,q]$ and $[0,j]$ are identified with the
corresponding faces of $\Delta[q]$ in the horizontal and vertical
directions, respectively, and the relations
\eqref{eq:staircase_relation} hold as both
$d_0^h(\psi_{\mathbb X}(z))_{j-1}$ and
$d_j^v(\psi_{\mathbb X}(z))_j$ are the restriction of $z$ to
$[j,q]\times[0,j-1]$. Under the isomorphisms
$\diagset\bisimFin\cong \simFin$ and
$\Stair(\bisimFin)\cong \simFin$, the map $\psi_{\bisimFin}$ is the
identity, and so, for a bisimplicial set $P:\mathbb E\to\bisimFin$
over $\bisimFin$, the map $\psi_{\mathbb E}$ is a map over
$\simFin$. On edges, the map $\psi_{\mathbb E}$ is given by
\[
 \psi_{\mathbb E}(z)=(d_1^v z,d_0^h z),
 \qquad\text{so that}\qquad
 \psi_{\mathbb E}(s_0^v h)=(h,s_0^v d_0^h h).
\]
Hence, if $P$ is a left double $(\infty,1)$-fibration, then, by
\cref{prop:marked_diagonal_description}, the map $\psi_{\mathbb E}$
sends the marked edges of $\Diag(\mathbb E)$ to pure horizontal
staircase edges, which are $\Stair(P)$-coCartesian, and so it
defines a map of marked simplicial sets over $\simFin^\sharp$
\[
 \psi_{\mathbb E}:\Diag(\mathbb E)
 \longrightarrow\Stair(\mathbb E)^\natural.
\]
\end{construction}

To show that $\psi_{\mathbb E}$ is a weak equivalence in
$\cocart$, we use the following criterion.

\begin{lemma}
\label{lem:enough_marked_cocartesian_lifts}
Let $p:\mathcal Y\to \simFin$ be a coCartesian fibration and let
$a:\mathcal A^+\to \mathcal Y^\natural$ be a map of marked
simplicial sets over $\simFin^\sharp$. Suppose that the underlying
map of $a$ is a weak equivalence in $\sset_{(\infty,1)}$ and that,
for every edge $f$ of $\simFin$, every $p$-coCartesian edge over
$f$ is equivalent in $\Map(\Delta[1]_f,\mathcal Y)$ to the image
under $a$ of a marked edge of $\mathcal A^+$ over $f$. Then $a$
is a weak equivalence in $\cocart$.
\end{lemma}

\begin{proof}
Let $p':\mathcal Z\to \simFin$ be a coCartesian fibration. By
\cref{prop:cocartesian_weak_equivalences}, it is enough to show that
the map
\begin{equation}
\label{eq:compat_marked_mapping_equivalence}
 \Map^\flat(\mathcal Y^\natural,\mathcal Z^\natural)
 \longrightarrow
 \Map^\flat(\mathcal A^+,\mathcal Z^\natural)
\end{equation}
induced by $a$ is a weak equivalence in $\sset_{(\infty,1)}$. By
the dual of \cite[Proposition~3.3.1.7]{lurie2009higher}, the map
$p'$ is a fibration in $\sset_{(\infty,1)}$, i.e., $\mathcal Z$
is a fibrant object of the slice $(\sset_{(\infty,1)})_{/\simFin}$,
which is enriched over $\sset_{(\infty,1)}$ with the mapping objects
of \eqref{eq:unmarked_relative_mapping}. Since every object of this
slice is cofibrant, the weak equivalence $a$ induces a weak
equivalence between quasi-categories
\begin{equation}
\label{eq:compat_ambient_mapping_equivalence}
 a^*:\Map(\mathcal Y,\mathcal Z)\longrightarrow\Map(\mathcal A,\mathcal Z);
\end{equation}
see \cite[Corollary~7.7.2\textup{(2)}]{hirschhorn2003model}. By
\cref{lem:cocartesian_edges_replete}, the source and target of
\eqref{eq:compat_marked_mapping_equivalence} are the replete full
sub-quasi-categories of the source and target of
\eqref{eq:compat_ambient_mapping_equivalence} spanned by the maps that
send marked edges to $p'$-coCartesian edges. We claim that a vertex
$F$ of $\Map(\mathcal Y,\mathcal Z)$ sends $p$-coCartesian edges
to $p'$-coCartesian edges if and only if the map $Fa$ sends the
marked edges of $\mathcal A^+$ to $p'$-coCartesian edges. The
forward implication holds as $a$ is a map of marked simplicial sets.
Conversely, let $e$ be a $p$-coCartesian edge over an edge $f$
of $\simFin$. By assumption, there is a marked edge $g$ of
$\mathcal A^+$ over $f$ such that $e$ is equivalent to $a(g)$
in $\Map(\Delta[1]_f,\mathcal Y)$. Then $F(e)$ is equivalent to
$Fa(g)$ in $\Map(\Delta[1]_f,\mathcal Z)$, and so $F(e)$ is
$p'$-coCartesian by \cref{lem:cocartesian_edges_replete}, which
proves the claim. Hence the map
\eqref{eq:compat_marked_mapping_equivalence} is the restriction of the
map $a^*$ to the full sub-quasi-category
$\Map^\flat(\mathcal Y^\natural,\mathcal Z^\natural)$, which is the
preimage under $a^*$ of
$\Map^\flat(\mathcal A^+,\mathcal Z^\natural)$. In particular, the
map \eqref{eq:compat_marked_mapping_equivalence} is fully faithful, as
$a^*$ is. Moreover, given a vertex $H$ of
$\Map^\flat(\mathcal A^+,\mathcal Z^\natural)$, there is a vertex
$F$ of $\Map(\mathcal Y,\mathcal Z)$ such that $Fa$ is equivalent
to $H$ in $\Map(\mathcal A,\mathcal Z)$, as $a^*$ is essentially
surjective. Since
$\Map^\flat(\mathcal A^+,\mathcal Z^\natural)$ is a replete full
sub-quasi-category, the vertex $Fa$ and the equivalence between
$Fa$ and $H$ lie in it, and so $F$ is a vertex of
$\Map^\flat(\mathcal Y^\natural,\mathcal Z^\natural)$ by the claim.
Hence the map \eqref{eq:compat_marked_mapping_equivalence} is also
essentially surjective, and so it is a weak equivalence in
$\sset_{(\infty,1)}$; see \cite[\S1.2.10]{lurie2009higher}.
\end{proof}

With this criterion, we can show that the staircase is a fibrant
replacement of the marked diagonal:

\begin{proposition}
\label{prop:diagonal_to_staircase_cocart_equivalence}
Let $P:\mathbb E\to\bisimFin$ be a left double
$(\infty,1)$-fibration. Then the map
\[
\psi_{\mathbb E}:\Diag(\mathbb E)\longrightarrow
\Stair(\mathbb E)^\natural
\]
is a weak equivalence in $\cocart$ with fibrant target, which is
natural in $P$. Moreover, the map $\psi_{\mathbb E}$ corresponds to
the map $\chi_{\mathbb E}$ of \cref{const:chi_roundtrip} under the
adjunction $\Diag\dashv\Int$, i.e., we have
$\chi_{\mathbb E}=\Int(\psi_{\mathbb E})\eta_{\mathbb E}$, where
$\eta$ denotes the unit of this adjunction.
\end{proposition}

\begin{proof}
The object $\Stair(\mathbb E)^\natural$ is fibrant in $\cocart$ by
\cref{cor:staircase_is_cocartesian_fibration}, and the map
$\psi_{\mathbb E}$ is natural in $P$ by
\cref{const:diagonal_to_staircase,prop:chi_is_simplicial}. By
\cite[Theorem~1.2]{Stevenson2025TotalSimplicialSet}, the comparison
map from the diagonal to the total simplicial set of a bisimplicial
set is a weak equivalence in $\sset_{(\infty,1)}$. Applying this
result to the bisimplicial set $\mathbb E^{\mathrm{tr}}$, we get
that the underlying map $\diagset\mathbb E\to\Stair(\mathbb E)$ of
$\psi_{\mathbb E}$ is a weak equivalence in $\sset_{(\infty,1)}$.
Moreover, let $e$ be a $\Stair(P)$-coCartesian edge over an edge
$f$ of $\simFin$. As in the proof of \cref{prop:chi_is_simplicial},
there is a vertex $h$ of $\mathbb E_1(f)$ such that $e$ is
equivalent in $\Map(\Delta[1]_f,\Stair(\mathbb E))$ to the edge
$\chi_f(h)=(h,s_0^v d_0^h h)=\psi_{\mathbb E}(s_0^v h)$, which is the
image of the marked edge $s_0^v h$ of $\Diag(\mathbb E)$ over
$f$. Hence $\psi_{\mathbb E}$ is a weak equivalence in $\cocart$
by \cref{lem:enough_marked_cocartesian_lifts}. It remains to show that
$\chi_{\mathbb E}$ corresponds to $\psi_{\mathbb E}$ under the
adjunction. By the proofs of
\cref{prop:marked_diagonal_description,prop:marked_diagonal_adjunction},
the adjunct of $\psi_{\mathbb E}$ sends an element
$x\in\mathbb E_{m,r}$, represented by a map
$x:\Delta[m]\boxtimes\Delta[r]\to\mathbb E$, to the composite
\[
 \Delta[r]\times\Delta[m]
 \xrightarrow{\mathrm{tw}}
 \diagset(\Delta[m]\boxtimes\Delta[r])
 \xrightarrow{\diagset(x)}\diagset\mathbb E
 \xrightarrow{\psi_{\mathbb E}}\Stair(\mathbb E),
\]
and, by naturality of $\psi$, we have
$\psi_{\mathbb E}\diagset(x)=
\Stair(x)\psi_{\Delta[m]\boxtimes\Delta[r]}$. Moreover, for an
external product $\mathcal K\boxtimes \mathcal A$, the map
$\psi_{\mathcal K\boxtimes\mathcal A}$ is given by
\[
 (x,y)\longmapsto
 \bigl(x|_{[j,q]},y|_{[0,j]}\bigr)_{j=0}^q
 \qquad(x\in\mathcal K_q,\ y\in\mathcal A_q),
\]
and so it is the inverse of the isomorphism
\eqref{eq:staircase_external_product}. Hence the above composite is
the map $\chi_m(x)$ of \eqref{eq:compat_chi_formula}, as desired.
\end{proof}

\subsubsection{The Quillen equivalence and the second round trip}
\label{subsubsec:quillen_equivalence}

We now prove that the adjunction $\Diag\dashv\Int$ is a Quillen
equivalence between $\ldblfib$ and $\cocart$. For this, we first
show that the functor $\Int$ reflects weak equivalences between
fibrant objects.

\begin{lemma}
\label{lem:internalization_reflects_fibrant_equivalences}
The right Quillen functor $\Int:\cocart\to\ldblfib$ reflects weak
equivalences between fibrant objects.
\end{lemma}

\begin{proof}
Let $g:\mathcal E^\natural\to \mathcal F^\natural$ be a map between
fibrant objects of $\cocart$ such that $\Int(g)$ is a weak
equivalence in $\ldblfib$. Since $\Int(\mathcal E^\natural)$ and
$\Int(\mathcal F^\natural)$ are fibrant in $\ldblfib$, the map
\[
 \Int(g)_0=\coprod_{a\in (\simFin)_0}g_a:
 \coprod_{a\in (\simFin)_0}\mathcal E_a\longrightarrow\coprod_{a\in (\simFin)_0}\mathcal F_a
\]
is a weak equivalence in $\sset_{(\infty,1)}$ by
\cref{thm:left_double_model_structure}, and so the map
$g_a:\mathcal E_a\to \mathcal F_a$ is a weak equivalence in
$\sset_{(\infty,1)}$ for every object $a$ of $\Fin$, by
\cref{lem:fixed_index_coproducts}. Hence $g$ is a weak equivalence
in $\cocart$ by \cref{prop:cocartesian_weak_equivalences}.
\end{proof}

We can now prove the first main result of this section:

\begin{theorem}
\label{thm:left_double_cocart_quillen_equivalence}
The adjunction $\Diag\dashv\Int$ is a Quillen equivalence
\[\adjunction{\ldblfib}{\cocart}{\Diag}{\Int}
\]
\end{theorem}

\begin{proof}
By \cref{lem:internalization_reflects_fibrant_equivalences} and
\cite[Corollary~1.3.16\textup{(c)}]{hovey1999model}, it is enough to
show that the derived unit of the Quillen adjunction
$\Diag\dashv\Int$ of \cref{prop:marked_diagonal_quillen_pair} is a
weak equivalence in $\ldblfib$ at every object. Let
$i_{\mathcal Z}:\mathcal Z\hookrightarrow R\mathcal Z$ denote a
functorial fibrant replacement in $\cocart$, with $i_{\mathcal Z}$
a trivial cofibration. Since every object of $\ldblfib$ is cofibrant,
the derived unit at an object $\mathbb A$ is represented by the
composite
\[
 \eta^R_{\mathbb A}:
 \mathbb A\xrightarrow{\eta_{\mathbb A}}
 \Int\Diag(\mathbb A)
 \xrightarrow{\Int(i_{\Diag(\mathbb A)})}
 \Int R\Diag(\mathbb A),
\]
where $\eta$ denotes the unit of the adjunction. Suppose first that
$\mathbb A=\mathbb E$ is fibrant in $\ldblfib$, with structure map
$P:\mathbb E\to\bisimFin$. By
\cref{prop:diagonal_to_staircase_cocart_equivalence}, the map
$\psi_{\mathbb E}:\Diag(\mathbb E)\to\Stair(\mathbb E)^\natural$ is
a weak equivalence in $\cocart$ with fibrant target, and so there is
a lift $\overline\psi_{\mathbb E}$ in the following commutative
square, where $b_{\mathbb E}$ denotes the structure map of
$R\Diag(\mathbb E)$.
\[
\begin{tikzcd}[ampersand replacement=\&,column sep=large]
 \Diag(\mathbb E)
   \arrow[r,"\psi_{\mathbb E}"]
   \arrow[d,hook,"i_{\Diag(\mathbb E)}"']
 \&\Stair(\mathbb E)^\natural
   \arrow[d,two heads,"\Stair(P)"]
 \\
 R\Diag(\mathbb E)
   \arrow[r,"b_{\mathbb E}"']
   \arrow[ur,dashed,"\overline\psi_{\mathbb E}"']
 \&\simFin^\sharp.
\end{tikzcd}
\]
By the two-out-of-three property, the map
$\overline\psi_{\mathbb E}$ is a weak equivalence between fibrant
objects of $\cocart$, and hence the map
$\Int(\overline\psi_{\mathbb E})$ is a weak equivalence in
$\ldblfib$ by \cref{prop:model_category_recollections}(2). Moreover,
by \cref{prop:diagonal_to_staircase_cocart_equivalence}, we have
\[
 \chi_{\mathbb E}
 =\Int(\psi_{\mathbb E})\eta_{\mathbb E}
 =\Int(\overline\psi_{\mathbb E})\eta^R_{\mathbb E}.
\]
Since $\chi_{\mathbb E}$ is a weak equivalence by
\cref{thm:first_roundtrip_comparison}, so is the map
$\eta^R_{\mathbb E}$, by the two-out-of-three property. For an
arbitrary object $\mathbb A$, let
$r:\mathbb A\to\mathbb A^{\mathrm f}$ be a fibrant replacement in
$\ldblfib$. By naturality of $\eta$ and of $i$, we have a
commutative square of bisimplicial sets
\[
\begin{tikzcd}[ampersand replacement=\&, sep=large]
 \mathbb A
   \arrow[r,"\eta^R_{\mathbb A}"]
   \arrow[d,"r"']
 \&\Int R\Diag(\mathbb A)
   \arrow[d,"\Int R\Diag(r)"]
 \\
 \mathbb A^{\mathrm f}
   \arrow[r,"\eta^R_{\mathbb A^{\mathrm f}}"']
 \&\Int R\Diag(\mathbb A^{\mathrm f}).
\end{tikzcd}
\]
The map $\Diag(r)$ is a weak equivalence in $\cocart$, as the
functor $\Diag$ is left Quillen and every object of $\ldblfib$ is
cofibrant, and hence so is the map $R\Diag(r)$, by the
two-out-of-three property. As $R\Diag(r)$ is a map between fibrant
objects, the map $\Int R\Diag(r)$ is a weak equivalence. Since $r$
and $\eta^R_{\mathbb A^{\mathrm f}}$ are also weak equivalences, we
conclude that $\eta^R_{\mathbb A}$ is a weak equivalence, by the
two-out-of-three property.
\end{proof}

\begin{rmk}[Other bases]
\label{rem:general_base}
The category $\Fin$ does not play a special role in
\cref{thm:left_double_cocart_quillen_equivalence}. Given a small
category $\mathrm C$, let $\mathcal B=N(\mathrm C)$ and
$\mathbb B=\mathcal B^h$, as in
\cref{subsec:prelim_simplicial_models}. Then the simplicial object
$\mathbb B$ is Reedy fibrant and its Segal maps are isomorphisms,
the staircase of a left double $(\infty,1)$-fibration over
$\mathbb B$ is a coCartesian fibration over $\mathcal B$, and the
functor $\Int$ sends a coCartesian fibration over $\mathcal B$ to a
left double $(\infty,1)$-fibration over $\mathbb B$. The proof of
\cref{thm:left_double_cocart_quillen_equivalence} then gives a Quillen
equivalence between the model structure for left double
$(\infty,1)$-fibrations over $\mathbb B$, defined as in
\cref{thm:left_double_model_structure}, and the coCartesian model
structure on marked simplicial sets over $\mathcal B^\sharp$;
compare \cite[Corollary~4.35]{Nuiten2023}. Only the spans of
\cref{subsec:internal_specialness,subsec:enriched_specialness} and the
corresponding localizations are specific to $\Fin$.
\end{rmk}

Finally, we compare a coCartesian fibration over $\simFin$ with the
staircase of its bisimplicial set of coCartesian strings. We denote by
$\varepsilon$ the counit of the adjunction $\Diag\dashv\Int$.

\begin{theorem}
\label{thm:second_roundtrip_comparison}
For every fibrant object $\mathcal X$ of $\cocart$, there is a
zigzag of weak equivalences in $\cocart$
\begin{equation}
\label{eq:compat_second_roundtrip}
 \mathcal X\xleftarrow{\ \varepsilon_{\mathcal X}\ }\Diag\Int(\mathcal X)
 \xrightarrow{\ \psi_{\Int(\mathcal X)}\ }\Stair(\Int(\mathcal X))^\natural,
\end{equation}
which is natural in $\mathcal X$.
\end{theorem}

\begin{proof}
Since $\Int(\mathcal X)\to\bisimFin$ is a left double
$(\infty,1)$-fibration by \cref{cor:left-double}, the map
$\psi_{\Int(\mathcal X)}$ is a weak equivalence in $\cocart$ by
\cref{prop:diagonal_to_staircase_cocart_equivalence}. Moreover, since
$\Int(\mathcal X)$ is cofibrant in $\ldblfib$ and $\mathcal X$ is
fibrant in $\cocart$, the counit $\varepsilon_{\mathcal X}$, whose
adjunct is the identity of $\Int(\mathcal X)$, is a weak equivalence
in $\cocart$ by \cref{thm:left_double_cocart_quillen_equivalence}.
The naturality of the zigzag follows from the naturality of
$\varepsilon$ and of $\psi$, using \cref{prop:chi_is_simplicial}.
\end{proof}

\begin{rmk}
\label{rem:second_roundtrip_single_map}
Since the objects $\mathcal X$ and
$\Stair(\Int(\mathcal X))^\natural$ are cofibrant and fibrant in
$\cocart$, the isomorphism
$[\psi_{\Int(\mathcal X)}][\varepsilon_{\mathcal X}]^{-1}$ in the
homotopy category of $\cocart$ is represented by a map
$\mathcal X\to\Stair(\Int(\mathcal X))^\natural$ of marked
simplicial sets over $\simFin^\sharp$, which is a weak equivalence
in $\cocart$ \cite[Theorem~1.2.10\textup{(ii), (iv)}]{hovey1999model}.
By \cref{prop:cocartesian_weak_equivalences}, the underlying map of
such a representative is a weak equivalence in $\sset_{(\infty,1)}$
that induces weak equivalences on fibers.
\end{rmk}

\subsubsection{The Quillen equivalence for symmetric monoidal objects}
\label{subsubsec:special_localizations}

Finally, we show that the Quillen equivalence of
\cref{thm:left_double_cocart_quillen_equivalence} induces a Quillen
equivalence between the model structures for special left double
$(\infty,1)$-fibrations and for special coCartesian fibrations.
Recall from
\cref{thm:internal_special_model_structure,thm:enriched_special_model_structure}
that
\[
 \ldblfibsp=L_{\mathsf S_{\mathrm{sp}}}\ldblfib,
 \qquad
 \cocartsp=L_{\mathsf S^+_{\mathrm{sp}}}\cocart,
\]
where the sets $\mathsf S_{\mathrm{sp}}$ and
$\mathsf S^+_{\mathrm{sp}}$ are defined in
\eqref{eq:internal_special_localizers} and
\eqref{eq:enriched_special_localizers}. The following result is the
precise form of the theorem stated in the introduction.

\begin{theorem}
\label{thm:special_quillen_equivalence}
The adjunction $\Diag\dashv\Int$ is a Quillen equivalence
\[\adjunction{\ldblfibsp}{\cocartsp}{\Diag}{\Int}\]
Moreover, for every special left double $(\infty,1)$-fibration
$P:\mathbb E\to\bisimFin$, the map $\chi_{\mathbb E}$ of
\cref{thm:first_roundtrip_comparison} is a weak equivalence between
fibrant objects of $\ldblfibsp$, and, for every fibrant object
$\mathcal X$ of $\cocartsp$, the zigzag
\eqref{eq:compat_second_roundtrip} of
\cref{thm:second_roundtrip_comparison} is a zigzag of weak
equivalences in $\cocartsp$ between fibrant objects.
\end{theorem}

\begin{proof}
By the naturality of the isomorphism
\eqref{eq:marked_diagonal_external_product}, for every $k\geq0$ and
$\mathcal A\in\HGen$, the functor $\Diag$ sends the map
$s_{k,\mathcal A}$ of \eqref{eq:internal_special_localizers} to a
map that is isomorphic, as a map over $\simFin^\sharp$, to the map
\[
\mathcal V_k^\sharp\times \mathcal A^\flat\longrightarrow
 \mathcal J_k^\sharp\times \mathcal A^\flat
\]
of \eqref{eq:enriched_special_localizers}. By
\cref{thm:left_double_cocart_quillen_equivalence} and
\cite[Theorem~3.3.20\textup{(1)}]{hirschhorn2003model}, the adjunction
$\Diag\dashv\Int$ is a Quillen equivalence between the left
Bousfield localization of $\ldblfib$ at $\mathsf S_{\mathrm{sp}}$
and the left Bousfield localization of $\cocart$ at the image of
$\mathsf S_{\mathrm{sp}}$ under the left derived functor of
$\Diag$. Since every object of $\ldblfib$ is cofibrant, this image
is the set $\Diag(\mathsf S_{\mathrm{sp}})$, and so the latter
localization is $\cocartsp$, which proves the first assertion. For
the second assertion, recall that the weak equivalences of
$\ldblfib$ and $\cocart$ are weak equivalences in $\ldblfibsp$
and $\cocartsp$, respectively. If $P$ is a special left double
$(\infty,1)$-fibration, then the coCartesian fibration $\Stair(P)$
is special by \cref{thm:staircase_preserves_reflects_specialness}, and
so the left double $(\infty,1)$-fibration $\Int(\Stair(P))$ is
special by \cref{cor:strings_specialness}. Similarly, if
$\mathcal X$ is a fibrant object of $\cocartsp$, then
$\Int(\mathcal X)\to\bisimFin$ is a special left double
$(\infty,1)$-fibration by \cref{cor:strings_specialness}, and so
$\Stair(\Int(\mathcal X))^\natural$ is fibrant in $\cocartsp$ by
\cref{thm:staircase_preserves_reflects_specialness}.
\end{proof}

\appendix
\section{Proof of the staircase extension lemma}
\label{app:staircase_extension}

In this appendix, we prove \cref{lem:staircase_extension}. We use the
notation of \cref{def:relative_staircase_horn}, and
$P:\mathbb E\to\bisimFin$ denotes a left double
$(\infty,1)$-fibration.

\begin{proof}[Proof of \cref{lem:staircase_extension}]
Recall that the relative horn $\eta$ amounts to the
$(q-1)$-simplices $\eta^r=(\eta^r_0,\ldots,\eta^r_{q-1})$ of
$\Stair(\mathbb E)$, for $0\leq r\leq q$ with $r\neq a$, which
satisfy the identities \eqref{eq:relative_staircase_horn_overlap}. We
construct the simplices $\xi_j\in\mathbb E_{q-j,j}$ by induction on
$0\leq j<q$. Suppose that, for every $i<j$, we have constructed a
simplex $\xi_i\in\mathbb E_{q-i,i}$ such that
$P_{q-i}(\xi_i)=c_{\sigma_{\geq i}}^{(i)}$ and such that
\eqref{eq:IH2-horizontal}, \eqref{eq:IH3-vertical}, and
\eqref{eq:IH3-staircase_gluing} hold with $i$ in place of $j$. We
obtain the simplex $\xi_j$ as a lift in a commutative square of
simplicial sets
\begin{equation}\label{eq:stage-j-square}
\begin{tikzcd}[sep=large]
  \mathcal L_j \ar[r,"\Phi_j"] \ar[d,hook] &
  \mathbb E_{q-j} \ar[d,two heads,"R_j","\sim"'] \\
  \Delta[j] \ar[r,"H_j"'] \ar[ru,dashed,"\xi_j"] &
  \mathbb E(\mathcal K_j)\times_{\bisimFin(\mathcal K_j)}(\bisimFin)_{q-j}
\end{tikzcd}
\end{equation}
which we now describe.

We first define the simplicial set $\mathcal K_j\subseteq\Delta[q-j]$
and the map $H_j$. Let
\begin{equation}\label{eq:Kj-cases}
  \mathcal K_j\coloneqq
  \bigcup_{\substack{1\le s\le q-j\\ j+s\ne a}}d_s\Delta[{q-j-1}]
  =
  \begin{cases}
    \mathcal K_{q-j,a-j}, & j<a,\\[1mm]
    \Lambda^0[{q-j}], & j\geq a,
  \end{cases}
\end{equation}
with $\mathcal K_{m,\ell}$ as in
\cref{lem:two_missing_horizontal_faces}; note that
$1\leq a-j<q-j$ if $j<a$. The map $R_j$ is the relative
restriction map along the inclusion
$\mathcal K_j\hookrightarrow\Delta[q-j]$, which is a trivial
fibration by
\cref{prop:left_horns_in_T,lem:two_missing_horizontal_faces}. The
component of $H_j$ in $(\bisimFin)_{q-j}$ is the degenerate
$j$-simplex $c_{\sigma_{\geq j}}^{(j)}$, and its component
$\Delta[j]\to\mathbb E(\mathcal K_j)$ is given on the face
$d_s\Delta[q-j-1]$ of $\mathcal K_j$ by the $j$-simplex
\begin{equation}\label{eq:Hj-face-definition}
  H_j|_{d_s\Delta[q-j-1]}\coloneqq\eta^{j+s}_j\in\mathbb E_{q-j-1,j},
  \qquad 1\le s\le q-j,\ j+s\ne a.
\end{equation}
These simplices are compatible on the intersections of the faces.
Indeed, given $1\le s<t\le q-j$ with $j+s\ne a$ and
$j+t\ne a$, the $j$-th component of the identity
$D_{j+s}\eta^{j+t}=D_{j+t-1}\eta^{j+s}$ is the identity
\[
  d_s^h\eta^{j+t}_j=d_{t-1}^h\eta^{j+s}_j,
\]
by \eqref{eq:staircase_faces}, as $j<j+s$ and $j<j+t-1$. Moreover,
the simplex $\eta^{j+s}_j$ lies over the degenerate $j$-simplex at
\[
 (d_0)^j d_{j+s}\sigma=d_s(d_0)^j\sigma
 =d_s\sigma_{\geq j},
\]
which is the restriction of $\sigma_{\geq j}$ to the face
$d_s\Delta[q-j-1]$, and so the two components of $H_j$ define a
map into the pullback.

Next, we define the simplicial set $\mathcal L_j\subseteq\Delta[j]$
and the map $\Phi_j$. Let $\mathcal L_0\coloneqq\varnothing$ and,
for $j>0$, let $\mathcal L_j$ be the union of the faces
$d_r\Delta[j-1]$ with $0\le r<j$ and $r\ne a$, and of the face
$d_j\Delta[j-1]$, so that
\begin{equation}\label{eq:Lj-cases}
  \mathcal L_j=
  \begin{cases}
    \partial\Delta[j], & 0<j\le a,\\[1mm]
    \Lambda^a[j], & j>a.
  \end{cases}
\end{equation}
For $j>0$, the map $\Phi_j$ is given by
\begin{equation}\label{eq:Phij-faces}
  \Phi_j|_{d_r\Delta[j-1]}\coloneqq
  \begin{cases}
    \eta^r_{j-1}, & 0\le r<j,\ r\ne a,\\[1mm]
    d_0^h\xi_{j-1}, & r=j.
  \end{cases}
\end{equation}
We check that these simplices are compatible on the intersections of
the faces. Given $0\le r<s<j$ with $r,s\ne a$, the
$(j-2)$-nd component of the identity $D_r\eta^s=D_{s-1}\eta^r$ is
the identity
\[
  d_r^v\eta^s_{j-1}=d_{s-1}^v\eta^r_{j-1},
\]
by \eqref{eq:staircase_faces}, as $r<s<j$. Given $0\le r<j$ with
$r\ne a$, we need to show that, for $j\geq2$,
\begin{equation}\label{eq:vertical-final-needed}
  d_r^v(d_0^h\xi_{j-1})=d_{j-1}^v\eta^r_{j-1}.
\end{equation}
If $r<j-1$, then, by \eqref{eq:IH3-vertical} for $\xi_{j-1}$ and by
\eqref{eq:staircase_relation} for $\eta^r$, we have
\[
  d_r^vd_0^h\xi_{j-1}
  =d_0^hd_r^v\xi_{j-1}
  =d_0^h\eta^r_{j-2}
  =d_{j-1}^v\eta^r_{j-1}.
\]
If $r=j-1\ne a$, then, by \eqref{eq:IH3-staircase_gluing} for
$\xi_{j-1}$, by \eqref{eq:IH2-horizontal} for $\xi_{j-2}$ and
$s=1$, and by \eqref{eq:staircase_relation} for $\eta^{j-1}$, we
have
\[
  d_{j-1}^vd_0^h\xi_{j-1}
  =d_0^hd_{j-1}^v\xi_{j-1}
  =d_0^hd_0^h\xi_{j-2}
  =d_0^hd_1^h\xi_{j-2}
  =d_0^h\eta^{j-1}_{j-2}
  =d_{j-1}^v\eta^{j-1}_{j-1}.
\]
Hence the map $\Phi_j$ is well defined. Moreover, it lies over the
degenerate simplex at $\sigma_{\geq j}$, as
$(d_0)^{j-1}d_r\sigma=(d_0)^j\sigma=\sigma_{\geq j}$ for
$0\leq r<j$, and $d_0\sigma_{\geq j-1}=\sigma_{\geq j}$.

We now show that the square \eqref{eq:stage-j-square} commutes, i.e.,
that $R_j\Phi_j=H_j|_{\mathcal L_j}$. The components in
$(\bisimFin)_{q-j}$ agree, as both maps lie over the degenerate
simplex at $\sigma_{\geq j}$. For the components in
$\mathbb E(\mathcal K_j)$, let $1\le s\le q-j$ with $j+s\ne a$.
Given $0\le r<j$ with $r\ne a$, we need to show that
\[
  d_s^h\eta^r_{j-1}=d_r^v\eta^{j+s}_j,
\]
which is the $(j-1)$-st component of the identity
$D_r\eta^{j+s}=D_{j+s-1}\eta^r$, by \eqref{eq:staircase_faces}, as
$r\leq j-1<j+s-1$. For the face $d_j\Delta[j-1]$, with $j>0$, we
have
\[
  d_s^hd_0^h\xi_{j-1}
  =d_0^hd_{s+1}^h\xi_{j-1}
  =d_0^h\eta^{j+s}_{j-1}
  =d_j^v\eta^{j+s}_j,
\]
by \eqref{eq:IH2-horizontal} for $\xi_{j-1}$ and $s+1$, and by
\eqref{eq:staircase_relation} for $\eta^{j+s}$. Hence the square
\eqref{eq:stage-j-square} commutes, and it admits a lift $\xi_j$, as
$R_j$ is a trivial fibration. The identity $R_j\xi_j=H_j$ gives
$P_{q-j}(\xi_j)=c_{\sigma_{\geq j}}^{(j)}$ and
\eqref{eq:IH2-horizontal}, and the identity
$\xi_j|_{\mathcal L_j}=\Phi_j$ gives \eqref{eq:IH3-vertical} and
\eqref{eq:IH3-staircase_gluing}. This completes the induction.

For the last assertion, let $\xi_0,\ldots,\xi_{q-1}$ be simplices as
in the statement. The verifications made above for
\eqref{eq:Phij-faces}, which for $j=q$ only involve the simplices
$\xi_{q-1}$ and $\xi_{q-2}$ and the relations in the statement,
show that the assignments $\Phi_q$ define a horn
$\Lambda^a[q]\to\mathbb E_0$ over the degenerate simplex at
$\sigma_{\geq q}$. Let $\xi_q:\Delta[q]\to\mathbb E_0$ be an
extension of $\Phi_q$. Since $(\bisimFin)_0$ is discrete, the
simplex $\xi_q$ lies over the degenerate simplex at
$\sigma_{\geq q}$. Moreover, we have
$d_q^v\xi_q=d_0^h\xi_{q-1}$, and so, by
\eqref{eq:IH3-staircase_gluing}, the tuple
$\xi=(\xi_0,\ldots,\xi_q)$ is a $q$-simplex of
$\Stair(\mathbb E)$ over $\sigma$. Finally, for $r\ne a$ and
$0\le i\le q-1$, we have
\begin{equation}\label{eq:total-faces-of-filler}
  (D_r\xi)_i=
  \begin{cases}
    d_{r-i}^h\xi_i=\eta^r_i, & i<r,\\[1mm]
    d_r^v\xi_{i+1}=\eta^r_i, & i\ge r,
  \end{cases}
\end{equation}
by \eqref{eq:staircase_faces}, where the first case follows from
\eqref{eq:IH2-horizontal}, and the second one from
\eqref{eq:IH3-vertical} if $i+1<q$ and from the definition of
$\Phi_q$ if $i+1=q$. Hence $D_r\xi=\eta^r$ for every
$r\ne a$, i.e., the simplex $\xi$ extends $\eta$.
\end{proof}

\bibliographystyle{amsalpha}
\bibliography{bibliography}
\end{document}